\documentclass[12pt]{article}
\usepackage{amssymb,amsmath,amsfonts}
\usepackage{graphicx}
\usepackage{relsize}
\usepackage{cite}
\usepackage{esvect}
\usepackage{txfonts}

\newtheorem{theorem}{Theorem}[section]
\newtheorem{definition}[theorem]{Definition}
\newtheorem{corollary}[theorem]{Corollary}
\newtheorem{proposition}[theorem]{Proposition}
\newtheorem{remark}[theorem]{Remark}
\newtheorem{lemma}[theorem]{Lemma}
\newtheorem{example}[theorem]{Example}

\newtheorem{claim}[theorem]{Claim}

\newtheorem{assumption}[theorem]{Assumption}
\newtheorem{fnassumption}[theorem]{Finiteness Assumption}
\newtheorem{wplemma}[theorem]{Whitney Partition Lemma}
\newtheorem{ptchlemma}[theorem]{Patching Lemma}
\newtheorem{mainlemma}[theorem]{Main Lemma}
\newtheorem{indhypothesis}[theorem]{Inductive Hypothesis}
\newtheorem{hmla}[theorem]{Hypotheses of the Main Lemma for the Label $\Ac$}
\newtheorem{lassumption}[theorem]{Large $A$ Assumption}
\newtheorem{lassum}[theorem]{Large $A$ Assumption for Whitney Partitions}
\newtheorem{smepsassumption}[theorem]{Small $\ve$ Assumption}

\newtheorem{statement}[theorem]{Statement}
\newcommand {\Ac}      {{\mathcal A}}
\newcommand {\Bc}      {{\mathcal B}}

\newcommand {\Ec}      {{\mathcal E}}
\newcommand {\Fc}      {{\mathcal F}}

\newcommand {\Kc}      {{\mathcal K}}

\newcommand {\Mc}      {{\mathcal M}}

\newcommand {\mA}      {{\mathbb A}}
\newcommand {\mB}      {{\mathbb B}}
\newcommand {\mC}      {{\mathbb C}}
\newcommand {\mD}      {{\mathbb D}}

\newcommand {\R}       {{\mathbb R}}
\newcommand {\N}       {{\mathbb N}}
\newcommand {\SP}       {{\mathbb S}}

\newcommand {\tB}      {\widetilde{B}}
\newcommand {\tC}      {\widetilde{C}}

\newcommand {\tE}      {\widetilde{E}}
\newcommand {\tF}      {\widetilde{F}}

\newcommand {\tK}      {\widetilde{K}}
\newcommand {\tL}      {\widetilde{L}}

\newcommand {\tS}      {\widetilde{S}}
\newcommand {\tT}      {\widetilde{T}}

\newcommand {\tW}      {\widetilde{W}}
\newcommand {\tX}      {\widetilde{X}}

\newcommand {\tf}      {\tilde{f}}

\newcommand {\ta}      {\tilde{a}}

\newcommand {\tv}      {\tilde{v }}

\newcommand {\tz}      {\tilde{z}}
\newcommand {\tl}      {\tilde{\ell}}

\newcommand {\trh}     {\tilde{\rho}}

\newcommand {\hF}      {\hat{F}}
\newcommand {\hf}      {\hat{f}}
\newcommand {\hS}      {\hat{S}}

\newcommand {\ha}      {\hat{a}}
\newcommand {\hb}      {\hat{b}}
\newcommand {\hc}      {\hat{c}}
\newcommand {\hx}      {\hat{x}}
\newcommand {\hxi}     {\hat{\xi}}
\newcommand {\ve}      {\varepsilon}

\newcommand {\gmh}     {\hat{\gamma}}

\newcommand {\RN}      {\R^n}
\newcommand {\RD}      {\R^D}
\newcommand {\MS}      {\Mc}
\newcommand {\BY}      {B_{\BS}}
\newcommand {\BL}      {B_{\BS}}
\newcommand {\dm}      {\rho}
\newcommand {\tdm}     {\widetilde{\dm}}
\newcommand {\dt}      {d}
\newcommand {\BS}      {Y}
\newcommand {\KM}      {\Kc_m(\BS)}
\newcommand {\KMO}     {\Kc_{m+1}(\BS)}
\newcommand {\KX}      {\Kc(\BS)}
\newcommand {\KMY}     {\Kc_m(\BS)}
\newcommand {\CNMY}    {\Conv_m(\BS)}
\newcommand {\CNY}     {\Conv^{(\Fc)}(\BS)}
\newcommand {\KH}      {\Kc(H)}
\newcommand {\SX}      {S_{\hspace{-0.5mm}Y}}
\newcommand {\SW}      {\tS}
\newcommand {\BR}      {\text{\sc Br}}
\newcommand {\AM}      {\AFF_m(\BS)}
\newcommand {\FN}      {N(m,\BS)}
\newcommand {\BU}      {\ell_\infty(U)}
\newcommand {\LTI}     {\ell^2_\infty}
\newcommand {\LTHI}    {\ell^3_\infty}
\newcommand {\REL}     {{\sc Rel}~}
\newcommand {\TRL}     {\text{\sc Rel}}
\newcommand {\ks}      {k^{\sharp}}
\newcommand {\kn}      {k^*}
\newcommand {\gnd}     {\gamma}
\newcommand {\gz}      {\gamma_0}
\newcommand {\GL}      {\Gamma_\ell}
\newcommand {\GLO}     {\Gamma_{\ell-1}}
\newcommand {\GH}      {\hat{\Gamma}}
\newcommand {\hAc}     {\hat{\Ac}}
\newcommand {\dl}      {\delta}
\newcommand {\Gm}      {\Gamma}
\newcommand {\hsg}     {\hat{\sigma}}
\newcommand {\hv}      {\hat{v}}
\newcommand {\hz}      {\hat{\zeta}}
\newcommand {\BXR}     {B(x_0,r_0)}
\newcommand {\GLA}     {\Gm_{\ell(\Ac)}}
\newcommand {\GLAO}    {\Gm_{\ell(\Ac)-1}}
\newcommand {\FF}      {F}

\newcommand {\ff}      {f}
\newcommand {\BO}      {\tB}
\newcommand {\tff}     {\tilde{\ff}}
\newcommand {\MPL}     {m+1}
\newcommand {\wtM}     {\widetilde{\Mc}}
\newcommand {\SB}      {S^{(Y)}}
\newcommand {\cng}      {c_{NC}}
\newcommand {\CNG}      {D_{NC}}
\newcommand {\CLS}      {C_{LS}}
\newcommand {\CWH}      {C_{Wh}}
\newcommand {\CETA}     {C_{\eta}}
\newcommand {\CAGR}     {C^\#}
\newcommand {\CLIP}     {C_{Lip}}
\newcommand {\DS}       {D^*}
\newcommand {\pl}       {10}
\newcommand {\emin}     {\ve^*}

\newcommand {\emp}     {\emptyset}
\newcommand {\ep}      {\varepsilon}
\newcommand {\vf}      {\varphi}
\newcommand {\je}      {\leftrightarrow}

\newcommand {\ip}[1]   {\langle{#1}\rangle}
\newcommand {\smed}    {\mathlarger{\sum}}

\newcommand {\cupsm}   {\mathsmaller{\bigcup}}
\newcommand {\capsm}   {\mathsmaller{\bigcap}}
\newcommand {\capbig}  {\mathlarger{\mathlarger{\cap}}}

\newcommand {\Lip}     {\operatorname{Lip}}
\newcommand {\diam}    {\operatorname{diam}}
\newcommand {\dist}    {\operatorname{dist}}

\newcommand {\Aff}     {\operatorname{Aff}}
\newcommand {\AFF}     {\operatorname{{\it Aff}}}
\newcommand {\Vect}    {\operatorname{Vect}}
\newcommand {\Conv}    {\operatorname{Conv}}
\newcommand {\cnv}     {\operatorname{conv}}
\newcommand {\dhf}     {\operatorname{d_H}}
\newcommand {\affspan} {\operatorname{{\it affhull\hspace*{0.5mm}}}}
\newcommand {\Orb}     {\operatorname{Orb}}
\newcommand {\Or}      {\operatorname{\text{\sc Or}}}
\newcommand {\RELX}    {\operatorname{RELX}}

\newcommand {\cl}      {\operatorname{\,cl}}
\newcommand {\cntr}    {\operatorname{center}}
\newcommand {\Sel}     {\operatorname{Selec}}
\newcommand {\DA}     {\operatorname{DA}}
\newcommand {\VST}     {\vskip 1mm}
\newcommand {\VSU}     {\vskip 1mm}
\newcommand {\bx}      {\hspace{10mm}$\Box$}
\newcommand {\rbx}     {\hspace{10mm}$\vartriangleleft$}

\newcommand {\BEND}    {\hspace{10mm}\Box}
\newcommand {\nn}      {\nonumber}
\newcommand {\rf}[1]    {(\ref{#1})}      
\newcommand {\reff}[1] {\ref{#1}}         
\newcommand{\lbl}[1]      {\label{#1}}       
\newcommand{\be}          {\begin{eqnarray}}
\newcommand{\bel}[1]      {\begin{eqnarray} \label{#1}}
\newcommand{\ee}           {\end{eqnarray}}
\newcommand {\SECT}[2] {\section*{\centerline{\normalsize
{\bf #1}}} \setcounter{section}{#2}
\setcounter{theorem}{0}\setcounter{equation}{0}}
\begin{document}
\parindent 1em
\parskip 0mm
\medskip
\centerline{\large{\bf\Large Sharp finiteness principles for Lipschitz selections: long version}}\vskip 9mm
\vspace*{7mm} \centerline{By~  {\sc Charles Fefferman}}\vspace*{3 mm}
\centerline {\it Department of Mathematics, Princeton University,}\vspace*{1 mm}
\centerline{\it Fine Hall Washington Road, Princeton, NJ 08544,  USA}\vspace*{1 mm}
\centerline{\it e-mail: cf@math.princeton.edu}
\vspace*{3mm}
\centerline{\sc and}\vspace*{4mm}
\centerline{{\sc Pavel Shvartsman}}\vspace*{3 mm}
\centerline {\it Department of Mathematics, Technion - Israel Institute of Technology,}\vspace*{1 mm}
\centerline{\it 32000 Haifa, Israel}\vspace*{1 mm}
\centerline{\it e-mail: pshv@technion.ac.il}
\vspace*{10 mm}
\renewcommand{\thefootnote}{ }
\footnotetext[1]{{\it\hspace{-6mm}Math Subject
Classification} 46E35\smallskip
\par This research was supported by Grant No 2014055 from the United States-Israel Binational Science Foundation (BSF). The first author was also supported in part by NSF grant DMS-1265524 and AFOSR grant FA9550-12-1-0425.}

\begin{abstract} Let $(\MS,\dm)$ be a metric space and let $\BS$ be a Banach space. Given a positive integer $m$, let $F$ be a set-valued mapping from $\MS$ into the family of all compact convex subsets of $\BS$ of dimension at most $m$. In this paper we prove a finiteness principle for the existence of a Lipschitz selection of $F$ with the sharp value of the finiteness number.

\end{abstract}
\renewcommand{\contentsname}{ }
\tableofcontents
\addtocontents{toc}{{\centerline{\sc{Contents}}}
\vspace*{10mm}\par}
\SECT{1. Introduction.}{1}

\addtocontents{toc}{~~~1. Introduction.\hfill \thepage\par\VST}
\indent\par {\bf 1.1. Main definitions and main results.}
\addtocontents{toc}{~~~~1.1. Main definitions and main results. \hfill \thepage\par}\medskip
\indent\par Let $(\MS,\dm)$ be a {\it pseudometric} space, i.e., $\dm: \MS\times\MS\to \R_+\cup\{+\infty\}$ is symmetric, non-negative, satisfies the triangle inequality, and $\rho(x,x)=0$ for all $x\in\Mc$, but $\rho(x,y)$ may be $0$ for $x\ne y$ and $\rho$ may admit the value $+\infty$.
We call a pseudometric space $(\MS,\dm)$ finite if $\Mc$ is finite, but we say that the pseudometric $\rho$ is finite if $\rho(x,y)$ is finite for every $x,y\in \Mc$.
\par Let $(\BS,\|\cdot\|)$ be a Banach space. Given a non-negative integer $m$ we let $\KM$ denote the family of all non-empty convex compact subsets of $\BS$ of dimension at most $m$. We recall that a (single-valued) mapping $f:\MS\to\BS$ is called a {\it selection} of a set-valued mapping $F:\MS\to\KM$ if $f(x)\in F(x)$ for all $x\in\MS$. A selection $f$ is said to be Lipschitz if there exists a constant $\lambda>0$ such that
\bel{LIP-N}
\|f(x)-f(y)\|\le\lambda\,\dm(x,y)~~~\text{for all} ~~~x,y\in\MS.
\ee
We let $\Lip(\MS,\BS)$ denote the space of all Lipschitz continuous mappings from $\MS$ into $\BS$ equipped with the seminorm $\|f\|_{\Lip(\MS,\BS)}=\inf\lambda$ where the infimum is taken over all constants $\lambda$ which satisfy \rf{LIP-N}.
\par Let
\bel{FNCN}
N(m,\BS)=2^{\min\{m+1,\dim Y\}}\,.
\ee
\par Our main result is the following ``Finiteness Principle for Lipschitz Selections''.
\begin{theorem}\lbl{MAIN-FP} Let $(\Mc,\rho)$ be a pseudometric space and let $F:\MS\to\KM$ be a set-valued mapping. Assume that for every subset $\Mc'\subset\MS$ consisting of at most $\FN$ points, the restriction $F|_{\Mc'}$ of $F$ to $\Mc'$ has a Lipschitz selection $f_{\Mc'}:\Mc'\to\BS$ whose seminorm satisfies $\|f_{\Mc'}\|_{\Lip(\Mc',\BS)}\le 1$.
\par Then $F$ has a Lipschitz selection $f:\MS\to\BS$ with $\|f\|_{\Lip(\MS,\BS)}$ bounded by a constant depending only on $m$.
\end{theorem}
\smallskip
\par In Section 6 we prove a variant of this result for {\it finite} pseudometric spaces $(\Mc,\dm)$. We show that in this case the family $\KM$ in the formulation of Theorem \reff{MAIN-FP} can be replaced by a wider family $\CNMY$ of {\it all non-empty convex (not necessarily compact) subsets of $Y$ of dimension at most $m$}. See Theorem \reff{FINITE}.
\par Before we discuss the main ideas of the proof of Theorem \reff{MAIN-FP} let us recall something of the history of the Lipschitz selection problem. The finiteness principle given in this theorem has been conjectured for $\BS=\RD$ in \cite{BS}, and, in full generality, in \cite{S02}.
\par Note that the sharp finiteness number for the case of the trivial distance function $\dm\equiv 0$ is equal to $n(m,Y)=\min\{m+2,\dim Y+1\}$. In fact, the finiteness principle for Lipschitz selections with respect to this trivial pseudometric coincides with the classical Helly's Theorem \cite{DGK}: there is a point common to all of the family of sets $\{F(x): x\in\MS\}\subset\KM$ provided for every $n(m,Y)$-point subset $\Mc'\subset\MS$ the family $\{F(x):x\in\Mc'\}$ has a common point. Thus Theorem \reff{MAIN-FP} can be considered as a certain generalization of Helly's Theorem to the case of arbitrary pseudometrics.
\par For the case $\BS=\R^2$ Theorem \reff{MAIN-FP} was proved in \cite{S02}. Fefferman, Israel and Luli \cite{FIL1} proved this theorem for $(\MS,\dm)=(\RN,\|\cdot\|)$ and $\BS=\RD$. An analog of Theorem \reff{MAIN-FP} for set-valued mappings into the family $\AM$ of all {\it affine subspaces} of $\BS$ of dimension at most $m$ has been proven by Shvartsman in \cite{S86} ($\BS=\RD$, see also \cite{S92}), \cite{S01} ($\BS$ is a Hilbert space), and \cite{S04} ($\BS$ is a Banach space).
\par The number $\FN$ from the formulation of Theorem \reff{MAIN-FP} is in general sharp.
\begin{theorem}\lbl{SH-FP}(\cite{S92}, \cite{S02}) Theorem \reff{MAIN-FP} is false in general if $\FN$ is replaced by some number $N$ with $N<\FN$.
\par Thus, for every non-negative integer $m$ and every Banach space $\BS$ of dimension $\dim Y\ge m$, there exist a metric space $(\MS,\dm)$ and a set-valued mapping
$F:\MS\to\KM$ such that the restriction $F|_{\Mc'}$ of $F$ to every subset $\Mc'\subset\MS$ consisting of at most
$\FN-1$ points has a Lipschitz selection $f_{\Mc'}$ with the seminorm $\|f_{\Mc'}\|_{\Lip(\Mc',\BS)}\le 1$,
but nevertheless $F$ does not have a Lipschitz selection.
\end{theorem}
\par See also Section 8 where we describe main ideas of the proof of this result for $m=1$ and $m=2$.\bigskip
\par {\bf 1.2. Main ideas of our approach.}
\addtocontents{toc}{~~~~1.2. Main ideas of our approach. \hfill \thepage\par\VSU}
\smallskip
\indent\par Let us briefly indicate the main ideas of the proof of Theorem \reff{MAIN-FP}.
\par One of the main ideas in this proof is to bring in the notion of Nagata dimension. The Nagata dimension (or Assouad-Nagata dimension) \cite{NGT,A} of a metric space is a certain metric version of the topological dimension. We recall one of the equivalent definitions of this notion. See, e.g., \cite{BDHM}.
\begin{definition}\lbl{NG-C} (``Nagata Condition'' and ``Nagata Dimension'') {\em We say that $(X,d)$ satisfies {\it the Nagata condition} if there exist a constant $\cng\in(0,1]$ and a non-negative integer $\CNG$ such that for every $s>0$ there exists a cover of $X$ by subsets of diameter at most $s$, at most $\CNG+1$ of which meet any given ball in $X$ of radius $\cng s$.
\smallskip
\par We refer to the smallest value of $\CNG$ as {\it the Nagata dimension}. More specifically, the Nagata dimension $\dim_N X$ of a metric space $(X,d)$ is the smallest integer $n$ for which there exists a constant $C\ge 1$ such that for all $s>0$, there exists a covering of $X$ by subsets of diameter at most $s$ with every ball in $X$ of diameter at most $s/C$ meeting at most $n+1$ elements of the covering.}
\end{definition}
\par We refer the reader to \cite{LS,BDHM,LDR}
and references therein for numerous results related to the Nagata condition and dimension.
\smallskip
\begin{theorem}\lbl{NG-FP} Let $(X,\dt)$ be a finite metric space satisfying the Nagata condition with constants $\cng$, $\CNG$. Given $m\in\N$ there exist a constant $\ks\in\N$ depending only on $m$, and a constant $\gnd>0$ depending only on $m$, $\cng$, $\CNG$, for which the following holds :
\smallskip
\par Let $\lambda$ be a positive constant and let $F:X\to\CNMY$ be a set-valued mapping such that, for every subset $X'\subset X$ consisting of at most $\ks$ points, the restriction $F|_{X'}$ of $F$ to $X'$ has a Lipschitz selection $f_{X'}:X'\to\BS$ whose seminorm satisfies $\|f_{X'}\|_{\Lip(X',\BS)}\le \lambda$.
\par Then $F$ has a Lipschitz selection $f:X\to\BS$ with $\|f\|_{\Lip(X,\BS)}\le \gnd\lambda$.
\end{theorem}
\par We consider the proof of Theorem \reff{NG-FP}, which we present in Sections 2-4, to be the most difficult technical part of this paper.
\smallskip
\par To establish Theorem \reff{NG-FP} we adapt the proof of a finiteness principle for $C^m$ selection \cite{FIL1} from $\RN$ to a metric space $X$ of bounded Nagata dimension. As in \cite{FIL1}, the geometry of certain convex sets $\Gamma_\ell(x)$, $(\ell\ge 0, x\in X)$ plays a crucial role. We refer the reader to the introduction of \cite{FIL1} and the website \cite{FIL-H}.
\smallskip
\par In Section 2.1 we consider an important family of metric spaces with finite Nagata dimension - the family of {\it finite metric trees}. We recall that a finite metric space $(X,\dt)$ is said to be a {\it metric tree} if $X$ is equipped with a structure of a (graph-theoretic) tree so that for every $x,y\in X$
$$
d(x,y)=\sum\limits_{i=0}^{k-1}d(z_i,z_{i+1})
$$
where $\{z_0,z_1,...,z_k\}$ is {\it the unique ``path'' in $X$ joining $x$ to $y$} (i.e., $z_0=x$, $z_k=y$, $z_i\neq z_j$ for $i\neq j$, and $z_j$ joined to $z_{j+1}$ by an edge).
\par It is proven in \cite{LS} that {\it every metric tree satisfies the Nagata condition} (with absolute constants $\cng,\CNG$), and the Nagata dimension of an arbitrary metric tree is $1$. See also Lemma \reff{MTR-ND}.
\par Hence we conclude that Theorem \reff{NG-FP} is true for every finite metric tree. See Corollary \reff{A-12}.\medskip
\par Let
$$
\Kc(\BS)=\cupsm\,\{\Kc_m(\BS): m\in\N\}
$$
be the family of all non-empty finite dimensional convex compact subsets of the Banach space $\BS$, and let $\dhf$ denote the {\it Hausdorff distance} between subsets of $\BS$. Basing on Corollary \reff{A-12}, in Section 5 we prove the following theorem which actually reduces the original problem to the case of the metric space $(\KM,\dhf)$.
\begin{theorem}\lbl{HDS-M} Let $(\Mc,\dm)$ be a metric space.  Let $m\ge 1$ and let $F:\MS\to\KM$ be a set-valued mapping. Suppose that for every subset $\Mc'\subset\MS$ consisting of at most $\ks$ points, the restriction $F|_{\Mc'}$ of $F$ to $\Mc'$ has a Lipschitz selection $f_{\Mc'}:\Mc'\to\BS$ whose seminorm satisfies $\|f_{\Mc'}\|_{\Lip(\Mc',\BS)}\le 1$.
\par Then there exists a mapping $G:\MS\to\KM$ satisfying the following conditions:\smallskip
\par (i). $G(x)\subset F(x)$ for every $x\in\MS$;
\medskip
\par (ii). For every $x,y\in \MS$ the following inequality
$$
\dhf(G(x),G(y))\le \gz\,\dm(x,y)
$$
holds. Here $\gz$ is a constant depending only on $m$, and  $\ks=\ks(m)$ is the constant from Theorem \reff{NG-FP}.
\end{theorem}
\par At the next step of the proof of Theorem \reff{MAIN-FP} we apply to the mapping $G$ the following Lipschitz selection theorem for the metric space $(\Kc(\BS),\dhf)$.
\begin{theorem}\lbl{ST-P}(\cite{S02}) There exists a mapping $S_{\BS}:\Kc(\BS)\to\BS$ such that \smallskip
\par (i). $S_{\BS}(K)\in K$ for each $K\in\Kc(\BS)$;
\smallskip
\par (ii). For every $K_1,K_2\in\Kc(\BS)$,
$$
\|S_{\BS}(K_1)-S_{\BS}(K_2)\|\le \gamma_1\dhf(K_1,K_2)
$$
where $\gamma_1=\gamma_1(\dim K_1,\dim K_2)$ is a constant depending only on dimensions of $K_1$ and $K_2$.
\end{theorem}
\par We refer to $\SX(K)$ as a {\it Steiner-type point} of a convex set $K\in\Kc(\BS)$. See Section 7 for more detail.
\par Finally, we put
$$
f(x)=S_{\BS}(G(x)), ~~~x\in\MS,
$$
where $G:\MS\to\Kc_m(\BS)$ is the set-valued mapping from Theorem \reff{HDS-M}.
\par Clearly, by part (i) of Theorem \reff{ST-P} and part (i) of Theorem \reff{HDS-M},
$$
f(x)=S_{\BS}(G(x))\in G(x)\subset F(x)~~~\text{for all}~~~x\in\Mc,
$$
proving that the function $f:\Mc\to \BS$ is a {\it selection} of $F$. In turn, by part (ii) of Theorem \reff{ST-P} and part (ii) of Theorem \reff{HDS-M}, for every $x,y\in\MS$
$$
\|f(x)-f(y)\|=\|S_{\BS}(G(x))-S_{\BS}(G(y))\|\le \gamma_1\,\dhf(G(x),G(y))\le
\gamma_0\,\gamma_1\,\dm(x,y).
$$
Here $\gamma_1$ is a constant depending only on $\dim G(x)$ and $\dim G(y)$. Since $\dim G(x),\dim G(y)\le m$, and  $\gamma_0$ depends only on $m$, the Lipschitz seminorm of $f$ on $\Mc$ is bounded by a constant depending only on $m$.
\par This proves a version of Theorem \reff{MAIN-FP} with $\FN$ replaced by $\ks$ for an arbitrary {\it metric} space $(\MS,\dm)$. See Corollary \reff{LS-MSP}. \smallskip
\par Using this result, in Section 6 we prove a similar version of Theorem \reff{MAIN-FP} for the general case, i.e., for an arbitrary {\it pseudometric} space $(\MS,\dm)$. See Proposition \reff{PM-SP}.
\par Finally, using Theorem \reff{D-2M} below, we obtain the statement of Theorem \reff{MAIN-FP} in its original form.
\begin{theorem}\lbl{D-2M} Let $(\Mc,\rho)$ be a finite pseudometric space with a finite pseudometric $\rho$, and let ${F:\Mc\to\CNMY}$ be a set-valued mapping. Suppose that for every subset $\Mc'\subset\Mc$ with $\#\Mc'\le \FN$ the restriction $F|_{\Mc'}$ has a Lipschitz selection $f_{\Mc'}:\Mc'\to \BS$ with $\|f_{\Mc'}\|_{\Lip(\Mc',\BS)}\le 1$.
\par Then $F$ has a Lipschitz selection $f:\Mc\to \BS$ with $\|f\|_{\Lip(\Mc,\BS)}\le \gamma$ where $\gamma$ is a positive constant depending only on $m$ and $\#\Mc$.
\end{theorem}
(Recall that $\CNMY$ denotes the family of all non-empty convex subsets of $Y$ of dimension at most $m$.)
\par We prove Theorem \reff{D-2M} in Section 6. This result  enables us to replace the finiteness number $\ks$ by the required sharp finiteness number $\FN$, completing the proof of Theorem \reff{MAIN-FP}.
\medskip
\par In Section 8 we present various remarks and comments related to the sharp finiteness principle proven in Theorem \reff{MAIN-FP}.
\bigskip
\par The existence of Lipschitz selections is closely related to {\it Whitney's Extension Problem} \cite{W1}:
\smallskip
\par Fix $m,n\ge 1$. Given $E\subset\RN$ and $\vf:E\to\R$, decide whether $\vf$ extends to a $C^m$ function $f:\RN\to\R$. If such an extension exists, then how small can we take its $C^m$-norm?
\medskip
\par There is a finiteness theorem for such problems and their relatives; see Brudnyi-Shvartsman
\cite{BS,BS1,BS3,S1,S02,S-Tr} and later papers by Fefferman, Klartag, Israel, Luli
\cite{F2,F-J,F4,F6,F7,F8,F-Bl,FIL1,FIL2}. See also A. Brudnyi, Yu. Brudnyi \cite{BB}.
\par In Brudnyi-Shvartsman \cite{BS1,BS3,S1,S02,S-Tr}, Lipschitz selection served as the main tool to attack Whitney's Problem. The later work \cite{F2,F-J,F4,F6,F7,F8,FIL1,FIL2} made no explicit mention of Lipschitz selection, but broadened Whitney's Problem to study $C^m$ functions $f:\RN\to\R$ that agree only approximately with a given function $\vf$ on $E$.
\par A Lipschitz selection problem can obviously be viewed as a search for a Lipschitz mapping $f:\Mc\to Y$ that agrees approximately with data.\smallskip
\par As in \cite{FIL1}, our present results lead to questions about efficient computation for Lipschitz selection problems on finite metric spaces. In connection with such issues, we ask whether the results of Har-Peled and Mendel \cite{HP-M} on the Well Separated Pairs Decomposition \cite{C-Kos} can be extended from doubling metrics to metrics of bounded Nagata dimension.
\medskip

\par {\bf Acknowledgments.} We are grateful to Alexander  Brudnyi, Arie Israel, Bo'az Klartag, Garving (Kevin) Luli and the participants of the 10th Whitney Problems Conference, Williamsburg, VA, for valuable conversations.
\par We are grateful also to the College of William and Mary, Williamsburg, VA, USA, the American Institute of Mathematics, San Jose, CA, USA, the Fields Institute, Toronto, Canada, the Banff International Research Station, Banff, Canada, the Centre International de Rencontres Math\'ematiques (CIRM), Luminy, Marseille, France, and the Technion, Haifa, Israel, for hosting and supporting workshops on the topic of this paper and closely related problems.
\par Finally, we thank the US-Israel Binational Science Foundation, the US National Science Foundation, the Office
of Naval Research and the Air Force Office of Scientific Research for generous support.
\bigskip
\SECT{2. Nagata dimension and Whitney partitions on metric trees.}{2}
\addtocontents{toc}{2. Nagata condition and Whitney partitions on metric spaces. \hfill\thepage\par\VST}

\indent\par {\bf 2.1. Metric trees and Nagata condition.}
\addtocontents{toc}{~~~~2.1. Metric trees and Nagata condition. \hfill \thepage\par}\medskip
\par Let $(X,d)$ be a metric space. We write $B(x,r)$ to denote the ball $\{y\in X:d(x,y)<r\}$
(strict inequality) in the metric space $(X,d)$. We also write $\diam A=\sup\,\{d(a,b):a,b\in A\}$ and
$$
\dist(A',A'')=\inf\{d(a',a''): a'\in A', a''\in A''\}
$$
to denote the diameter of a set $A\subset X$ and the distance between sets $A',A''\subset X$ respectively.
\smallskip
\par Let us consider an important example of a metric space with finite Nagata dimension.
\par Let $T=(X,E)$ be a {\it finite tree}. Here $X$ denotes the set of nodes and $E$ denotes the set of edges of $T$. We write $x\je y$ to indicate that nodes $x,y\in X$, $x\ne y$, are joined by an edge; we denote that edge by $[xy]$.
\par Suppose we assign a positive number $\Delta(e)$ to each edge $e\in E$. Then we obtain a notion of distance $d(x,y)$ for any $x,y\in X$, as follows.
\par We set
\bel{D-XX}
d(x,x)=0~~~\text{for every}~~~x\in X.
\ee
\par Because $T$ is a tree, any two distinct nodes $x,y\in X$ are joined by one and only ``path''
$$
x=x_0\je x_1\je...\je x_L=y~~~\text{with all the}~~x_i ~~\text{distinct}.
$$
\par We define
\bel{D-TR}
d(x,y)=\smed_{i=1}^L\, \Delta([x_{i-1}x_i]).
\ee
We call the resulting metric space $(X,d)$ a {\it metric tree}.
\medskip
\begin{lemma}\lbl{MTR-ND} Every metric tree satisfies the Nagata condition with $\CNG=1$ and $\cng=1/16$.
\end{lemma}
\par {\it Proof.} Given a metric tree $(X,d)$, we fix an origin $0\in X$ and make the following definition:
\bigskip
\par Every point $x\in X$ is joined to the origin by one and only one ``path''
$$
0=x_0\je x_1\je...\je x_L=x,
~~~\text{with all the}~~x_i~~\text{\it distinct.}
$$
\par We call $x_0,x_1,..., x_L$ the {\it ancestors} of $x$. We define the {\it distinguished ancestor} of $x$, denoted  $\DA(x)$, to be $x_i$ for the smallest $i\in \{0,...,L\}$ for which
\bel{ANC}
d(0,x_i)>\lfloor d(0,x)\rfloor-1,
\ee
where $\lfloor \cdot\rfloor$ denotes the greatest integer function. (Note that there is at least one $x_i$ satisfying \rf{ANC}, namely $x_L=x$. Thus, every $x\in X$ has a distinguished ancestor.)
\medskip
\par We note two simple properties of $\DA(x)$, namely,
\medskip
\par (1) $d(x,\DA(x))\le 2$;
\medskip
\par (2) $\DA(x)$ is an ancestor of any ancestor $y$ of $x$ that satisfies  $d(0,y)>\lfloor d(0,x)\rfloor-1$.
\medskip
\par We now exhibit a Nagata covering of $X$ for the lengthscale $s=4$.\smallskip
\par For $q=0,1$ and $z\in X$, let
$$
X_q(z)=\{x\in X:
z=\DA(x)~~\text{and}~~\lfloor d(0,x)\rfloor=q\mod 2\}.
$$
Clearly, the $X_q(z)$ cover $X$. Moreover, (1) tells us that each $X_q(z)$ has diameter at most $4$.
\par We assert the following\smallskip
\par {\sc Claim:} If $z\ne z'$ and $q=q'$, then the distance from $X_q(z)$ to $X_{q'}(z')$ is at least $1/2$.
\medskip
\par The {\sc Claim} immediately implies that any given ball $B\subset X$ of radius $1/4$ meets at most one of the $X_0(z)$ and at most one of the $X_1(z)$, hence at most two of the $X_q(z)$.\smallskip
\par Let us establish the {\sc Claim}; if it were false, then we could find
$$
z\ne z',~q\in\{0,1\},~x\in X_q(z),~x'\in X_{q}(z')~~~
\text{with}~~~d(x,x')\le 1/2.
$$
\par We will derive a contradiction from these conditions as follows.
\par Because $d(x,x')\le 1/2$, we have
$$
|\,\lfloor d(0,x)\rfloor-\lfloor d(0,x')\rfloor\,|\le 1.
$$
On the other hand, $\lfloor d(0,x)\rfloor\equiv \lfloor d(0,x')\rfloor \mod 2$. Hence, $\lfloor d(0,x)\rfloor= \lfloor d(0,x')\rfloor$.
\par Next, let $\tz$ be the closest common ancestor of $x,x'$. Because $d(x,x')\le 1/2$, we have $d(x,\tz)\le 1/2$ and $d(x',\tz)\le 1/2$, and therefore the ancestor $\tz$ of $x$ satisfies
$$
d(0,\tz)>\lfloor d(0,x)\rfloor-1.
$$
Hence, (2) implies that $z$ is an ancestor of $\tz$. Similarly, $z'$ is an ancestor of $\tz$.
\par It follows that either $z$ is an ancestor of $z'$, or $z'$ is an ancestor of $z$. Without loss of generality, we may suppose that $z$ is an ancestor of $z'$. Consequently, $z$ is an ancestor of $x'$; moreover,
$$
d(0,z)>\lfloor d(0,x)\rfloor-1=\lfloor d(0,x')\rfloor-1.
$$
\par Thanks to (2), we now know that $z'$ is an ancestor of $z$. Thus, each of the points $z,z'$ is an ancestor of the other, and therefore $z=z'$, contradicting an assumption that the {\sc Claim} is false.
\smallskip
\par We have produced a covering of an arbitrary metric tree by subsets $X_i$ of diameter at most $4$, such that no ball of radius $1/4$ intersects more than two of the $X_i$.
\par Applying the above result to the metric tree $(X,\frac{4}{s}\,d)$ for given $s>0$, we produce a covering of $X$ by $X_i$ such that, with respect to $d$, each $X_i$ has diameter at most $s$, and no ball of radius $s/16$ meets more than two of the $X_i$. Thus, we have verified the Nagata condition for metric trees.\bx

\indent\par {\bf 2.2. Whitney Partitions.}
\addtocontents{toc}{~~~~2.2. Whitney Partitions.\hfill \thepage\par\VSU}\medskip
\par In this section, we prove the following result.
\begin{wplemma}\lbl{WPL} Let $(X,d)$ be a metric space,
and let $r(x)>0$ be a positive function on $X$. We assume the following, for constants $\cng\in(0,1]$, $\CNG\in\N\cup\{0\}$ and $\CLS\ge 1$:\vskip 1mm
\par \textbullet~({\sc Nagata Condition}) Given $s>0$ there exists a covering of $X$ by subsets $X_i$ $(i\in I)$ of diameter at most $s$, such that every ball of radius
$\cng s$ in $X$ meets at most $\CNG+1$ of the $X_i$.
\smallskip
\par \textbullet~({\sc Consistency of the Lengthscale}) Let $x,y\in X$. If $d(x,y)\le r(x)+r(y)$, then
\bel{C-LSC}
\CLS^{-1}\,r(x)\le r(y)\le \CLS r(x).
\ee
\par Let $a>0$.
\medskip
\par Then there exist functions $\vf_\nu:X\to\R$, and points $x_\nu\in X$, with the following properties: \smallskip
\par \textbullet~Each $\vf_\nu\ge 0$, and each
$\vf_\nu=0$ outside $B(x_\nu,ar_\nu)$. Here $r_\nu=r(x_\nu)$.
\par \textbullet~ Any given $x\in X$ satisfies $\vf_\nu(x)\ne 0$ for at most $C$ distinct $\nu$.
\par \textbullet~ $\sum\limits_\nu\,\vf_\nu=1$ on $X$.
\par \textbullet~For each $\nu$ and for all $x,y\in X$, we have
$$
|\vf_\nu(x)-\vf_\nu(y)|\le C\,d(x,y)/r_\nu.
$$
\par Here $C$ is a constant depending only on $\cng$, $\CNG$, $\CLS$ and $a$.
\end{wplemma}
\par {\it Proof.} We write $c,C,C'$, etc. to denote constants determined by $\cng$, $\CNG$, $\CLS$ and $a$. These symbols may denote different constants in different occurrences.
\par We introduce a large constant $A$ to be fixed later. We make the following
\begin{lassum}\lbl{WP-LA} {\em $A$ exceeds a large enough constant determined by $\cng$, $\CNG$, $\CLS$, $a$.}
\end{lassum}
\par We write $c(A), C(A),C'(A)$, etc. to denote constants determined by $A$, $\cng$, $\CNG$, $\CLS$, $a$. These symbols may denote different constants in different occurrences.
\par Let $P$ denote the set of all integer powers of $2$. For $s\in P$ let $(X(i,s))_{i\in I(s)}$ be a covering of $X$ given by the Nagata condition. Thus,
\medskip
$$
\diam X(i,s)\le s;
$$
and, for fixed $s\in P$,
\bel{2-C}
\text{any given}~~~x\in X~~~\text{lies in at most}~~C~~\text{of the sets}~~X^{++}(i,s).
\ee
Here
$$
X^{++}(i,s)=\{y\in X: d(y,X(i,s))<\cng s/64\}~~~~(i\in I(s)).
$$
\par We also define
$$
X^{+}(i,s)=\{y\in X: d(y,X(i,s))<\cng s/128\} ~~~\text{for}~~~(i\in I(s)).
$$
\par Let
$$
\theta_{i,s}(x)=\max\{0,(1-256\,d(x,X(i,s))/(\cng s))\}
$$
for $x\in X$, $i\in I(s)$, $s\in P$.
\par Then
\bel{3a-C}
0\le\theta_{i,s}\le 1,
\ee
\bel{3b-C}
\|\theta_{i,s}\|_{\Lip(X,\R)}\le C\,s^{-1},
\ee
and
$$
\theta_{i,s}=0~~~\text{outside}~~~X^+(i,s),
$$
but
$$
\theta_{i,s}=1~~~\text{on}~~~X(i,s).
$$
\par For each $s\in P$ and $i\in I(s)$, we pick a representative point $x(i,s)\in X(i,s)$. (We may assume that the $X(i,s)$ are all nonempty.) We let \REL denote
the set of all $(i,s)$ such that
\bel{4-C}
A^{-3}r(x(i,s))\le s\le A^{-1}r(x(i,s)).
\ee
\par We establish the basic properties of the set \REL.
\begin{lemma}\lbl{WP1} Given $x_0\in X$ there exists $(i,s)\in\TRL$ such that $x_0\in X(i,s)$ and therefore $\theta_{i,s}(x_0)=1$.
\end{lemma}
\par {\it Proof.} Pick $s_0\in P$ such that $s_0/2\le r(x_0)/A^2\le 2s_0$. Because the $X(i,s_0)$ $(i\in I(s_0))$ cover $X$, we may fix $i_0\in I(s_0)$ such that $x_0\in X(i_0,s_0)$. The points $x_0$ and $x(i_0,s_0)$ both belong to $X(i_0,s_0)$, hence
$$
d(x_0,x(i_0,s_0))\le \diam X(i_0,s_0)\le s_0\le 2r(x_0)/A^2\,.
$$
\par The Large $A$ Assumption \reff{WP-LA} and the {\sc Consistency of the Lengthscale} together now imply that
$$
cr(x_0)\le r(x(i_0,s_0))\le Cr(x_0),
$$
and therefore
$$
cs_0\le r(x(i_0,s_0))/A^2\le Cs_0\,.
$$
Thanks to the Large $A$ Assumption \reff{WP-LA}, we therefore have \rf{4-C} for $(i_0,s_0)$. Thus, $(i_0,s_0)\in\TRL$ and $x_0\in X(i_0,s_0)$.\bx
\begin{lemma}\lbl{WP2} If $(i,s)\in\TRL$ and $x_0\in X^{++}(i,s)$, then
$$
cA^{-3}r(x_0)\le s\le CA^{-1}r(x_0),
$$
and therefore
$$
\|\theta_{i,s}\|_{\Lip(X,\R)}\le CA^{3}/r(x_0).
$$
\end{lemma}
\par {\it Proof.} Both $x_0$ and $x(i,s)$ lie in $X^{++}(i,s)$, hence
$$
d(x_0,x(i,s))\le \diam X^{++}(i,s)\le 2\cng s/64+
\diam X(i,s)\le Cs\le Cr(x(i,s))/A
$$
thanks to \rf{4-C}.
\par The Large $A$ Assumption \reff{WP-LA} and {\sc Consistency of the Lengthscale} now tell us that
$$
cr(x_0)\le r(x(i,s))\le Cr(x_0),
$$
and therefore \rf{4-C} and \rf{3b-C} imply the conclusion of Lemma \reff{WP2}.\bx
\begin{corollary}\lbl{WP-C1} Any given point $x_0\in X$ lies in $X^{++}(i,s)$ for at most  $C(A)$ distinct $(i,s)\in \TRL$. Consequently, $\theta_{i,s}(x_0)$ is nonzero for at most $C(A)$ distinct $(i,s)\in\TRL$.
\end{corollary}
\par {\it Proof.} There are at most $C(A)$ distinct $s\in P$ satisfying the conclusion of Lemma \reff{WP2}.
For each such $s$ there are at most $C$ distinct $i$ such that $x_0\in X^{++}(i,s)$; see \rf{2-C}.\bx
\begin{corollary}\lbl{WP-C2} Suppose $X^{++}(i,s)\cap X^{++}(i_0,s_0)\ne\emp$ with $(i,s),(i_0,s_0)\in\TRL$. Then 
$$
c(A)s_0\le s\le C(A)s_0.
$$
\end{corollary}
\par {\it Proof.} Pick $x_0\in X^{++}(i,s)\cap X^{++}(i_0,s_0)$. Lemma \reff{WP2} gives
$$
c(A)r(x_0)\le s\le C(A)r(x_0)~~~\text{and}~~~
c(A)r(x_0)\le s_0\le C(A)r(x_0).\BEND
$$
\begin{lemma}\lbl{WP3} Let $(i_0,s_0), (i,s)\in\TRL$. If $x\in X^+(i_0,s_0)$, then for any $y\in X$
\bel{5-C}
|\theta_{i,s}(x)-\theta_{i,s}(y)|\le C(A)\,d(x,y)/s_0.
\ee
\end{lemma}
\par {\it Proof.} We proceed by cases.\medskip
\par {\it Case 1:} $d(x,y)<\cng s_0/128$. \medskip
\par Then $x,y\in X^{++}(i_0,s_0)$. If $x$ or $y$ belongs to $X^{++}(i,s)$, then Corollary \reff{WP-C2} tells us that 
$$
c(A)s_0\le s\le C(A)s_0\,;
$$
hence, \rf{3b-C} yields the desired estimate \rf{5-C}.
\par If instead neither $x$ nor $y$ belongs to $X^{++}(i,s)$, then $\theta_{i,s}(x)=\theta_{i,s}(y)=0$, hence \rf{5-C} holds trivially.\medskip
\par {\it Case 2:} $d(x,y)\ge \cng s_0/128$. Then \rf{3a-C} gives
$$
|\theta_{i,s}(x)-\theta_{i,s}(y)|\le 1\le C\,d(x,y)/s_0.
$$
Thus, \rf{5-C} holds in all cases.\bx
\medskip
\par Now define
\bel{6-C}
\Theta(x)=\smed_{(i,s)\in\TRL}\,\theta_{i,s}(x)~~~~
\text{for all}~~~x\in X.
\ee
\par Corollary \reff{WP-C1} shows that there are at most $C(A)$ nonzero summands in \rf{6-C} for any fixed $x$. Moreover, each summand is between $0$ and $1$ (see \rf{3a-C}), and for each fixed $x$, one of the summands is equal to $1$ (see Lemma \reff{WP1}). Therefore,
\bel{7-C}
1\le \Theta(x)\le C(A)~~~~\text{for all}~~~~x\in X.
\ee
\begin{lemma}\lbl{WP4} Let $x,y\in X$ and $(i_0,s_0)\in\TRL$. If $x\in X^+(i_0,s_0)$, then
$$
|\Theta(x)-\Theta(y)|\le C(A)\,d(x,y)/s_0\,.
$$
\end{lemma}
\par {\it Proof.} There are at most $C(A)$ distinct $(i,s)\in\TRL$ for which $\theta_{i,s}(x)$ or $\theta_{i,s}(y)$ is nonzero. For each such $(i,s)$ we apply Lemma \reff{WP3}, then sum over $(i,s)$.\bx
\medskip
\par Now, for $(i_0,s_0)\in\TRL$, we set
\bel{8-C}
\vf_{i_0,s_0}(x)=\theta_{i_0,s_0}(x)/\Theta(x)\,.
\ee
\par This function is defined on all of $X$, and it is zero outside $X^+(i_0,s_0)$. Moreover,
\bel{9-C}
\vf_{i_0,s_0}\ge 0~~~~~\text{and}~~~~~
\smed_{(i_0,s_0)\in\TRL}\vf_{i_0,s_0}=1~~~\text{on}~~~X.
\ee
\par Note that because
$$
\diam X^+(i_0,s_0)\le Cs_0\le C\,A^{-1}r(x(i_0,s_0))
$$
(see \rf{4-C}), the function $\vf_{i_0,s_0}$ is zero outside the ball $B(x(i_0,s_0),C\,A^{-1}r(x(i_0,s_0)))$.
Thanks to our Large $A$ Assumption \reff{WP-LA}, it follows that
\bel{10-C}
\vf_{i,s}~~~\text{is identically zero outside the ball} ~~~B(x(i,s),ar(x(i,s)))\,.
\ee
\begin{lemma}\lbl{WP5} For $x,y\in X$ and $(i_0,s_0)\in\TRL$, we have
$$
|\vf_{i_0,s_0}(x)-\vf_{i_0,s_0}(y)|\le C(A)\,d(x,y)/s_0.
$$
\end{lemma}
\par {\it Proof.} Suppose first that $x\in X^+(i_0,s_0)$. Then
$$
|\vf_{i_0,s_0}(x)-\vf_{i_0,s_0}(y)|=
\left|\frac{\theta_{i_0,s_0}(x)}{\Theta(x)}-
\frac{\theta_{i_0,s_0}(y)}{\Theta(y)}\right|\nn\\
\le
\frac{|\theta_{i_0,s_0}(x)-\theta_{i_0,s_0}(y)|}{\Theta(x)}+
\theta_{i_0,s_0}(y)\,
\frac{|\Theta(x)-\Theta(y)|}{\Theta(x)\Theta(y)}\,.
$$
The first term on the right is at most $C(A)\,d(x,y)/s_0$ by \rf{3b-C} and \rf{7-C}; the second term on the right is at most $C(A)\,d(x,y)/s_0$ thanks to \rf{3a-C}, Lemma \reff{WP4} and \rf{7-C}.
Thus,
\bel{11-C}
|\vf_{i_0,s_0}(x)-\vf_{i_0,s_0}(y)|\le C(A)\,d(x,y)/s_0
~~~~\text{if}~~~x\in X^+(i_0,s_0).
\ee
Similarly, \rf{11-C} holds if $y\in X^+(i_0,s_0)$.
\par Finally, if neither $x$ nor $y$ belongs to $X^+(i_0,s_0)$, then $\vf_{i_0,s_0}(x)=\vf_{i_0,s_0}(y)=0$, so \rf{11-C} is obvious.
\par Thus, \rf{11-C} holds in all cases.\bx
\begin{corollary}\lbl{WP-C3} For $x,y\in X$ and $(i_0,s_0)\in\TRL$, we have
$$
|\vf_{i_0,s_0}(x)-\vf_{i_0,s_0}(y)|\le C(A)\,d(x,y)/r(x(i_0,s_0)).
$$
\end{corollary}
\par {\it Proof.} Immediate from Lemma \reff{WP5} and inequalities \rf{4-C}.\bx
\bigskip
\par We can now finish the proof of the Whitney Partition Lemma \reff{WPL}. We pick $A$ to be a constant determined by $\cng$, $\CNG$, $\CLS$, $a$, taken large enough to satisfy the Large $A$ Assumption \reff{WP-LA}. We then take our functions $\vf_\nu$ to be the $\vf_{(i,s)}$ $((i,s)\in\TRL)$, and we take our $x_\nu$ to be the points $x(i,s)$  $((i,s)\in\TRL)$. We set $r_\nu=r({x_\nu})$.
\smallskip
\par The following holds: \smallskip
\par \textbullet~Each $\vf_\nu\ge 0$, and each
$\vf_\nu=0$ outside $B(x_\nu,ar_\nu)$; see \rf{9-C} and \rf{10-C}.
\smallskip
\par \textbullet~ Any given $x\in X$ satisfies $\vf_\nu(x)\ne 0$ for at most $C$ distinct $\nu$.
This follows from Corollary \reff{WP-C1}, equation \rf{8-C}, and the fact that $A$ is now determined by $\cng$, $\CNG$, $\CLS$, $a$.\smallskip
\par \textbullet~ $\sum\limits_\nu\,\vf_\nu=1$ on $X$; see \rf{9-C}.
\smallskip
\par \textbullet~For each $\nu$ and for all $x,y\in X$, we have
$$
|\vf_\nu(x)-\vf_\nu(y)|\le C\,d(x,y)/r_\nu;
$$
see Corollary \reff{WP-C3}, and note that $A$ is now determined by $\cng$, $\CNG$, $\CLS$ and $a$.
\smallskip
\par The proof of the Whitney Partition Lemma \reff{WPL} is complete.\bx
\smallskip
\begin{remark} {\em Later on there will be another Large $A$ Assumption, and another definition of the set {\sc Rel}, different from those in this section.\rbx}
\end{remark}
\bigskip

\indent\par {\bf 2.3. Patching Lemma.}
\addtocontents{toc}{~~~~2.3. Patching Lemma.\hfill \thepage\par\VSU}\medskip
\begin{ptchlemma}\lbl{PTHM} Let $(X,d)$ be a metric space, and let $Y$ be a Banach space. For each $\nu$ in some index set, assume we are given the following objects:\smallskip
\par \textbullet~ A point $x_\nu\in X$ and a positive number $r_\nu>0$ (a ``lengthscale'').
\smallskip
\par \textbullet~ A function $\theta_\nu:X\to \R$\,.
\smallskip
\par \textbullet~ A vector $\eta_\nu\in Y$ and a vector-valued function $F_\nu:X\to Y$.
\medskip
\par We make the following assumptions: We are given positive constants $\CLS\ge 1$, $\CWH$, $\CETA$, $\CAGR$, $\CLIP$, $\DS$, such that the following conditions are satisfied
\smallskip
\par \textbullet~ ({\sc Consistency of the Lengthscale})
\bel{C-LS}
\CLS^{-1}\le r_\nu/r_\mu\le \CLS~~~
\text{whenever}~~~d(x_\mu,x_\nu)\le r_\mu+r_\nu.
\ee
\par ({\sc Whitney Partition Assumptions})
\smallskip
\par \textbullet~ $\theta_\nu\ge 0$ on $X$ and
$\theta_\nu= 0$ outside $B(x_\nu,a\,r_\nu)$, where
\bel{A-SMALL}
a=(4\,\CLS)^{-1}.
\ee
\par \textbullet~
$|\theta_\nu(x)-\theta_\nu(y)|\le \CWH\cdot d(x,y)/r_\nu$~
for~ $x,y\in X$.
\smallskip
\par \textbullet~ Any given $x\in X$ satisfies $\theta_\nu(x)\ne 0$~ for at most $\DS$ distinct $\nu$.
\smallskip
\par \textbullet~ $\smed\limits_\nu\,\theta_\nu=1$ on $X$.
\smallskip
\par \textbullet~ ({\sc Consistency of the $\eta_\nu$})~
$\|\eta_\mu-\eta_\nu\|\le \CETA\cdot [r_\nu+r_\nu+d(x_\mu,x_\nu)]$~ for each~ $\mu,\nu$.
\smallskip
\par \textbullet~ ({\sc Agreement of $F_\nu$ with $\eta_\nu$})~
$\|F_\nu(x)-\eta_\nu\|\le \CAGR\, r_\nu$~ for~
$x\in B(x_\nu,r_\nu)$.
\smallskip
\par \textbullet~ ({\sc Lipschitz continuity of $F_\nu$})~
$\|F_\nu(x)-F_\nu(y)\|\le \CLIP\cdot d(x,y)$~ for~
$x,y\in B(x_\nu,r_\nu)$.
\bigskip
\par Define
$$
F(x)=\sum_\nu\theta_\nu(x)\,F_\nu(x)~~~\text{for}~~x\in X.
$$
Then $F$ satisfies
$$
\|F(x)-F(y)\|\le C\,d(x,y)~~~\text{for}~~~x,y\in X,
$$
where $C$ is determined by $\CLS$, $\CWH$, $\CETA$, $\CAGR$, $\CLIP$, $\DS$.
\end{ptchlemma}
\medskip
\par To start the proof of the {\sc Patching Lemma} \reff{PTHM}, we define
$$
\TRL(x)=\{\nu: \theta_\nu(x)\ne 0\}, ~~~x\in X.
$$
Then $1\le\#(\TRL(x))\le D^*$, and
\bel{D-NU}
d(x,x_\nu)\le a\,r_\nu~~~\text{for}~~~v\in\TRL(x).
\ee
\par We also recall that $\CLS\ge 1$ and $a=(4\,\CLS)^{-1}$ so that
\bel{CLS-A}
\CLS\cdot a= 1/4~~~\text{and}~~~a\le 1/4.
\ee
\vskip 2mm
\par We will use the following result.
\begin{lemma}\lbl{STR} Let $\nu,\nu_0\in\TRL(x)$, $\mu_0\in\TRL(y)$, and suppose that $d(x,y)\le a\cdot[r_{\nu_0}+r_{\mu_0}]$. Then
$$
x,y\in B(x_\nu,r_\nu)\cap B(x_{\nu_0},r_{\nu_0})\cap B(x_{\mu_0},r_{\mu_0})
$$
and the ratios
$$
r_{\nu_0}/r_{\mu_0},~r_{\mu_0}/r_{\nu_0},~
r_{\nu}/r_{\nu_0},~ r_{\nu_0}/r_{\nu},~
r_{\nu}/r_{\mu_0},~ r_{\mu_0}/r_{\nu}
$$
are at most $\CLS$.
\end{lemma}
\par {\it Proof.} We have the following inequalities
\smallskip
\par $(\bigstar 1)$~
$d(x_\nu,x_{\nu_0})\le  d(x_\nu,x)+d(x,x_{\nu_0})\le a\,r_\nu+a\,r_{\nu_0}$,
\smallskip
\par $(\bigstar 2)$~
$d(x_{\nu_0},x_{\mu_0})\le  d(x_{\nu_0},x)+d(x,y)+d(y,x_{\mu_0})\le a\,r_{\nu_0}+[a\,r_{\nu_0}+a\,r_{\mu_0}]+a\,r_{\mu_0}$,
\smallskip
\par $(\bigstar 3)$~
$d(x_{\nu},x_{\mu_0})\le  d(x_{\nu},x)+d(x,y)+d(y,x_{\mu_0})\le a\,r_{\nu}+[a\,r_{\nu_0}+a\,r_{\mu_0}]+a\,r_{\mu_0}$.
\medskip
\par From $(\bigstar 1)$, $(\bigstar 2)$, \rf{CLS-A}, and {\sc Consistency of the Lengthscale} \rf{C-LS}, we have
$$
r_{\nu}/r_{\nu_0},~ r_{\nu_0}/r_{\nu},~
r_{\nu_0}/r_{\mu_0},~r_{\mu_0}/r_{\nu_0}\le \CLS.
$$
Therefore, $(\bigstar 3)$ and \rf{CLS-A} imply that
$$
d(x_{\nu},x_{\mu_0})
\le a\,r_{\nu}+\CLS a\,r_{\nu}+2a\,r_{\mu_0}
\le r_{\nu}+r_{\mu_0},
$$
and, consequently, another application of {\sc Consistency of the Lengthscale} \rf{C-LS} gives
$$
r_{\nu}/r_{\mu_0},~ r_{\mu_0}/r_{\nu}\le \CLS.
$$
\par Next, note that, by \rf{D-NU} and \rf{CLS-A},
$$
d(x,x_{\nu})\le a\,r_{\nu}<r_{\nu}
$$
and
$$
d(y,x_{\nu})\le d(y,x)+ d(x,x_{\nu})\le  [a\,r_{\nu_0}+a\,r_{\mu_0}]+a\,r_{\nu}\le
(3\CLS\,a)r_\nu< r_\nu.
$$
Hence,
$$
x,y\in B(x_\nu,r_\nu).
$$
\smallskip
\par Similarly,
$$
d(x,x_{\nu_0})\le  a\,r_{\nu_0}<r_{\nu_0}
$$
and
$$
d(y,x_{\nu_0})\le d(y,x)+ d(x,x_{\nu_0})\le  [a\,r_{\mu_0}+a\,r_{\nu_0}]+a\,r_{\nu_0}\le
(3\CLS\,a)r_{\nu_0}< r_{\nu_0}.
$$
Hence,
$$
x,y\in B(x_{\nu_0},r_{\nu_0}).
$$
\par Finally,
$$
d(x,x_{\mu_0})\le d(x,y)+ d(y,x_{\mu_0})\le  [a\,r_{\mu_0}+a\,r_{\nu_0}]+a\,r_{\mu_0}\le
(3\CLS\,a)r_{\mu_0}< r_{\mu_0}
$$
and
$$
d(y,x_{\mu_0})\le  a\,r_{\mu_0}< r_{\mu_0}.
$$
Hence,
$$
x,y\in B(x_{\mu_0},r_{\mu_0}).
$$
\par The proof of the lemma is complete.\bx

\bigskip
\par {\it Proof of the Patching Lemma \reff{PTHM}.}
\smallskip
\par We write $c,C,C'$, etc. to denote constants determined by $\CLS$, $\CWH$, $\CETA$, $\CAGR$, $\CLIP$, $\DS$. These symbols may denote different constants in different occurrences.
\par Let $x,y\in X$ be given. We must show that
$$
\|F(x)-F(y)\|\le C\,d(x,y).
$$
\par Fix $\mu_0,\nu_0$, with $x\in\TRL(\nu_0)$ and $y\in\TRL(\mu_0)$. We distinguish two cases.
\medskip
\par {\it CASE 1:} Suppose
$$
d(x,y)\le a\cdot[r_{\nu_0}+r_{\mu_0}]~~~\text{with}~~~
a=(4\,\CLS)^{-1}.
$$
Then Lemma \reff{STR} yields
\bel{E-L2}
x,y\in B(x_\nu,r_\nu)\cap B(x_{\nu_0},r_{\nu_0})\cap B(x_{\mu_0},r_{\mu_0})
\ee
for all $\nu\in \TRL(x)\cup\TRL(y)$. (If $\nu\in\TRL(y)$,
we apply Lemma \reff{STR} with $y,x,\mu_0,\nu_0$ in place of $x,y,\nu_0,\mu_0$.) Also, for such $\nu$, Lemma \reff{STR} gives
\bel{RMN}
c\,r_{\nu_0}\le r_{\nu}\le C\,r_{\nu_0}
~~~\text{and}~~~c\,r_{\nu_0}\le r_{\mu_0}\le C\,r_{\nu_0}. \ee
\par For $v\in\TRL(x)$, we have
\bel{Y-BL}
\|F_\nu(y)-\eta_{\nu_0}\|\le \|F_\nu(y)-\eta_{\nu}\|+
\|\eta_\nu-\eta_{\nu_0}\|\le C\,r_\nu+C\,[r_\nu+r_{\nu_0}+d(x_\nu,x_{\nu_0})].
\ee
(Here, we may apply {\sc Agreement of $F_\nu$ with $\eta_\nu$}, because $y\in B(x_\nu,r_{\nu})$.) Also, by \rf{E-L2},
$$
d(x_\nu,x_{\nu_0})\le d(x_\nu,x)+d(x,x_{\nu_0})\le r_{\nu}+r_{\nu_0}~~~\text{for}~~~
\nu\in\TRL(x).
$$
\par The above estimates and \rf{RMN} tell us that
$$
\|F_\nu(y)-\eta_{\nu_0}\|\le C\,r_{\nu_0}~~~\text{if}~~~\nu\in\TRL(x).
$$
\par Similarly, suppose $v\in\TRL(y)$. Then \rf{Y-BL} holds. (We may apply {\sc Agreement of $F_\nu$ with $\eta_\nu$}, because $y\in B(x_\nu,r_{\nu})$.)
Also, by \rf{E-L2},
$$
d(x_\nu,x_{\nu_0})\le d(x_\nu,y)+d(y,x_{\nu_0})\le r_{\nu}+\,r_{\nu_0}
~~\text{for all}~~\nu\in\TRL(y).
$$
\par The above estimates and \rf{RMN} tell us that
$$
\|F_\nu(y)-\eta_{\nu_0}\|\le C\,r_{\nu_0}~~~\text{for all}~~~\nu\in\TRL(y).
$$
Thus,
$$
\|F_\nu(y)-\eta_{\nu_0}\|\le C\,r_{\nu_0}~~~\text{for all}~~~\nu\in\TRL(x)\cup\TRL(y).
$$
\par We now write
$$
F(x)-F(y)=\smed_{\nu\in\TRL(x)\cup\TRL(y)}
\theta_\nu(x)\cdot[F_\nu(x)-F_\nu(y)]+
\smed_{\nu\in\TRL(x)\cup\TRL(y)}
[\theta_\nu(x)-\theta_\nu(y)]\cdot[F_\nu(y)-\eta_{\nu_0}]
\equiv I+II.
$$
\par We note that
$$
\|I\|\le \smed_{\nu\in\TRL(x)\cup\TRL(y)}
\theta_\nu(x)\cdot[C\,d(x,y)]=C\,d(x,y).
$$
\par Each summand in $II$ satisfies
$$
|\theta_\nu(x)-\theta_\nu(y)|\le C\,d(x,y)/r_\nu
~~~\text{and}~~~
\|F_\nu(y)-\eta_{\nu_0}\|\le C\,r_{\nu_0}\le C'r_\nu,
$$
hence
$$
\|[\theta_\nu(x)-\theta_\nu(y)]
\cdot[F_\nu(y)-\eta_{\nu_0}]\|
\le C''d(x,y).
$$
Because there are at most $2\DS$ summands in $II$, it follows that
$$
\|II\|\le C\,d(x,y).
$$
Combining our estimates for terms $I$ and $II$, we find that
$$
\|F(x)-F(y)\|\le C\,d(x,y)~~~\text{in CASE 1.}
$$
\bigskip
\par {\it CASE 2:} Suppose
$$
d(x,y)> a\cdot[r_{\nu_0}+r_{\mu_0}]~~~\text{with}~~~
a=(4\,\CLS)^{-1}.
$$
For $\nu\in\TRL(x)$, we have
$$
d(x_\nu,x_{\nu_0})\le  d(x_\nu,x)+d(x,x_{\nu_0}) \le a\cdot r_{\nu}+a\cdot r_{\nu_0},
$$
hence
$$
c\,r_{\nu_0}\le r_{\nu}\le C\,r_{\nu_0}
$$
and
$$
\|F_\nu(x)-\eta_{\nu_0}\|\le \|F_\nu(x)-\eta_{\nu}\|+
\|\eta_\nu-\eta_{\nu_0}\|\le C\,r_\nu+[C\,r_\nu+C\,r_{\nu_0}+Cd(x_\nu,x_{\nu_0})]
\le
C\,r_{\nu_0}.
$$
\par Consequently,
$$
\|F(x)-\eta_{\nu_0}\|=\left\|\smed_{v\in\TRL(x)}
\theta_\nu(x)\cdot[F_\nu(x)-\eta_{\nu_0}]\right\|\le
C\,r_{\nu_0}
\smed_{v\in\TRL(x)}\theta_\nu(x)=C\,r_{\nu_0}.
$$
\par Similarly,
$$
\|F(y)-\eta_{\mu_0}\|\le C\,r_{\mu_0}.
$$
Therefore,
\be
\|F(x)-F(y)\|&\le& C\,r_{\nu_0}+C\,r_{\mu_0}+\|\eta_{\nu_0}-\eta_{\mu_0}\|
\le
C'\,r_{\nu_0}+C'\,r_{\mu_0}+C'd(x_{\nu_0},x_{\mu_0})\nn\\
&\le&
C'\,r_{\nu_0}+C'\,r_{\mu_0}+C'[d(x_{\nu_0},x)+d(x,y)+
d(y,x_{\mu_0})]\nn\\
&\le&
C''\,r_{\nu_0}+C''\,r_{\mu_0}+C''d(x,y).\nn
\ee
\par Moreover, because we are in CASE 2, we have
$$
r_{\nu_0}+r_{\mu_0}\le\tfrac1a\,d(x,y)=4\,\CLS\,d(x,y).
$$
\par It now follows that
$$
\|F(x)-F(y)\|\le C'''d(x,y)~~~\text{in\, CASE}~2.
$$
\par Thus, the conclusion of the {\sc Patching Lemma} holds in all cases.\bx

\SECT{3. Basic Convex Sets, Labels and Bases.}{3}
\addtocontents{toc}{3. Basic Convex Sets, Labels and Bases. \hfill\thepage\par\VST}

\indent\par {\bf 3.1. Main properties of Basic Convex Sets.}
\addtocontents{toc}{~~~~3.1. Main properties of Basic Convex Sets.\hfill \thepage\par}\medskip
\indent\par We recall that $(Y,\|\cdot\|)$ denotes a Banach space. Given a set $S\subset\BS$ we let $\affspan(S)$ denote the affine hull of $S$, i.e., the smallest (with respect to inclusion) affine subspace of $\BS$ containing  $S$. We define the affine dimension $\dim S$ of $S$ as the dimension of its affine hull, i.e., $$\dim S=\dim\affspan (S).$$
\par Given $y\in\BS$ and $r>0$ we let
$$
\BL(y,r)=\{z\in \BS:~\|z-y\|\le r\}
$$
denote a {\it closed} ball in $\BS$ with center $y$ and radius $r$. By $\BY=\BL(0,1)$ we denote the unit ball in $\BS$.
\smallskip
\par Let $(\Mc,\rho)$ be a {\it finite} pseudometric space
with a finite pseudometric $\rho$. Let us fix a constant $\lambda>0$ and a set-valued mapping $F:\Mc\to\CNMY$. Recall that $\CNMY$ denotes the family of all non-empty convex subsets of $Y$ of dimension at most $m$.
\smallskip
\par In this section we introduce a family of convex sets $\GL(x)\subset Y$ parametrized by $x\in \Mc$ and a non-negative integer $\ell$. To do so, we first define integers $k_0,k_1,k_2,...$ by the formula
\bel{KEL}
k_\ell=(m+2)^\ell~~~(\ell\ge 0).
\ee
\begin{definition}\lbl{GML-D} {\em Let $x\in \Mc$ and let $S\subset \Mc$. A point $\xi\in Y$ belongs to the set $\Gamma(x,S)$ if there exists a mapping $f:S\cup\{x\}\to Y$ such that:
\par (i) $f(x)=\xi$ and $f(z)\in F(z)$ for all $z\in S\cup\{x\}$;
\par (ii) For every $z,w\in S\cup\{x\}$ the following inequality
$$
\|f(z)-f(w)\|\le \lambda\,\rho(z,w)
$$
holds.}
\end{definition}
\par We then define
\bel{BC-3}
\GL(x)=\bigcap_{\substack{S\subset \Mc\\\#S\le k_\ell}}
\Gamma(x,S)~~~\text{for}~~~x\in \Mc,~\ell\ge 0.
\ee
\par For instance, given $x\in\Mc$ let us present an explicit formula for $\Gamma_0(x)$. By \rf{BC-3} for $\ell=0$,
$$
\Gamma_0(x)=\bigcap_{S\subset\Mc,\,\,\#S\le 1}
\Gamma(x,S).
$$
Clearly, by Definition \reff{GML-D},
$$
\Gamma(x,\{z\})=F(x)\,\capsm
\left(F(z)+\lambda\,\rho(x,z)\BY\right)~~\text{for every}~~z\in\Mc\,,
$$
and $\Gamma(x,\emp)=F(x)$, so that
\bel{GM-ZR}
\Gamma_0(x)=\bigcap_{z\in\Mc}
\left(F(z)+\lambda\,\rho(x,z)\BY\right)\,.
\ee
\begin{remark}\lbl{GL-DFL} {\em Of course, the sets $\GL(x)$ also depend on the set-valued mapping $F$, the constant $\lambda$ and $m$. However, we use $\Gamma$'s only in this section, Sections 3-4 and Section 6.1 where these objects, i.e., $F$, $\lambda$ and $m$, are clear from the context. Therefore we omit $F$, $\lambda$ and $m$ in the notation of $\Gamma$'s.
\rbx}
\end{remark}
\smallskip
\par The above $\Gamma'\text{s}$ are (possibly empty) {\it convex subsets} of $Y$. Note that
\bel{BC-4a}
\Gamma(x,S)\subset F(x)~~~\text{for all}~~~x\in\Mc
~~~\text{and}~~~S\subset\Mc.
\ee
Hence,
\bel{AFF-4}
\Gamma(x,S)\subset \affspan(F(x)) ~~~~~~x\in\Mc,~S\subset\Mc.
\ee
\par From \rf{BC-4a} and \rf{BC-3} we obtain
\bel{GL-F2}
\GL(x)\subset F(x)~~~~\text{for}~~~x\in \Mc,~ \ell\ge 0.
\ee
Also, obviously,
\bel{BC-5}
\GL(x)\subset \GLO(x)~~~~\text{for}~~~x\in \Mc,~ \ell\ge 1.
\ee
\par We describe main properties of the sets $\GL$ in Lemma \reff{G-AB} below. The proof of this lemma relies on
Helly's intersection theorem \cite{DGK}, a classical result from the Combinatorial Geometry of convex sets.
\begin{theorem}\lbl{HT-2} Let $\Kc$ be a finite family of non-empty convex subsets of $Y$ lying in an affine subspace of $Y$ of dimension $m$. Suppose that every subfamily of $\Kc$ consisting of at most $m+1$ elements has a common point. Then there exists a point common to all of the family $\Kc$.
\end{theorem}
\begin{lemma}\lbl{G-AB} Let $\ell\ge 0$. Suppose that the restriction $F|_{\Mc'}$ of $F$ to an arbitrary subset $\Mc'\subset\Mc$ consisting of at most $k_{\ell+1}$ points has a Lipschitz selection $f_{\Mc'}:\Mc'\to Y$ with  $\|f_{\Mc'}\|_{\Lip(\Mc',Y)}\le\lambda$. Then for all $x\in \Mc$ \smallskip
\par (a)\, $\GL(x)\ne\emp$\,;\medskip
\par (b)\, $\GL(x)\subset \GLO(y)+\lambda\,\rho(x,y)\,\BY$~ for all $y\in\Mc$\,, provided $\ell\ge 1$.
\end{lemma}
\par {\it Proof.} Thanks to \rf{BC-3}, \rf{AFF-4} and Helly's Theorem \reff{HT-2}, conclusion (a) will follow if we can show that
\bel{BC-7}
\Gamma(x,S_1)\,\capsm...\capsm\,\Gamma(x,S_{m+1})\ne\emp
\ee
for every $S_1,...,S_{m+1}\subset \Mc$ such that $\#S_i\le k_\ell$ (each $i$). (We note that, by \rf{AFF-4}, each set $\Gamma(x,S)$ is a subset of the affine space $\affspan(F(x))$ of dimension at most $m$. We also use the fact that there are only finitely many $S\subset \Mc$ because $\Mc$ is finite.)\medskip
\par However, $S_1\cup...\cup S_{m+1}\cup\{x\}\subset \Mc$ has cardinality at most
$$
(m+1)\cdot k_\ell+1\le k_{\ell+1}.
$$
\par The lemma's hypothesis therefore produces a function $\tff:S_1\cup...\cup S_{m+1}\cup\{x\}\to Y$ such that
$\tff(z)\in \FF(z)$ for all $z\in S_1\cup...\cup S_{m+1}\cup\{x\}$, and
$$
\|\tff(z)-\tff(w)\|\le \lambda\,\rho(z,w)~~~\text{for all}~~~z,w\in S_1\cup...\cup S_{m+1}\cup\{x\}.
$$
Then $\tff(x)$ belongs to $\Gamma(x, S_i)$ for $i=1,...,m+1$, proving \rf{BC-7} and thus also proving (a). \medskip
\par To prove (b), let $x,y \in \Mc$, and let $\xi\in \Gamma_\ell(x)$ with $\ell\ge 1$. We must show that there exists $\eta\in \Gamma_{\ell-1}(y)$ such that
$\|\xi-\eta\|\le\lambda\cdot \rho(x,y)$. To produce such an $\eta$, we proceed as follows.
\par Given a set $S\subset \Mc$ we introduce a set $\GH(x,y,\xi,S)$ consisting of all points $\eta\in Y$ such that there exists a mapping $f:S\cup\{x,y\}\to Y$ satisfying the following conditions:\smallskip
\par (i) $f(x)=\xi$, $f(y)=\eta$, and $f(z)\in F(z)$ for all $z\in S\cup\{x,y\}$;
\par (ii) For every $z,w\in S\cup\{x,y\}$ the following inequality
$$
\|f(z)-f(w)\|\le \lambda\,\rho(z,w)
$$
holds.
\smallskip
\par Clearly, $\GH(x,y,\xi,S)$ is a {\it convex} subset of $F(y)$. Let us show that
\bel{BC-8}
\bigcap_{\substack{S\subset \Mc\\\#S\le k_{\ell-1}}}
\GH(x,y,\xi,S)\,\ne\emp\,.
\ee
Thanks to Helly's Theorem \reff{HT-2}, \rf{BC-8} will follow if we can show that
\bel{BC-9}
\GH(x,y,\xi,S_1)\cap...\cap \GH(x,y,\xi,S_{\MPL})\,\ne\emp
\ee
for all $S_1,...,S_{\MPL}\subset \Mc$ with $\#S_i\le k_{\ell-1}$ (each $i$).
\par We set $S=S_1\cup...\cup S_{\MPL}\cup\{y\}$. Then $S\subset \Mc$ with
$$
\#S\le (\MPL)\cdot k_{\ell-1}+1\le k_\ell.
$$
Because $\xi\in \GL(x)$, there exists $\tff:S_1\cup...\cup S_{\MPL}\cup\{x,y\}\to Y$ such that
$$
\tff(x)=\xi,~ \tff(z)\in \FF(z)~~~\text{for all}~~~
z\in S_1\cup...\cup S_{\MPL}\cup\{x,y\},
$$
and
$$
\|\tff(z)-\tff(w)\|\le \lambda\,\rho(z,w)~~~\text{for}~~~z,w\in S_1\cup...\cup S_{\MPL}\cup\{x,y\}.
$$
We then have $\tff(y)\in\GH(x,y,\xi,S_i)$ for $i=1,...,\MPL$, proving \rf{BC-9} and therefore also proving \rf{BC-8}.
\par Let
$$
\eta\in \bigcap_{\substack{S\subset \Mc\\\#S\le k_{\ell-1}}}
\GH(x,y,\xi,S)\,.
$$
Taking $S=\emp$, we obtain a function $\ff:\{x,y\}\to Y$ with $\ff(x)=\xi$, $\ff(y)=\eta$ and
$$
\|\ff(z)-\ff(w)\|\le \lambda\,\rho(z,w)~~~\text{for}~~~z,w\in \{x,y\}.
$$
Therefore,
\bel{BC-10}
\|\eta-\xi\|\le \lambda\,\rho(z,w).
\ee
Moreover, because $\GH(x,y,\xi,S)\subset\Gamma(y,S)$ for any $S\subset \Mc$ (see Definition \reff{GML-D}), we have
\bel{BC-11}
\eta\in \bigcap_{\substack{S\subset \Mc\\\#S\le k_{\ell-1}}}
\Gamma(y,S)=\GLO(y)\,.
\ee
\par Our results \rf{BC-10}, \rf{BC-11} complete the proof of (b).\bx

\bigskip
\indent\par {\bf 3.2. Statement of the Finiteness Theorem for bounded Nagata Dimension.}
\addtocontents{toc}{~~~~3.2. Statement of the Finiteness Theorem for bounded Nagata Dimension.\hfill \thepage\par}\medskip

\par We place ourselves in the following setting. \smallskip
\par \textbullet~We fix a positive integer $m$.\smallskip
\par \textbullet~$(X,d)$ is a finite metric space satisfying the Nagata condition (see Definition \reff{NG-C}). \smallskip
\par \textbullet~$\BS$ is a Banach space. We write $\|\cdot\|$ for the norm in $Y$, and $\|\cdot\|_{Y^*}$ for the norm in the dual space $Y^*$. We write $\ip{e,y}$ to denote the natural pairing between vectors $y\in Y$ and dual vectors $e\in Y^*$.\smallskip
\par \textbullet~For each $x\in X$ we are given a convex set
$$
\FF(x)\subset \Aff_{\FF}(x)\subset Y,
$$
where
$$
\Aff_F(x)~~~\text{is an affine subspace of}~~Y,~~\text{of dimension at most}~~m.
$$
Say, $\Aff_F(x)$ is a translate of the vector subspace $\Vect_{\FF}(x)\subset Y$.\smallskip
\par \textbullet~We make the following assumption for a large enough $\ks$ determined by $m$. \medskip
\begin{fnassumption}\lbl{FNA} Given $S\subset X$ with $\#S\le \ks$, there exists $f^S:S\to Y$ with Lipschitz constant at most $1$, such that $f^S(x)\in F(x)$ for all $x\in S$.
\end{fnassumption}
\smallskip
\par The above assumption implies the existence of a Lipschitz selection with a controlled Lipschitz constant. More precisely, we have the following result.
\begin{theorem}\lbl{PFT-FT}(Finiteness Theorem for bounded Nagata Dimension) Let $(X,\dt)$ be a finite metric space satisfying the Nagata condition with constants  $\cng$ and $\CNG$.
\par Given $m\in\N$ there exist a constant $\ks\in\N$ depending only on $m$, and a constant $\gnd>0$ depending only on $m$, $\cng$, $\CNG$, for which the following holds: Let $Y$ be a Banach space. For each $x\in X$, let $F(x)\subset Y$ be a convex set of (affine) dimension at most $m$.\smallskip
\par Suppose that for each $S\subset X$ with $\#S\le \ks$ there exists $f^S:S\to Y$ with Lipschitz constant at most $1$, such that $f^S(x)\in F(x)$ for all $x\in S$.
\par Then there exists $f:X\to Y$ with Lipschitz constant at most $\gnd$, such that $f(x)\in F(x)$ for all $x\in X$.
\end{theorem}
\par We place ourselves in the above setting until the end of the proof of Theorem \reff{PFT-FT}. See Section 4.9.
\bigskip

\par In this setting we define Basic Convex Sets following the approach suggested in Section 3.1. More specifically, let $(\Mc,\rho)=(X,\dt)$, $\lambda=1$ and let $F:X\to \Conv_m(Y)$ be the set-valued mapping from Theorem \reff{PFT-FT}. We apply Definition \reff{GML-D} and formulae \rf{KEL}, \rf{BC-3} to these objects and obtain a family $\{\GL(x):x\in X,\ell=0,1,...\}$ of convex subsets of $Y$.
\par Thus,
\bel{GM-MT}
\GL(x)=\bigcap_{\substack{S\subset X\\\#S\le k_\ell}}
\Gamma(x,S)~~~\text{for}~~~x\in X,~\ell\ge 0,
\ee
where
\smallskip
\par (i)~ $k_0,k_1,k_2,...$ is a sequence of positive integers defined by the formula
$$
k_\ell=(m+2)^\ell~~~(\ell\ge 0);
$$
\par (ii)~ $\Gamma(x,S)$ for $S\subset X$ is a subset of $Y$ defined as follows: A point $\xi\in \Gamma(x,S)$ if there exists a mapping $f:S\cup\{x\}\to Y$ such that:
\medskip
\par (a) $f(x)=\xi$ and $f(z)\in F(z)$ for all $z\in S\cup\{x\}$;
\smallskip
\par (b) For every $z,w\in S\cup\{x\}$ the following inequality
$$
\|f(z)-f(w)\|\le \dt(z,w)
$$
holds.\smallskip
\par We note that for every $x\in X$
$$
\GL(x)\subset F(x)~~\text{for}~~\ell\ge 0
$$
and $\GL(x)\subset \GLO(x)$ for $\ell\ge 1$. See \rf{GL-F2} and \rf{BC-5}.
\smallskip
\par Finally, we apply Lemma \reff{G-AB} to the setting of this section. The Finiteness Assumption \reff{FNA} enables us to replace the hypothesis of this lemma with the requirement $\ks\ge k_{\ell+1}$, which leads us to the following statement.
\begin{lemma}\lbl{BCS-L1} Let $\ell\ge 0$ and let $\ks\ge k_{\ell+1}$. Then
\smallskip
\par (A) $\GL(x)$ is nonempty for all $x\in X$;\medskip
\par (B) If $\ell\ge 1$, $\xi\in\GL(x)$ and $y\in X$, then there exists $\eta\in\GLO(y)$ such that $\|\xi-\eta\|\le d(x,y)$.
\end{lemma}
\bigskip

\indent\par {\bf 3.3. Labels and Bases.}
\addtocontents{toc}{~~~~3.3. Labels and Bases.\hfill \thepage\par\VSU}\medskip
\par A ``label'' is a finite sequence $\Ac=(e_1,e_2,...,e_s)$ of functionals $e_i\in Y^*$, with  $s\le m$.
\par We write $\#\Ac$ to denote the number $s$ of functionals $e_i$ appearing in $\Ac$. We allow the case $\#\Ac=0$, in which case $\Ac$ is the empty sequence $\Ac=(~\,)$.
\par Let $\Gamma\subset Y$ be a convex set, let $\Ac=(e_1,e_2,...,e_s)$ be a label, and let $r, C_B$ be positive real numbers. Finally, let $\zeta\in Y$.
\begin{definition}\lbl{LB-DB} {\em An {\it $(\Ac,r,C_B)$-basis} for $\Gamma$ at $\zeta$ is a sequence of $s$ vectors $v_1,...,v_s\in Y$, with the following properties:\medskip
\par (B0)~ $\zeta\in\Gamma$\,.
\medskip
\par (B1)~ $\ip{e_a,v_b}=\dl_{ab}$ (Kronecker delta) for $a,b=1,...,s$\,.
\medskip
\par (B2)~ $\|v_a\|\le C_B$ and $\|e_a\|_{Y^*}\le C_B$ for $a=1,...,s$\,.
\medskip
\par (B3)~ $\zeta+\frac{r\,v_a}{C_B}$~ and~ $\zeta-\frac{r\,v_a}{C_B}$ belong to $\Gamma$ for $a=1,...,s$\,.
}
\end{definition}
\smallskip
\par If $s\ge 1$, then of course (B3) implies (B0).
\par Let us note several elementary properties of $(\Ac,r,C_B)$-bases.
\begin{remark}\lbl{LBL-R}{\em (i) If $s=0$ then (B1), (B2), (B3) hold vacuously, so the assertion that $\Gamma$ has an $\left((~\,),r,C_B\right)$-basis at $\zeta$ means simply that $\zeta\in\Gamma$;
\smallskip
\par (ii) If $r'\le r$ and $C'_B\ge C_B$, then any $(\Ac,r,C_B)$-basis for $\Gamma$ at $\zeta$ is also an
$(\Ac,r',C'_B)$-basis for $\Gamma$ at $\zeta$;\smallskip
\par (iii) If $K\ge 1$, then any $(\Ac,r,C_B)$-basis for $\Gamma$ at $\zeta$ is also an $(\Ac,Kr,KC_B)$-basis for $\Gamma$ at $\zeta$;\smallskip
\par (iv) If $\Gamma\subset \Gamma'$, then every  $(\Ac,r,C_B)$-basis for $\Gamma$ at $\zeta$ is also an $(\Ac,r,C_B)$-basis for $\Gamma'$ at $\zeta$.}\rbx
\end{remark}
\begin{lemma}\lbl{LB-L1} (``Adding a Vector'')\, Suppose $\Gamma\subset Y$ (convex) has an $(\Ac,r,C_B)$-basis at $\xi$, where $\Ac=(e_1,e_2,...,e_s)$ and $s\le m-1$.
\par Let $\eta\in\Gm$, and suppose that
$$
\|\eta-\xi\|\ge r
$$
and
$$
\ip{e_a,\eta-\xi}=0~~~\text{for}~~~a=1,...,s.
$$
\par Then there exist $\zeta\in\Gm$ and $e_{s+1}\in Y^*$ with the following properties:\medskip
\par \textbullet~ $\|\zeta-\xi\|=\tfrac12 r$.
\medskip
\par \textbullet~ $\ip{e_a,\zeta-\xi}=0$~ for $a=1,...,s$ (not necessarily for $a=s+1$).
\medskip
\par \textbullet~ $\Gm$ has an
$(\Ac^+,r,C'_B)$-basis at $\zeta$, where $\Ac^+=(e_1,...,e_s,e_{s+1})$ and $C'_B$ is determined by $C_B$ and $m$.
\end{lemma}
\par {\it Proof.} In this proof, we write $c,C,C'$ etc. to denote constants determined by $C_B$ and $m$. These symbols may denote different constants in different occurrences.
\par Let $(v_1,...,v_s)$ be an $(\Ac,r,C_B)$-basis for $\Gm$ at $\xi$. Thus, $\xi\in\Gamma$,
\bel{LB-2}
\ip{e_a,v_b}=\dl_{ab}~~~\text{for}~~~a,b=1,...,s,
\ee
\bel{LB-3}
\|e_a\|_{Y^*}\le C_B,~~\|v_a\|\le C_B~~~\text{for}~~~ a=1,...,s,
\ee
\bel{LB-4}
\xi+\frac{r}{C_B}v_a,~~\xi-\frac{r}{C_B}\,v_a\in\Gamma~~~
\text{for}~~~a=1,...,s\,.
\ee
\smallskip
\par Let
$$
\zeta=\tau\,\eta+(1-\tau)\,\xi~~~\text{with}~~~
\tau=\tfrac12\,r\,\|\xi-\eta\|^{-1}\in(0,\tfrac12].
$$
Our hypotheses on $\xi$ and $\eta$ tell us that
\bel{LB-5}
\zeta\in\Gm,~~\|\zeta-\xi\|=\tfrac12 r,~~ \ip{e_a,\zeta-\xi}=0
~~~\text{for}~~~a=1,...,s.
\ee
Because $\eta\in\Gm$, $\Gm$ is convex, and $\tau\in(0,\tfrac12]$,~ \rf{LB-4} implies
\bel{LB-6}
\zeta+\tfrac12\tfrac{r}{C_B}\,v_a,~ \zeta-\tfrac12\tfrac{r}{C_B}\,v_a\in\Gamma~~~\text{for}~~~ a=1,...,s\,.
\ee
\par Let
\bel{LB-7}
v_{s+1}=\frac{\zeta-\xi}{\|\zeta-\xi\|}\,.
\ee
(The denominator is nonzero, by \rf{LB-5}.) Then
$$
\zeta+\|\zeta-\xi\|\,v_{s+1}=\zeta+(\zeta-\xi)=2\zeta-\xi=
2\tau\eta+(1-2\tau)\xi\in\Gm
$$
because $\xi,\eta\in\Gm$ and $\tau\in(0,\tfrac12]$.
\par Also,
$$
\zeta-\|\zeta-\xi\|\,v_{s+1}=\zeta-(\zeta-\xi)=\xi\in\Gm
\,.
$$
Recall that $\|\zeta-\xi\|=\frac12\,r$, hence the above remarks and \rf{LB-6} together yield
\bel{LB-8}
\zeta+crv_a,~\zeta-crv_a\in\Gm~~~\text{for}~~~a=1,...,s+1.
\ee
\par Also, because $\ip{e_a,\zeta-\xi}=0$ for $a=1,...,s$,
the definition of $v_{s+1}$, together with \rf{LB-2}, tells us that
\bel{LB-9}
\ip{e_a,v_b}=\dl_{ab}~~~\text{for}~~~a=1,...,s~~~\text{and}
~~~b=1,...,s+1.
\ee
\par We prepare to define a functional $e_{s+1}\in Y^*$. To do so, we first prove the estimate
\bel{LB-10}
\sum_{a=1}^{s+1}\,|\lambda_a|\le C\,\left\|\sum_{a=1}^{s+1}\lambda_a\,v_a\right\|
~~~\text{for all}~~~\lambda_1,...,\lambda_{s+1}\in\R.
\ee
To see this, we first note that \rf{LB-9} yields, for any $b=1,...,s$, the estimate
\bel{LB-11}
|\lambda_b|=\left|\ip{e_b,\sum_{a=1}^{s+1}\lambda_a\,v_a}
\right|\le
\|e_b\|_{Y^*}\cdot \left\|\sum_{a=1}^{s+1}\lambda_a\,v_a\right\|
\le C_B\,\left\|\sum_{a=1}^{s+1}\lambda_a\,v_a\right\|\,.
\ee
Consequently,
\be
|\lambda_{s+1}|&=&\left\|\lambda_{s+1}v_{s+1}\right\|
\le \left\|\sum_{a=1}^{s+1}\lambda_a\,v_a\right\|+
\sum_{a=1}^{s}|\lambda_a|\,\|v_a\|
\nn\\
&\le&
\left\|\sum_{a=1}^{s+1}\lambda_a\,v_a\right\|+
C_B\sum_{a=1}^{s}|\lambda_a|
\le (1+s\,C_B^2)\left\|\sum_{a=1}^{s+1}\lambda_a\,v_a
\right\|\,.\nn
\ee
Together with \rf{LB-11}, this completes the proof of \rf{LB-10}.\medskip
\par By \rf{LB-10} and the Hahn-Banach theorem, the linear functional
$$
\sum_{a=1}^{s+1}\lambda_a\,v_a\to\lambda_{s+1}
$$
on the span of $v_1,...,v_{s+1}$ extends to a linear functional $e_{s+1}\in Y^*$, with
\bel{LB-12}
\|e_{s+1}\|_{Y^*}\le C
\ee
and
\bel{LB-13}
\ip{e_{s+1},v_a}=\dl_{s+1,a}~~~\text{for}~~~a=1,...,s+1.
\ee
\par From \rf{LB-3}, \rf{LB-5}, \rf{LB-7}, \rf{LB-9}, \rf{LB-12}, \rf{LB-13} we have
\bel{LB-14}
\zeta\in\Gm,
\ee
\bel{LB-15}
\|e_a\|_{Y^*},~ \|v_a\|\le C~~~\text{for}~~~ a=1,...,s+1,
\ee
\bel{LB-16}
\ip{e_a,v_b}=\dl_{ab}~~~\text{for}~~~ a,b=1,...,s+1.
\ee
From  \rf{LB-8}, \rf{LB-14}, \rf{LB-15}, \rf{LB-16}, we see that $v_1,...,v_{s+1}$ form an
$((e_1,...,e_{s+1}), r, C)$-basis for $\Gm$ at $\zeta$.
\par Together with \rf{LB-5}, this completes the proof of Lemma \reff{LB-L1}.\bx
\begin{lemma}\lbl{LB-L2} (``Transporting a Basis'')\,
Given $m\in\N$ and $C_B>0$ there exists a constant $\ve_0\in(0,1]$ depending only on $m$, $C_B$, for which the following holds:
\par Suppose $\Gamma\subset Y$ (convex) has an $(\Ac,r,C_B)$-basis at $\xi_0$, where $\Ac=(e_1,e_2,...,e_s)$. Suppose $\Gm'\subset Y$ (convex) satisfies:\medskip
\par (*)~ Given any $\xi\in\Gm$ there exists $\eta\in\Gm'$ such that $\|\xi-\eta\|\le\ve_0r$.
\medskip
\par Then there exists $\eta_0\in\Gm'$ with the following properties:\medskip
\par \textbullet~ $\|\eta_0-\xi_0\|\le C\,r\,.$
\medskip
\par \textbullet~ $\ip{e_a,\eta_0-\xi_0}=0$~ for $a=1,...,s$.
\medskip
\par \textbullet~ $\Gm'$ has an $(\Ac,r,C)$-basis at $\eta_0$.
\medskip
\par Here, C is determined by $C_B$ and $m$.
\end{lemma}
\par {\it Proof.} In the trivial case $s=0$, Lemma \reff{LB-L2} holds because it simply asserts that there exists $\eta_0\in\Gm'$ such that $\|\eta_0-\xi_0\|\le C\,r$, which is immediate from (*). We suppose $s\ge 1$.
\par We take
\bel{EP0-A}
\ve_0~~\text{to be less than a small enough positive constant determined by}~~C_B~~\text{and}~~m.
\ee
At the end of our proof we can take $\ve_0$ to be, say, $\frac12$ times that small positive constant.
\smallskip
\par We write $c,C,C'$ etc. to denote constants determined by $C_B$ and $m$. These symbols may denote different constants in different occurrences.
\par Let $(v_1,...,v_s)$ be an $(\Ac,r,C_B)$-basis for $\Gm$ at $\xi_0$. Thus, $\xi_0\in\Gamma$,
\bel{LB-18}
\ip{e_a,v_b}=\dl_{ab}~~~\text{for}~~~a,b=1,...,s,
\ee
\bel{LB-19}
\|e_a\|_{Y^*}\le C_B,~~\|v_a\|\le C_B~~~\text{for}~~~ a=1,...,s,
\ee
and
\bel{LB-20}
\xi_0+c_1\sigma r\,v_a\in\Gamma~~~
\text{for}~~~a=1,...,s~~~\text{and}~~~\sigma=\pm 1\,.
\ee
\smallskip
\par Applying our hypothesis (*) to the vectors in \rf{LB-20}, we obtain vectors
$$
\zeta_{a,\sigma}\in Y~~~~(a=1,...,s,~~\sigma=\pm 1)
$$
such that
\bel{LB-21}
\xi_0+c_1\sigma r\,v_a+\zeta_{a,\sigma}\in\Gamma'~~~
\text{for}~~~a=1,...,s,~~\sigma=\pm 1,
\ee
and
\bel{LB-22}
\|\zeta_{a,\sigma}\|\le \ve_0\,r~~~
\text{for}~~~a=1,...,s,~~\sigma=\pm 1\,.
\ee
\par We define vectors
\bel{LB-23}
\eta_{00}=\frac{1}{2s}\,\sum_{a=1}^s\,\sum_{\sigma=\pm 1}
\,(\xi_0+c_1\sigma r v_a+\zeta_{a,\sigma})=
\xi_0+\frac{1}{2s}\,\sum_{a=1}^s\,\sum_{\sigma=\pm 1}
\zeta_{a,\sigma}
\ee
and
\bel{LB-24}
\tv_a=\frac{[\xi_0+c_1rv_a+\zeta_{a,1}]
-[\xi_0-c_1 r v_a+\zeta_{a,-1}]}{2c_1r}=
v_a+\left(\frac{\zeta_{a,1}-\zeta_{a,-1}}{2c_1r}\right)
\ee
for $ a=1,...,s$.
\par From \rf{LB-21} and the first equality in \rf{LB-23}, we have
$$
\eta_{00}\in \Gm'.
$$
From \rf{LB-22} and the second equality in \rf{LB-23}, we have
\bel{LB-26}
\|\eta_{00}-\xi_0\|\le \ve_0 r.
\ee
From \rf{LB-22} and the second equality in \rf{LB-24}, we have
\bel{LB-27}
\|\tv_a-v_a\|\le C\,\ve_0~~~\text{for}~~~a=1,...,s.
\ee
Also, for $b=1,...,s$ and $\hsg=\pm 1$, the first equalities in \rf{LB-23}, \rf{LB-24} give
$$
\eta_{00}+\frac{1}{s}c_1r\hsg\,\tv_b
=\frac{1}{2s}\,\sum_{a=1}^s\,\sum_{\sigma=\pm 1}
\,(\xi_0+c_1\sigma r v_a+\zeta_{a,\sigma})+
\frac{\hsg}{2s}
[(\xi_0+c_1rv_b+\zeta_{b,1})-
(\xi_0-c_1 r v_b+\zeta_{b,-1})],
$$
which exhibits $\eta_{00}+\tfrac{1}{s}c_1r\hsg\,\tv_b$ as a convex combination of the vectors in \rf{LB-21}. Consequently,
$$
\eta_{00}+c_2r\,\tv_b,~~\eta_{00}-c_2r\,\tv_b\in\Gamma'
~~~\text{for}~~~b=1,...,s,
$$
which implies that
\bel{LB-28}
\eta_{00}+c_2r\sum_{a=1}^s\,\tau_a\tv_a\in\Gm'
~~~\text{for any}~~~\tau_1,...,\tau_s\in\R~~~\text{with}~~  \sum_{a=1}^s\,|\tau_a|\le 1.
\ee
Here we use the following trivial remark on convex sets: Suppose $\xi+\eta_i,\xi-\eta_i$, $(i=1,...,I)$ belong to a convex set $\Gamma$. Then
$$
\xi+\sum_{i=1}^I\,\tau_i\eta_i\in\Gamma~~~\text{for all}~~~\tau_1,...,\tau_I\in\R~~~\text{with}~~
\sum_{i=1}^I\,|\tau_i|\le 1.
$$
\smallskip
\par From \rf{LB-18}, \rf{LB-19}, \rf{LB-27}, we have
\bel{ID-CL}
|\ip{e_a,\tv_b}-\dl_{ab}|\le C\ve_0 ~~~\text{for}~~~a,b=1,...,s.
\ee
\par We let $A$ denote the $s\times s$ matrix $A=(\ip{e_a,\tv_b})_{a,b=1}^s$. Let $I=(\delta_{ab})_{a,b=1}^s$ be the identity matrix. Given an $s\times s$ matrix $T$, we let $\|T\|_{op}$ denote the operator norm of $T$ as an operator from $\ell^2_s$ into $\ell^2_s$. Clearly, $\|T\|_{op}$ is equivalent (with constants depending only on $s$) to $\max\{|t_{ab}|:1\le a,b\le s\}$ provided $T=(t_{ab})_{a,b=1}^s$.
\par Hence, by \rf{ID-CL},
\bel{MN-S}
\|A-I\|_{op}\le C\ve_0\,.
\ee
\par We recall the standard fact from matrix algebra which states that an $s\times s$ matrix $T$ is invertible and the inequality $\|T^{-1}-I\|_{op}\le \|T-I\|_{op}/(1-\|T-I\|_{op})$ is satisfied provided $\|T-I\|_{op}<1$. Therefore, by \rf{MN-S}, for $\ve_0$ small enough (see \rf{EP0-A}), the matrix $A$ is invertible, and the following inequality
\bel{A-KL}
\|A^{-1}-I\|_{op}\le 2\,\|A-I\|_{op}
\ee
holds.
\par Let $(A^{\bf T})^{-1}=(M_{gb})_{g,b=1,...,s}$ where $A^{\bf T}$ denotes the transpose of $A$. Then
\bel{LB-29}
\ip{e_a,\sum_{b=1}^s\,M_{gb}\,\tv_b}=\dl_{ag}
~~~\text{for}~~~a,g=1,...,s.
\ee
Moreover, by \rf{MN-S} and \rf{A-KL},
\bel{LB-30}
|M_{gb}-\dl_{gb}|\le C\,\ve_0
~~~\text{for}~~~g,b=1,...,s.
\ee
\par We set
\bel{LB-31}
\hv_g=\sum_{b=1}^s\,M_{gb}\,\tv_b
~~~\text{for}~~~g=1,...,s.
\ee
\par Then \rf{LB-19}, \rf{LB-27}, \rf{LB-30}, \rf{LB-31} yield
\bel{LB-32}
\|\hv_g\|\le C~~~\text{for}~~~g=1,...,s,
\ee
while \rf{LB-29}, \rf{LB-31} give
\bel{LB-33}
\ip{e_a,\hv_g}=\dl_{ag}~~~\text{for}~~~a,g=1,...,s.
\ee
\par Moreover, \rf{LB-28}, \rf{LB-30}, \rf{LB-31} together imply that
\bel{LB-34}
\eta_{00}+c_3r\sum_{g=1}^s\,\tau_g\hv_g\in\Gm'
~~~\text{for all}~~~\tau_1,...,\tau_s~~~\text{such that each}~~  |\tau_g|\le 1.
\ee
\par To see this, we simply write the linear combination
of the $\hv_g$ in \rf{LB-34} as a linear combination of the $\tv_b$ using \rf{LB-31}, and then recall \rf{LB-28}.
\par From \rf{LB-19}, \rf{LB-26} we have
\bel{LB-35}
|\ip{e_a,\eta_{00}-\xi_0}|\le C\ve_0\, r ~~~\text{for}~~~a=1,...,s.
\ee
\par We set
\bel{LB-36}
\eta_0=\eta_{00}-
\sum_{g=1}^s\,\ip{e_g,\eta_{00}-\xi_0}\,\hv_g,
\ee
so that by \rf{LB-33},
\bel{LB-37}
\ip{e_a,\eta_0-\xi_0}=\ip{e_a,\eta_{00}-\xi_0}-
\sum_{g=1}^s\,\ip{e_g,\eta_{00}-\xi_0}\ip{e_a,\hv_g}=0
~~~\text{for}~~~a=1,...,s.
\ee
Also,
\bel{LB-38}
\|\eta_0-\xi_0\|\le \|\eta_{00}-\xi_0\|+
\sum_{g=1}^s\,|\ip{e_g,\eta_{00}-\xi_0}|\cdot\|\hv_g\|
\le C\ve_0r
\ee
by \rf{LB-26}, \rf{LB-32}, \rf{LB-35}.
\par From \rf{LB-35} and our {\it small $\ve_0$ assumption} \rf{EP0-A}, we have
$$
|\ip{e_a,\eta_{00}-\xi_0}|\le \tfrac12 c_3 r
~~~\text{for}~~~a=1,...,s,
$$
with $c_3$ as in \rf{LB-34}.
\par Therefore \rf{LB-34} and \rf{LB-36} tell us that
$$
\eta_0+c_3r
\sum_{g=1}^s\,\tau_g\hv_g\in\Gm' ~~~\text{for any}~~~
\tau_1,...,\tau_s~~~\text{such that}~~~
|\tau_g|\le\tfrac12~~~\text{for each}~~g.
$$
In particular,
\bel{LB-39}
\eta_0\in\Gm'
\ee
and
$$
\eta_0+\tfrac12 c_3r\,\hv_g,~\eta_0-\tfrac12 c_3r\,\hv_g\in\Gm'~~~\text{for}~~~g=1,...,s.
$$
\par Also, recalling \rf{LB-19}, \rf{LB-32}, \rf{LB-33}, we note that
$$
\|e_a\|_{Y^*},~\|\hv_a\|\le C~~~\text{for}~~~a=1,...,s
$$
and
\bel{LB-42}
\ip{e_a,\hv_g}=\dl_{ag}~~~\text{for}~~~a,g=1,...,s.
\ee
\par Our results \rf{LB-39},...,\rf{LB-42} tell us that $\hv_1,...,\hv_s$ form an $(\Ac,r,C)$-basis for $\Gm'$ at $\eta_0$, with $\Ac=(e_1,...,e_s)$. That's the third bullet point in the statement of Lemma \reff{LB-L2}. The other two bullet points are immediate from our results  \rf{LB-38} and \rf{LB-37}.
\par The proof of Lemma \reff{LB-L2} is complete.\bx
\SECT{4. The Main Lemma.}{4}
\addtocontents{toc}{4. The Main Lemma.\hfill\thepage \par\VST}

\indent\par {\bf 4.1. Statement of the Main Lemma.}
\addtocontents{toc}{~~~~4.1. Statement of the Main Lemma.\hfill \thepage\par}\medskip
\par Recall that $(X,d)$ is a (finite) metric space satisfying the Nagata condition with constants $\cng$ and $\CNG$.
\par For any label $\Ac=(e_1,...,e_s)$, we define
\bel{ELL-2}
\ell(\Ac)=2+3\cdot(m-\#\Ac)=2+3\cdot(m-s).
\ee
\par Note that
$$
\ell(\Ac)\ge\ell(\Ac^+)+3~~~\text{whenever}~~~ \#\Ac^+>\#\Ac.
$$
\par We also recall the definition and properties of the sets $\Gamma_\ell(x)$ introduced in Section 3.2. See \rf{GM-MT} and Lemma \reff{BCS-L1}.
\bigskip
\par We now choose the constant $\ks$ in our Finiteness Assumption \reff{FNA}. We take
\bel{DK-1}
\ks=k_{\ell^{\#}+1}=(m+2)^{\ell^{\#}+1}
\ee
as in equation \rf{KEL}, with
\bel{DK-2}
\ell^{\#}=2+3m.
\ee
\par Together with Lemma \reff{BCS-L1} and our definition of $\ell(\Ac)$, this yields the following result.
\begin{lemma}\lbl{BCS-L1P} Let $\Ac$ be a label. Then
\smallskip
\par (A) $\GL(x)\ne\emp$ for any $x\in X$ and any $\ell\le\ell(\Ac)$.
\medskip
\par (B) Let $1\le\ell\le\ell(\Ac)$, let $x,y\in X$, and let $\xi\in\GL(x)$. Then there exists $\eta\in\GLO(y)$ such that
$$
\|\xi-\eta\|\le d(x,y).
$$
\end{lemma}
\par In Sections 4.2-4.9 we will prove the following result.
\begin{mainlemma}\lbl{SML-ML} Let $x_0\in X$, $\xi_0\in Y$, $r_0>0$, $C_B\ge 1$ be given, and let $\Ac$ be a label.
\par Suppose that $\Gm_{\ell(\Ac)}(x_0)$ has an $(\Ac,\ve^{-1}r_0,C_B)$-basis at $\xi_0$, where $\ve>0$ is less than a small enough constant $\emin>0$ determined by $m$, $C_B$, $\cng$, $\CNG$.
\par Then there exists $f:\BXR\to Y$ with the following properties:\medskip
\bel{SML-A1}
\|\ff(z)-\ff(w)\|\le C(\ve)\,d(z,w)~~~\text{for all}~~~z,w\in \BXR,
\ee
\bel{SML-A2}
\|\ff(z)-\xi_0\|\le C(\ve)\,r_0~~~\text{for all}~~~z\in \BXR,
\ee
\bel{SML-A3}
\ff(z)\in\Gm_0(z)~~~\text{for all}~~~z\in \BXR.
\ee
Here $C(\ve)$ is determined by $\ve$, $m$, $C_B$, $\cng$, $\CNG$.
\end{mainlemma}
\medskip
\par We will prove the Main Lemma \reff{SML-ML} by downward induction on $\#\Ac$, starting with the case ${\#\Ac=m}$, and ending with the case $\#\Ac=0$.

\indent\par {\bf 4.2. Proof of the Main Lemma in the Base Case $\#\Ac=m$.}
\addtocontents{toc}{~~~~4.2. Proof of the Main Lemma in the Base Case $\#\Ac=m$.\hfill \thepage\par}\medskip
\par In this section, we assume the hypothesis of the Main Lemma \reff{SML-ML} in the base case ${\Ac=(e_1,...,e_m)}$. Thus, in this case $\#\Ac=m$ and $\ell(\Ac)=2$, (see \rf{ELL-2}).
\par We recall that for each $x\in X$ we have $\GL(x)\subset \FF(x)\subset \Aff_{\FF}(x)$ (all $\ell\ge0$), where $\Aff_{\FF}(x)$ is a translate of the vector space $\Vect_{\FF}(x)$ of dimension $\le m$. We write $c, C,C'$, etc. to denote constants determined by $m$, $C_B$, $\cng$, $\CNG$. These symbols may denote different constants in different occurrences.
\begin{lemma}\lbl{SML-L1} For each $z\in \BXR$, there exists
\bel{SML-1}
\eta^z\in\Gm_{1}(z)
\ee
such that
\bel{SML-2}
\|\,\eta^z-\xi_0\,\|\le C\,\ve^{-1}r_0,
\ee
\bel{SML-3}
\ip{e_a,\eta^z-\xi_0}=0~~~\text{for}~~~a=1,...,m,
\ee
\bel{SML-4}
\Gm_{1}(z)~~~\text{has an}~~~(\Ac,\ve^{-1}r_0,C)\text{-basis at}~~~\eta^z.
\ee
\end{lemma}
\par {\it Proof.} We apply Lemma \reff{LB-L2}, taking $\Gm$ to be $\Gm_{2}(x_0)$, $\Gm'$ to be
$\Gm_{1}(z)$, and $r$ to be $\ve^{-1}r_0$. To apply that lemma, we must check the key hypothesis (*), which asserts in the present case that
\bel{SML-5}
\text{Given}~~\xi\in \Gm_{2}(x_0)~~\text{there exists}~~\eta\in\Gm_{1}(z)~~\text{such that}~~\|\xi-\eta\|\le\ve_0\cdot(\ve^{-1}r_0),
\ee
where $\ve_0$ is a small enough constant determined by $C_B$ and $m$.\medskip
\par To check \rf{SML-5}, we recall Lemma \reff{BCS-L1P} (B). Given $\xi\in\Gm_{2}(x_0)$ there exists $\eta\in\Gm_{1}(z)$ such that
$$
\|\xi-\eta\|\le d(z,x_0)\le r_0~~(\text{because}~ z\in\BXR)<\ve_0\cdot(\ve^{-1}r_0);
$$
here, the last inequality holds thanks to our assumption that $\ve$ is less than a small enough constant determined by $m$, $C_B$, $\cng$, $\CNG$.
\par Thus, \rf{SML-5} holds, and we may apply Lemma \reff{LB-L2}. That lemma   provides a vector $\eta^z$ satisfying \rf{SML-1},...,\rf{SML-4}, completing the proof of Lemma \reff{SML-L1}.\bx \smallskip
\par For each $z\in\BXR$, we fix a vector $\eta^z$ as in Lemma \reff{SML-L1}. Repeating the idea of the proof of Lemma \reff{SML-L1}, we establish the following result.
\begin{lemma}\lbl{SML-L2} Given $z,w\in \BXR$, there exists a vector
\bel{SML-6}
\eta^{z,w}\in\Gm_{0}(w)
\ee
such that
\bel{SML-7}
\|\eta^{z,w}-\eta^{z}\|\le C\,\ve^{-1}d(z,w)
\ee
and
\bel{SML-8}
\ip{e_a,\eta^{z,w}-\eta^{z}}=0~~~\text{for}~~~a=1,...,m.
\ee
\end{lemma}
\par {\it Proof.} If $z=w$, we can just take $\eta^{z,w}=\eta^{z}$. Suppose $z\ne w$. Because $z,w\in\BXR$, we have $0<d(z,w)\le 2r_0$. Therefore, \rf{SML-4} tells us that
\bel{SML-9}
\Gm_{1}(z)~~~\text{has an}~~~(\Ac,\tfrac12\ve^{-1}d(z,w), C)\text{-basis at}~~\eta^z.
\ee
\par We prepare to apply Lemma \reff{LB-L2}, this time taking
$$
\Gm=\Gm_{1}(z),~~~\Gm'=\Gm_{0}(w), ~~~r=\tfrac12 \ve^{-1}d(z,w).
$$
We must verify the key hypothesis (*), which asserts in the present case that:
\medskip
\par Given any $\xi\in\Gm_{1}(z)$ there exists $\eta\in\Gm_{0}(w)$ such that
\bel{SML-10}
\|\xi-\eta\|\le\ve_0\cdot(\tfrac12 \ve^{-1}d(z,w)),
\ee
where $\ve_0$ arises from the constant $C$ in \rf{SML-9}
as in Lemma \reff{LB-L2}. In particular, $\ve_0$ depends only on $m$, $C_B$, $\cng$, $\CNG$. Therefore, our assumption that $\ve$ is less than a small enough constant determined by $m$, $C_B$, $\cng$, $\CNG$ tells us that
$$
d(z,w)<\ve_0\cdot(\tfrac12 \ve^{-1}d(z,w)).
$$
\par Consequently, Lemma \reff{BCS-L1P} (B) produces for each $\xi\in \Gm_{1}(z)$ an $\eta\in \Gm_{0}(w)$ such that
$$
\|\xi-\eta\|\le d(z,w)<\ve_0\cdot(\tfrac12 \ve^{-1}d(z,w)),
$$
which proves \rf{SML-10}.
\par Therefore, we may apply Lemma \reff{LB-L2}. That lemma provides a vector $\eta^{z,w}$ satisfying  \rf{SML-6}, \rf{SML-7}, \rf{SML-8}, and additional properties that we don't need here.
\par The proof of Lemma \reff{SML-L2} is complete.\bx
\begin{lemma}\lbl{SML-L3} Let $w\in\BXR$. Then any vector $v\in\Vect_{\FF}(w)$ satisfying $\ip{e_a,v}=0$ for $a=1,...,m$ must be the zero vector.
\end{lemma}
\par {\it Proof.} Applying \rf{SML-4}, we obtain an $(\Ac,\ve^{-1}r_0,C)$-basis $(v_1,...,v_m)$ for $\Gm_{1}(w)$ at $\eta^w$. From the definition of an $(\Ac,\ve^{-1}r_0,C)$-basis, see Definition \reff{LB-DB}, we have
\bel{SML-11}
\ip{e_a,v_b}=\dl_{ab}~~~\text{for}~~~a,b=1,...,m,
\ee
and
$$
\eta^w+c\ve^{-1}r_0v_a,~\eta^w-c\ve^{-1}r_0v_a\in
\Gm_{1}(w)\subset \FF(w)\subset\Aff_{\FF}(w)
~~~\text{for}~~~a=1,...,m,
$$
from which we deduce that
\bel{SML-12}
v_a\in \Vect_{\FF}(w)~~~\text{for}~~~a=1,...,m.
\ee
From \rf{SML-11}, \rf{SML-12} we see that
$$
v_1,...,v_m\in\Vect_{\FF}(w)
$$
are linearly independent. However, $\Vect_{\FF}(w)$ has dimension at most $m$. Therefore, $v_1,...,v_m$ form a basis for $\Vect_{\FF}(w)$. Lemma \reff{SML-L3} now follows at once from \rf{SML-11}.\bx
\medskip
\par Now let $z,w\in\BXR$. From Lemmas \reff{SML-L1} and \reff{SML-L2} we have
$$
\eta^w,\, \eta^{z,w}\in \Gm_{0}(w)\subset \FF(w)\subset\Aff_{\FF}(w),
$$
and consequently
\bel{SML-13}
\eta^w-\eta^{z,w}\in\Vect_{\FF}(w).
\ee
\par On the other hand, \rf{SML-3} and \rf{SML-8} tell us that
$$
\ip{e_a,\eta^w-\xi_0}=0,~~\ip{e_a,\eta^z-\xi_0}=0,~~
\ip{e_a,\eta^z-\eta^{z,w}}=0~~~\text{for}~~~a=1,...,m.
$$
Therefore,
\bel{SML-14}
\ip{e_a,\eta^w-\eta^{z,w}}=0~~~\text{for}~~~a=1,...,m.
\ee
\par From \rf{SML-13}, \rf{SML-14} and Lemma \reff{SML-L3}, we conclude that $\eta^{z,w}=\eta^w$. Therefore, from
\rf{SML-7}, we obtain the estimate
\bel{SML-15}
\|\eta^z-\eta^w\|\le C\ve^{-1}\,d(z,w)
~~~\text{for}~~~z,w\in\BXR.
\ee
\par We now define
$$
\ff(z)=\eta^z~~~~\text{for}~~~z\in\BXR.
$$
Then \rf{SML-1}, \rf{SML-2}, \rf{SML-15} tell us that
\bel{SML-17}
\ff(z)\in\Gm_0(z)~~~\text{for all}~~~z\in\BXR,
\ee
\bel{SML-18}
\|\ff(z)-\xi_0\|\le C\ve^{-1}r_0~~~\text{for}~~~z\in\BXR,
\ee
and
\bel{SML-19}
\|\ff(z)-\ff(w)\|\le C\ve^{-1}\,d(z,w)~~~\text{for}~~~z,w\in\BXR.
\ee
\par Our results \rf{SML-17}, \rf{SML-18}, \rf{SML-19} immediately imply the conclusions of the Main Lemma \reff{SML-ML}.
\par This completes the proof of the Main Lemma \reff{SML-ML} in the base case $\#\Ac=m$.\bx
\bigskip
\vskip 3mm
\indent\par {\bf 4.3. Setup for the Induction Step.}
\addtocontents{toc}{~~~~4.3. Setup for the Induction Step. \hfill \thepage\par}\medskip
\par Fix a label $\Ac=(e_1,...,e_s)$ with $0\le s\le m-1$.
We assume the
\begin{indhypothesis}\lbl{SIS-IH} {\em Let $x_0^+\in X$, $\xi_0^+\in Y$, $r_0^+>0$, $C_B^+\ge 1$ be given, and let $\Ac^+$ be a label such that $\#\Ac^+>\#\Ac$.
\par Then the Main Lemma \reff{SML-ML} holds, with  $x_0^+$, $\xi_0^+$, $r_0^+$, $C_B^+$, $\Ac^+$, in place of $x_0$, $\xi_0$, $r_0$, $C_B$, $\Ac$, respectively.}
\end{indhypothesis}
\par We assume the
\begin{hmla}\lbl{SIS-HMLA} {\em $x_0\in X$, $\xi_0\in Y$, $r_0>0$, $C_B\ge 1$, $\Gm_{\ell(\Ac)}(x_0)$ has an $(\Ac,\ve^{-1}r_0, C_B)$-basis at $\xi_0$.}
\end{hmla}
\par We introduce a positive constant $A$, and we make the following assumptions.
\begin{lassumption}\lbl{SIS-LAA} {\em $A$ exceeds a large enough constant determined by $m$, $C_B$, $\cng$, $\CNG$.}
\end{lassumption}
\begin{smepsassumption}\lbl{SIS-SEA} {\em $\ve$ is less than a small enough constant determined by $A$, $m$, $C_B$, $\cng$, $\CNG$.}
\end{smepsassumption}
\smallskip
\par We write $c$, $C$, $C'$, etc. to denote constants determined by $m$, $C_B$, $\cng$, $\CNG$; we write $c(A)$, $C(A)$, $C'(A)$, etc. to denote constants determined by $A$, $m$, $C_B$, $\cng$, $\CNG$; we write $c(\ve)$, $C(\ve)$, $C'(\ve)$, etc. to denote constants determined by $\ve$, $m$, $A$, $C_B$, $\cng$, $\CNG$. These symbols may denote different constants in different occurrences.
\par Note that $C(\ve)$ now has a meaning different from that in the Main Lemma \reff{SML-ML}, because $C(\ve)$ may now depend on $A$.
\medskip
\par Under the above assumptions, we will prove that there exists $\ff:\BXR\to Y$ satisfying
\bel{SIS-A1*}
\|\ff(z)-\ff(w)\|\le C(\ve)\,d(z,w)~~~\text{for all} ~~~z,w\in\BXR,
\ee
\bel{SIS-A2*}
\|\ff(z)-\xi_0\|\le C(\ve)\,r_0~~~\text{for all}~~~z\in\BXR,
\ee
\bel{SIS-A3*}
\ff(z)\in\Gm_0(z)~~~\text{for all}~~~z\in\BXR.
\ee
These conclusions differ from the conclusions \rf{SML-A1},
\rf{SML-A2}, \rf{SML-A3} of the Main Lemma \reff{SML-ML} only in that here, $C(\ve)$ may depend on $A$.
\par Once we have proven the existence of such an $f$ under the above assumptions, we then pick $A$ to be a constant determined by $m$, $C_B$, $\cng$, $\CNG$, taken large enough to satisfy the Large $A$ Assumption \reff{SIS-LAA}.
\par Once we do so, our present Small $\ve$ Assumption \reff{SIS-SEA} will follow from the small $\ve$ assumption made in the Main Lemma \reff{SML-ML}. Moreover, the conclusions \rf{SIS-A1*}, \rf{SIS-A2*}, \rf{SIS-A3*} will then imply conclusions \rf{SML-A1}, \rf{SML-A2}, \rf{SML-A3}. Consequently, we will have proven
the Main Lemma \reff{SML-ML} for $\Ac$. That will complete our downward induction on $\#\Ac$, thereby proving
the Main Lemma \reff{SML-ML} for all labels.
\bigskip
\par To recapitulate:\smallskip
\par We assume the Inductive Hypothesis \reff{SIS-IH} and the Hypotheses of the Main Lemma for the Label $\Ac$ \reff{SIS-HMLA}, and we make the Large $A$ Assumption \reff{SIS-LAA} and the Small $\ve$ Assumption \reff{SIS-SEA}.
\par Under the above assumptions, our task is to prove that there exists  $\ff:\BXR\to Y$ satisfying \rf{SIS-A1*}, \rf{SIS-A2*}, \rf{SIS-A3*}. Once we do that, the Main Lemma \reff{SML-ML} will follow.\medskip
\par We keep the assumptions and notation of this section in force until the end of the proof of the Main Lemma \reff{SML-ML}.
\bigskip\bigskip

\indent\par {\bf 4.4. A Family of Useful Vectors.}
\addtocontents{toc}{~~~~4.4. A Family of Useful Vectors. \hfill \thepage\par}\medskip
\par Recall that $\GLA(x_0)$ has an $(\Ac,\ve^{-1}r_0,C_B)$-basis at $\xi_0$.
\par Let $z\in B(x_0,\pl r_0)$. Then, thanks to our Small $\ve$ Assumption \reff{SIS-SEA}, we have
\bel{FUV-1}
d(z,x_0)\le \pl r_0<\ve_0\cdot(\ve^{-1}r_0),
\ee
where $\ve_0$ arises from $C_B,m$ as in Lemma \reff{LB-L2}.
\medskip
\par We apply that lemma, taking $\Gm=\GLA(x_0)$ and $\Gm'=\GLAO(z)$, and using \rf{FUV-1} and Lemma \reff{BCS-L1P} (B) to verify the key hypothesis (*) in Lemma \reff{LB-L2}. Thus, we obtain a vector $\eta^z\in Y$, with the following properties:
\bel{FUV-UV1}
\GLAO(z)~~~\text{has an}~~~(\Ac,\ve^{-1}r_0,C)\text{-basis
at}~~~\eta^z,
\ee
\bel{FUV-UV2}
\|\eta^z-\xi_0\|\le C\ve^{-1}r_0,
\ee
and
\bel{FUV-UV3}
\ip{e_a,\eta^z-\xi_0}=0~~~\text{for}~~~a=1,...,s.
\ee
We fix such a vector $\eta^z$ for each
$z\in B(x_0,\pl r_0)$.
\bigskip\medskip
\indent\par {\bf 4.5. The Basic Lengthscales.}
\addtocontents{toc}{~~~~4.5. The Basic Lengthscales. \hfill \thepage\par}\medskip
\begin{definition} {\em Let $x\in B(x_0,5r_0)$, and let $r>0$. We say that $(x,r)$ is OK if conditions (OK1) and (OK2) below are satisfied.\medskip
\par (OK1)~ $d(x_0,x)+5r\le 5r_0$.\smallskip
\par (OK2)~ Either condition (OK2A) or condition (OK2B) below is satisfied.\medskip
\par~~~~(OK2A)~ $\#B(x,5r)\le 1$ (i.e., $B(x,5r)$ is the singleton $\{x\}$).
\smallskip
\par~~~~(OK2B)~ For some label $\Ac^+$ with $\#\Ac^+>\#\Ac$, the following holds:\medskip
\par\hspace*{19mm} For each $w\in B(x,5r)$ there exists a vector $\zeta^w\in Y$ satisfying conditions (OK2Bi),
\par\hspace*{19mm} (OK2Bii), (OK2Biii) below:\medskip
\par\hspace*{19mm} (OK2Bi)~ $\Gm_{\ell(\Ac)-3}(w)$ has an $(\Ac^+,\ve^{-1}r,A)$-basis at $\zeta^w$.\smallskip
\par\hspace*{19mm} (OK2Bii)~ $\|\zeta^w-\xi_0\|
\le A\ve^{-1}r_0$.\smallskip
\par\hspace*{19mm} (OK2Biii)~ $\ip{e_a,\zeta^w-\xi_0}=0$ for $a=1,...,s$.
}
\end{definition}
\medskip
\par Of course (OK1) guarantees that $B(x,5r)\subset B(x_0,5r_0)$.\medskip
\par Note that $(x,r)$ cannot be OK if $r>r_0$, because then (OK1) cannot hold. On the other hand, if $x\in B(x_0,5x_0)$, then $d(x_0,x)<5r_0$, hence (OK1) holds for small enough $r$, and (OK2) holds as well (because $B(x,5r)=\{x\}$ for small enough $r$; recall that $(X,d)$ is a finite metric space). Thus, for fixed $x\in B(x_0,5r_0)$, we find that $(x,r)$ is OK if $r$ is small enough, but not if $r$ is too big.
\par For each $x\in B(x_0,5r_0)$ we may therefore
\bel{BLSC}
\text{fix a {\it basic lengthscale}}~~~r(x)>0,
\ee
such that
\bel{BL-1}
(x,r(x))~~\text{is OK, but}~~(x,2r(x))~~\text{is not OK.}
\ee
\par Indeed, we may just take $r(x)$ to be any $r'$ such that $(x,r')$ is OK and
$$
r'>\tfrac12\sup\,\{r:(x,r)~~\text{is OK}\}.
$$
\par We let $\RELX$ denote the set of all $x\in B(x_0,5r_0)$ such that
\bel{BL-RELX}
B(x,r(x))\cap B(x_0,r_0)\ne\emp.
\ee
\par Clearly,
\bel{BXR-R}
B(x_0,r_0)\subset \RELX.
\ee
\par From \rf{BL-1} and (OK1), we have
$$
d(x_0,x)+5r(x)\le 5r_0~~~\text{for each}~~~x\in B(x_0,5r_0).
$$
\begin{lemma}\lbl{BL-L1}(``Good Geometry'') Let $z_1,z_2\in B(x_0,5r_0)$. If
\bel{BL-A2}
d(z_1,z_2)\le r(z_1)+r(z_2),
\ee
then
$$
\tfrac14r(z_1)\le r(z_2)\le 4r(z_1).
$$
\end{lemma}
\par {\it Proof.} Suppose not. After possibly interchanging $z_1$ and $z_2$, we have
\bel{BL-A1}
r(z_1)<\tfrac14 r(z_2).
\ee
\par Now $(z_2,r(z_2))$ is OK (see \rf{BL-1}). Therefore it satisfies (OK1), i.e.,
$$
d(x_0,z_2)+5r(z_2)\le 5r_0.
$$
Therefore, by \rf{BL-A2},
\be
d(x_0,z_1)+5\cdot (2r(z_1))
&\le& d(x_0,z_2)+d(z_1,z_2) +10r(z_1)
\le d(x_0,z_2)+r(z_1)+r(z_2) +10r(z_1)
\nn\\
&\le& d(x_0,z_2)+\tfrac{11}{4}r(z_2)+r(z_2)
<d(x_0,z_2)+5r(z_2)\le 5r_0,\nn
\ee
i.e., $(z_1,2r(z_1))$ satisfies (OK1).
\par Moreover,
\bel{BL-A3}
B(z_1,10r(z_1))\subset B(z_2,5r(z_2)).
\ee
Indeed, if $w\in B(z_1,10r(z_1))$, then \rf{BL-A1} and \rf{BL-A2} give
$$
d(w,z_2)\le  d(w,z_1)+d(z_1,z_2)
\le 10r(z_1)+r(z_1)+r(z_2)
\le \tfrac{11}{4}r(z_2)+r(z_2)<5r(z_2),
$$
proving \rf{BL-A3}.
\par Because $(z_2,r(z_2))$ is OK, it satisfies (OK2A) or (OK2B). If $(z_2,r(z_2))$ satisfies (OK2A), then so does
$(z_1,2r(z_1))$, thanks to  \rf{BL-A3}. In that case,
$(z_1,2r(z_1))$ satisfies (OK1) and (OK2A), hence
$(z_1,2r(z_1))$ is OK, contradicting  \rf{BL-1}.
\smallskip
\par On the other hand, suppose $(z_2,r(z_2))$ satisfies (OK2B). Fix $\Ac^+$ with $\#\Ac^+>\#\Ac$ such that for every $w\in B(z_2,5r(z_2))$ there exists $\zeta^w$ satisfying
\medskip
\par \textbullet~ $\Gm_{\ell(\Ac)-3}(w)$ has an $(\Ac^+,\ve^{-1}r(z_2),A)$-basis at $\zeta^w$.
\medskip
\par \textbullet~ $\|\zeta^w-\xi_0\|\le A\ve^{-1}r_0$.
\medskip
\par \textbullet~ $\ip{e_a,\zeta^w-\xi_0}=0$~~ for $a=1,...,s$.
\medskip
\par Thanks to \rf{BL-A3} there exists such a\, $\zeta^w$ for every $w\in B(z_1,5\cdot(2r(z_1)))$. Note that, by \rf{BL-A1}, the $(\Ac^+,\ve^{-1}r(z_2),A)$-basis in the first bullet point above is also an $(\Ac^+,\ve^{-1}\cdot(2r(z_1)),A)$-basis.
\par It follows that $(z_1,2r(z_1))$ satisfies (OK2B). We have seen that $(z_1,2r(z_1))$ satisfies (OK1), so again
$(z_1,2r(z_1))$ is OK, contradicting \rf{BL-1}.
\par Thus, in all cases, our assumption that Lemma \reff{BL-L1} fails leads to a contradiction.\bx
\bigskip
\indent\par {\bf 4.6. Consistency of the Useful Vectors.}
\addtocontents{toc}{~~~~4.6. Consistency of the Useful Vectors. \hfill \thepage\par}\medskip
\par Recall the useful vectors $\eta^z$ $(z\in B(x_0,\pl r_0))$, see \rf{FUV-UV1}, \rf{FUV-UV2}, \rf{FUV-UV3}, and the set $\RELX$, see \rf{BL-RELX}. In this section we establish the following result.
\begin{lemma}\lbl{CUV-L1} Let $z_1,z_2\in \RELX$. Then
$$
\|\eta^{z_1}-\eta^{z_2}\|\le C\ve^{-1}[r(z_1)+r(z_2)+d(z_1,z_2)].
$$
\end{lemma}
\par {\it Proof.} If
$$
r(z_1)+r(z_2)+d(z_1,z_2)\ge r_0/\pl,
$$
then the lemma follows from \rf{FUV-UV2} applied to $z=z_1$ and to $z=z_2$.
\par Suppose
\bel{MN-910}
r(z_1)+r(z_2)+d(z_1,z_2)< r_0/\pl.
\ee
\par Because $z_1\in\RELX$, we have $d(z_1,x_0)\le r_0+r(z_1)$, hence
$$
d(z_1,x_0)+5\cdot(2r(z_1))\le r_0+11r(z_1)<5r_0.
$$
Thus $(z_1,2r(z_1))$ satisfies (OK1), and $B(z_1,10r(z_1))\subset B(x_0,5r_0)$.
\par Recall from \rf{FUV-UV1} that $\GLAO(z_2)$ has an $(\Ac,\ve^{-1}r_0,C)$-basis at $\eta^{z_2}$. By \rf{MN-910}, it follows that
\bel{CUV-1}
\GLAO(z_2)~~~\text{has an}~~~(\Ac,\ve^{-1}[r(z_1)+r(z_2)+d(z_1,z_2)],C)
\text{-basis}~~\text{at}~~
\eta^{z_2}.
\ee
\par Our Small $\ve$ Assumption \reff{SIS-SEA} shows that
$$
d(z_1,z_2)\le \ve_0\cdot
\ve^{-1}[r(z_1)+r(z_2)+d(z_1,z_2)],
$$
for the $\ve_0$ arising from Lemma \reff{LB-L2}, where we use the constant $C$ in \rf{CUV-1} as the constant $C_B$ in Lemma \reff{LB-L2}.
\par Therefore, by Lemma \reff{LB-L2} and Lemma \reff{BCS-L1P} (B), with
$$
\Gm=\GLAO(z_2),~~\Gm'=\Gm_{\ell(\Ac)-2}(z_1),~~
r=\ve^{-1}[r(z_1)+r(z_2)+d(z_1,z_2)],
$$
we obtain a vector
$$
\zeta\in\Gm_{\ell(\Ac)-2}(z_1)
$$
such that
\bel{444}
\|\zeta-\eta^{z_2}\|\le C\ve^{-1}[r(z_1)+r(z_2)+d(z_1,z_2)]
\ee
and
$$
\ip{e_a,\zeta-\eta^{z_2}}=0~~~\text{for}~~~a=1,...,s,
$$
hence
\bel{CUV-5}
\ip{e_a,\zeta-\eta^{z_1}}=0~~~\text{for}~~~a=1,...,s.
\ee
\par We will prove that
$$
\|\zeta-\eta^{z_1}\|\le \ve^{-1}r(z_1);
$$
\rf{444} will then imply the conclusion of Lemma \reff{CUV-L1}.
\par Suppose instead that
\bel{CUV-6}
\|\zeta-\eta^{z_1}\|> \ve^{-1}r(z_1).
\ee
\par We will derive a contradiction.
\par By \rf{FUV-UV1}, and because $r(z_1)<r_0/\pl$ (see \rf{MN-910}), we know that
\bel{CUV-7}
\Gm_{\ell(\Ac)-2}(z_1)~~~\text{has an}~~~ (\Ac,\ve^{-1}r(z_1),C)\text{-basis at}~~\eta^{z_1}.
\ee
Our results \rf{CUV-5}, \rf{CUV-6}, \rf{CUV-7} are the hypotheses of Lemma \reff{LB-L1} (``Adding a Vector''). Applying that lemma, we obtain a vector
$$
\hz\in\Gm_{\ell(\Ac)-2}(z_1),
$$
with the following properties:
\bel{CUV-9}
\|\hz-\eta^{z_1}\|=\tfrac12 \ve^{-1}r(z_1),
\ee
$$
\ip{e_a,\hz-\eta^{z_1}}=0~~~\text{for}~~~a=1,...,s;
$$
also
\bel{CUV-11}
\Gm_{\ell(\Ac)-2}(z_1)~~~\text{has an}~~~ (\Ac^+,\ve^{-1}r(z_1),C)\text{-basis at}~~\hz,
\ee
for a label of the form $\Ac^+=(e_1,...,e_s,e_{s+1})$; and 
\bel{CUV-10}
\ip{e_a,\hz-\xi_0}=0~~~\text{for}~~~a=1,...,s.
\ee
See \rf{FUV-UV3}.
\par In particular,
$$
\#\Ac^+=\#\Ac+1.
$$
\par From \rf{CUV-11} we have
\bel{CUV-13}
\Gm_{\ell(\Ac)-2}(z_1)~~~\text{has an}~~~ (\Ac^+,\ve^{-1}\cdot(2r(z_1)),\tC)\text{-basis at}~~\hz.
\ee
\par Now let $w\in B(z_1,5\cdot(2r(z_1)))$. Let $\ve_0$ arise from Lemma \reff{LB-L2} where we use $\tC$ from \rf{CUV-13} as the constant $C_B$ in Lemma \reff{LB-L2}. We have
$$
d(z_1,w)< 10r(z_1)<\ve_0\cdot(\ve^{-1}\cdot(2r(z_1))),
$$
thanks to our Small $\ve$ Assumption \reff{SIS-SEA}. Therefore, Lemma \reff{BCS-L1P} (B) allows us to verify the key hypothesis (*) in  Lemma \reff{LB-L2}, with $\Gm=\Gm_{\ell(\Ac)-2}(z_1)$, $\Gm'=\Gm_{\ell(\Ac)-3}(w)$, $r=\ve^{-1}\cdot(2r(z_1))$.
\par Applying Lemma \reff{LB-L2}, we obtain a vector
$$
\zeta^w\in\Gm_{\ell(\Ac)-3}(w)
$$
with the following properties:
\bel{CUV-15}
\|\zeta^w-\hz\|\le C\ve^{-1}\cdot(2r(z_1)),
\ee
$$
\ip{e_a,\zeta^w-\hz}=0~~~\text{for}~~~a=1,...,s+1;
$$
hence by \rf{CUV-10},
\bel{CUV-17}
\ip{e_a,\zeta^w-\xi_0}=0~~~\text{for}~~~a=1,...,s.
\ee
Also,
\bel{CUV-18}
\Gm_{\ell(\Ac)-3}(w)~~~\text{has an}~~~ (\Ac^+,\ve^{-1}\cdot(2r(z_1)),C)\text{-basis at}~~\zeta^w.
\ee
\par We have
$$
\|\zeta^w-\xi_0\|\le \|\zeta^w-\hz\|+\|\hz-\eta^{z_1}\|+
\|\eta^{z_1}-\xi_0\|
\le  C\ve^{-1}r(z_1)+\tfrac12\ve^{-1}r(z_1)+C\ve^{-1}r_0
$$
by \rf{CUV-15}, \rf{CUV-9} and \rf{FUV-UV2}.
\par Recalling that $r(z_1)<r_0/\pl$, we conclude that
\bel{CUV-19}
\|\zeta^w-\xi_0\|\le C\ve^{-1}\cdot r_0.
\ee
\par Thus, for every $w\in B(z_1,5\cdot(2r(z_1)))$, our vector $\zeta^w$ satisfies \rf{CUV-17}, \rf{CUV-18}, \rf{CUV-19}. Comparing \rf{CUV-18}, \rf{CUV-19}, \rf{CUV-17} with (OK2Bi), (OK2Bii), (OK2Biii), and recalling our Large $A$ Assumption \reff{SIS-LAA},
we conclude that (OK2B) holds for $(z_1,2r(z_1))$. We have already seen that (OK1) holds for $(z_1,2r(z_1))$. Thus $(z_1,2r(z_1))$ is OK, contradicting the defining property \rf{BL-1} of $r(z_1)$.
\par This contradiction proves that \rf{CUV-6} cannot hold, completing the proof of Lemma \reff{CUV-L1}.\bx
\bigskip\medskip
\indent\par {\bf 4.7. Additional Useful Vectors.}
\addtocontents{toc}{~~~~4.7. Additional Useful Vectors. \hfill \thepage\par}\medskip
\begin{lemma}\lbl{AUV-L1} Let $x\in B(x_0,5r_0)$, and suppose that $\#B(x,5r(x))\ge 2$.
\par Then there exist a vector $\zeta^x\in Y$ and a label $\Ac^+$ with the following properties:\medskip
\bel{AUV-1}
\#\Ac^+>\#\Ac,
\ee
\bel{AUV-2}
\Gm_{\ell(\Ac)-3}(x)~~~\text{has an}~~~ (\Ac^+,\ve^{-1}r(x),A)\text{-basis at}~~\zeta^x,
\ee
\bel{AUV-3}
\|\zeta^x-\eta^x\|\le \ve^{-1}r(x),
\ee
\bel{AUV-4}
\ip{e_a,\zeta^x-\eta^x}=0~~~\text{for}~~~a=1,...,s.
\ee
\end{lemma}
\par {\it Proof.} Recall that $(x,r(x))$ is OK. We are assuming that (OK2A) fails for $(x,r(x))$, hence (OK2B)
holds. Fix $\Ac^+$ as in (OK2B), and let $\zeta^x$ be as in (OK2B) with $w=x$. Then \rf{AUV-1}, \rf{AUV-2}, \rf{AUV-4}
hold, thanks to (OK2B); however, \rf{AUV-3} may fail in case $r(x)$ is much smaller than $r_0$. If \rf{AUV-3} holds, we are done. 
\par Suppose instead that \rf{AUV-3} fails, i.e.,
\bel{AUV-5}
\|\zeta^x-\eta^x\|> \ve^{-1}r(x).
\ee
We recall from \rf{FUV-UV1} that $\Gm_{\ell(\Ac)-1}(x)$ has an  $(\Ac,\ve^{-1}r_0,C)$-basis at $\eta^x$. We have also $r(x)\le r_0$ because $(x,r(x))$ is OK; and $\Gm_{\ell(\Ac)-1}(x)\subset\Gm_{\ell(\Ac)-3}(x)$. Therefore
\bel{AUV-6}
\Gm_{\ell(\Ac)-3}(x)~~~\text{has an}~~~ (\Ac,\ve^{-1}r(x),C)\text{-basis at}~~\eta^x.
\ee
\par From \rf{AUV-4}, \rf{AUV-5}, \rf{AUV-6} and Lemma \reff{LB-L1} (``Adding a Vector''), we obtain a vector $\hz\in Y$ and a label $\hAc$ with the following properties:
\bel{AUV-7}
\#\hAc>\#\Ac,
\ee
\bel{AUV-8}
\|\hz-\eta^x\|= \tfrac12 \ve^{-1}r(x),
\ee
\bel{AUV-9}
\ip{e_a,\hz-\eta^x}=0~~~\text{for}~~~a=1,...,s,
\ee
\bel{AUV-10}
\Gm_{\ell(\Ac)-3}(x)~~~\text{has an}~~~ (\hAc,\ve^{-1}r(x),C')\text{-basis at}~~\hz.
\ee
\par Comparing \rf{AUV-7},...,\rf{AUV-10} with \rf{AUV-1},...,\rf{AUV-4}, and recalling our Large $A$ Assumption \reff{SIS-LAA}, we see that $\hz$ and $\hAc$ have all the properties asserted for $\zeta^x$ and $\Ac^+$ in the statement of Lemma \reff{AUV-L1}.
\par Thus, Lemma \reff{AUV-L1} holds in all cases.\bx
\bigskip\medskip
\indent\par {\bf 4.8. Local Selections.}
\addtocontents{toc}{~~~~4.8. Local Selections. \hfill \thepage\par}\medskip
\begin{lemma}\lbl{LS-L1} (``Local Selections'') Given $x\in\RELX$, there exists $\ff:B(x,r(x))\to Y$ with the following properties:
\medskip
\par (I)~ $\|\ff(z)-\ff(w)\|\le C(\ve)\,d(z,w)$~~ for~ $z,w\in B(x,r(x))$.
\medskip
\par (II)~ $\ff(z)\in\Gm_0(z)$~~ for~ $z\in B(x,r(x))$.
\medskip
\par (III)~ $\|\ff(z)-\eta^x\|\le C(\ve)\,r(x)$~~ for~ $z\in B(x,r(x))$.
\medskip
\par (IV)~ $\|\ff(z)-\xi_0\|\le C(\ve)\,r_0$~~ for~ $z\in B(x,r(x))$.
\medskip
\end{lemma}
\par {\it Proof.} We proceed by cases.\medskip
\par {\it Case 1.} Suppose $\#B(x,5r(x))>1$.
\par Then Lemma \reff{AUV-L1} applies. Let $\Ac^+$, $\zeta^x$ be as in that lemma. Thus,
\bel{LS-1}
\#\Ac^+>\#\Ac,
\ee
\bel{LS-2}
\|\zeta^x-\eta^x\|\le \ve^{-1}r(x)
\ee
and
$$
\Gm_{\ell(\Ac)-3}(x)~~~\text{has an}~~~ (\Ac^+,\ve^{-1}r(x),A)\text{-basis at}~~\zeta^x\,;
$$
hence
\bel{LS-4}
\Gm_{\ell(\Ac^+)}(x)~~~\text{has an}~~~ (\Ac^+,\ve^{-1}r(x),A)\text{-basis at}~~\zeta^x,
\ee
because $\ell(\Ac)-3\ge \ell(\Ac^+)$ whenever $\#\Ac^+>\#\Ac$.
\par We recall from our Small $\ve$ Assumption \reff{SIS-SEA} that
\bel{LS-5}
\ve~~\text{is less than a small enough constant determined by}~A,\cng,\CNG,m.
\ee
\par Thanks to \rf{LS-4}, \rf{LS-5}, the Hypotheses of the Main Lemma \reff{SIS-HMLA} are satisfied, with $\Ac^+$, $x$, $\zeta^x$, $r(x)$, $A$, in place of
$\Ac$, $x_0$, $\xi_0$, $r_0$, $C_B$, respectively. Moreover, thanks to \rf{LS-1} and the Inductive Hypothesis \reff{SIS-IH}, we are assuming the validity of the Main Lemma \reff{SML-ML} for $\Ac^+,...,A$.
\par Therefore, we obtain a function $f:B(x,r(x))\to Y$ satisfying (I), (II) and the inequality
$$
\|\ff(z)-\zeta^x\|\le C(\ve)\,r(x),~~~~ z\in B(x,r(x)).
$$
This inequality together with \rf{LS-2} implies
(III).
\par Moreover, (IV) follows from (III) because, for $z\in B(x,r(x))\subset B(x_0,5r_0)$, we have
$$
\|\ff(z)-\xi_0\|\le \|\ff(z)-\eta^x\|+\|\eta^x-\xi_0\|
\le C(\ve)r(x)+C\ve^{-1}r_0\le C'(\ve)r_0;
$$
here we use \rf{FUV-UV2} and the fact that $(x,r(x))$ satisfies (OK1).
\par This completes the proof of Lemma \reff{LS-L1} in Case 1.\bigskip
\par {\it Case 2.} Suppose $\#B(x,5r(x))\le 1$.
\par Then, $B(x,5r(x))=\{x\}$ and $\eta^x\in\GLAO(x)\subset\Gamma_0(x)$. Hence the function
$f(x)=\eta^x$ satisfies (I),(II),(III), and also (IV) thanks to \rf{FUV-UV2}.
\par Thus, Lemma \reff{LS-L1} holds in all cases.\bx
\bigskip\medskip
\par {\bf 4.9. Proof of the Main Lemma: the final step.}
\addtocontents{toc}{~~~~4.9. Proof of the Main Lemma: the final step. \hfill \thepage\par\VSU}\medskip
\indent\par Let $\Bc_0$ be the metric space
$$
\Bc_0=\left(B(x_0,r_0),d|_{B(x_0,r_0)\times B(x_0,r_0)}\right),
$$
i.e., the ball $B(x_0,r_0)$ supplied with the metric $d$.
\par For the rest of this section, we work in the metric space $\Bc_0$. Given $x\in B(x_0,r_0)$ and $r>0$, we write $\BO(x,r)$ to denote the ball in $\Bc_0$ with center $x$ and radius $r$; thus $\BO(x,r)= B(x,r)\cap B(x_0,r_0)$.
\par Note that the Nagata condition for $\Bc_0$ holds with the same constants $\cng$ and $\CNG$ as for $(X,d)$. See Definition \reff{NG-C}.
\par Let $r:X\to\R_+$ be {\it the basic lengthscale} constructed in Section 4.5 (see \rf{BLSC}), and let
\bel{A-CL}
\CLS=4 ~~~\text{and}~~~a=(4\,\CLS)^{-1}.
\ee
Note that, by Lemma \reff{BL-L1}, {\sc Consistency of the Lengthscale} (see \rf{C-LSC}) holds for the lengthscale $r(x)$ on $B(x_0,r_0)$ with the constant $\CLS$ given by \rf{A-CL}.
\par We apply the Whitney partition Lemma \reff{WPL} to the metric space $\Bc_0$, the lengthscale
$$\{r(x):x\in B(x_0,r_0)\}$$
and the constants $\CLS$, $a$ determined by \rf{A-CL}, and obtain a partition of unity ${\{\theta_\nu:
B(x_0,r_0)\to\R_+\}}$ and points 
\bel{XV-BR}
x_\nu\in B(x_0,r_0)
\ee
with the following properties.
\medskip
\par \textbullet~ Each $\theta_\nu\ge 0$ and for each $\nu$, $\theta_\nu=0$ outside $\BO(x_\nu,a r_\nu)$; here
$a$ is determined by \rf{A-CL}, and $r_\nu=r(x_\nu)$.
\smallskip
\par \textbullet~ Any given $x$ satisfies $\theta_\nu(x)\ne 0$ for at most $\DS$ distinct $\nu$, where $\DS$ depends only on $\cng$, $\CNG$.
\smallskip
\par \textbullet~ $\smed\limits_\nu\,\theta_\nu(x)=1$ for all $x\in B(x_0,r_0)$.
\smallskip
\par \textbullet~ Each $\theta_\nu$ satisfies
$$
|\theta_\nu(x)-\theta_\nu(y)|\le \frac{C}{r_\nu}\,d(x,y)
$$
for all $x,y\in B(x_0,r_0)$; here again $r_\nu=r(x_\nu)$.
\smallskip
\par From Lemma \reff{BL-L1} (``Good Geometry''), we know that
\smallskip
\par \textbullet~ For each $\mu,\nu$, if $d(x_\mu,x_\nu)\le r_\mu+r_\nu$, then $\frac14 r_\nu\le r_\mu\le 4r_\nu$.
\smallskip
\par Moreover, by \rf{BXR-R} and \rf{XV-BR},
\bel{RELX-N}
x_\nu\in\RELX~~~\text{for each}~~\nu,
\ee
so that, by Lemma \reff{LS-L1}, there exists a function $\hf_\nu:B(x_\nu,r_\nu)\to Y$ satisfying the following conditions
\medskip
\par \textbullet~ $\|\hf_\nu(z)-\hf_\nu(w)\|\le C(\ve)\,d(z,w)$~~ for~ $z,w\in B(x_\nu,r_\nu)$.
\smallskip
\par \textbullet~ $\hf_\nu(z)\in\Gm_0(z)$~~ for~ $z\in B(x_\nu,r_\nu)$.
\smallskip
\par \textbullet~ $\|\hf_\nu(z)-\eta_\nu\|\le C(\ve)\,r_\nu$~~ for~ $z\in B(x_\nu,r_\nu)$, where $\eta_\nu\equiv \eta^{x_\nu}$.
\smallskip
\par \textbullet~ $\|\hf_\nu(z)-\xi_0\|\le C(\ve)\,r_0$~~ for~ $z\in B(x_\nu,r_\nu)$.
\medskip
\par Let $f_\nu=\hf|_{\BO(x_\nu,r_\nu)}$. We extend $f_\nu$ from $\BO(x_\nu,r_\nu)=B(x_\nu,r_\nu)\cap B(x_0,r_0)$ to all of $B(x_0,r_0)$ by setting $f_\nu=0$ outside $\BO(x_\nu,r_\nu)$.
\smallskip
\par Since each $x_\nu\in\RELX$ (see \rf{RELX-N}), from Lemma \reff{CUV-L1}, we have
\smallskip
\par \textbullet~ $\|\eta_\nu-\eta_\mu\|\le C(\ve)\cdot
[r_\nu+r_\mu+d(x_\nu,x_\mu)]$ for each $\mu,\nu$.
\medskip
\par The above conditions on the $\theta_\nu$, $\eta_\nu$, $\hf_\nu$, $f_\nu$, $r_\nu$ and $a$ (cf. \rf{A-SMALL} with \rf{A-CL}) allow us to apply the Patching Lemma
\reff{PTHM}. We conclude that
$$
f(x)=\smed_\nu\,\theta_\nu(x)\,f_\nu(x)~~~
(\text{all}~~x\in B(x_0,r_0))
$$
satisfies
$$
\|f(x)-f(y)\|\le C(\ve)\,d(x,y)~~~\text{for}~~
x,y\in B(x_0,r_0).
$$
\par Moreover, for fixed $x\in B(x_0,r_0)$, we know that $f(x)$ is a convex combination of finitely many values $f_\nu(x)$ with $\BO(x_\nu,ar_\nu)\ni x$; for those $\nu$ we have $f_\nu(x)\in \Gamma_0(x)$ and $\|f_\nu(x)-\xi_0\|\le C(\ve)\,r_0$. Therefore, $f(x)\in \Gamma_0(x)$ and
$\|f(x)-\xi_0\|\le C(\ve)\,r_0$ for all $x\in B(x_0,r_0)$.
\medskip
\par Thus, $f$ satisfies \rf{SIS-A1*}, \rf{SIS-A2*} and \rf{SIS-A3*}, completing the proof of the Main Lemma \reff{SML-ML}.\bx
\medskip
\medskip
\par {\it Proof of the Finiteness Theorem \reff{PFT-FT} for Bounded Nagata Dimension.} Let $x_0\in X$, $r_0=\diam X+1$, $C_B=1$, and $\Ac=(~)$. Let $\ve=\tfrac12\,\emin$ where $\emin$ is as in the Main Lemma \reff{SML-ML} for $m$, $C_B=1$, $\cng$ and $\CNG$. Thus, $\ve$ depends only on $m$, $\cng$ and $\CNG$.
\par By Lemma \reff{BCS-L1P} (A), $\Gm_{\ell(\Ac)}(x_0)\ne\emp$ so that there exists $\xi_0\in \Gm_{\ell(\Ac)}(x_0)$. Since $\#\Ac=0$, the set
$\Gm_{\ell(\Ac)}(x_0)$ has an $(\Ac,\ve^{-1}r_0,C_B)$-basis at $\xi_0$. See Remark \reff{LBL-R}, (i).
\par Hence, by the Main Lemma \reff{SML-ML}, there exists a mapping $f:\BXR\to Y$ such that
$$
\|f(z)-f(w)\|\le C\,d(z,w)~~~\text{for all}~~~z,w\in \BXR,
$$
and
$$
f(z)\in\Gm_0(z)~~~\text{for all}~~~z\in \BXR.
$$
Here $C$ is a constant determined by $\ve$, $m$, $C_B$, $\cng$, $\CNG$. Thus, $C$ depends only on $m$, $\cng$, $\CNG$. 
\par Clearly, $\BXR=X$. Furthermore, $\Gm_0(z)\subset F(z)$ for every $z\in X$ (see \rf{GL-F2}), so that $f(z)\in F(z)$, $z\in X$. Thus, $f$ is a selection of $F$ on $X$
with Lipschitz constant at most a certain constant depending only on $m$, $\cng$, $\CNG$.
\par The proof of Theorem \reff{PFT-FT} is complete.\bx
\medskip
\smallskip
\par {\it Proof of Theorem \reff{NG-FP}.} The proof is immediate from Theorem \reff{PFT-FT} applied to the metric space $(X,\lambda d)$.\bx
\smallskip
\par Let us apply Theorem \ref{NG-FP} to {\it metric trees}. We recall that, by Lemma \reff{MTR-ND}, each {\it metric tree} is a finite metric space satisfying the Nagata condition with $\cng=1/16$ and $\CNG=1$. Thus, we obtain the following
\begin{corollary}\lbl{A-12} Let $m\in\N$, let $(X,d)$ be a metric tree and let $\lambda$ be a positive constant. Let ${F:X\to\CNMY}$ be a set-valued mapping such that, for every subset $X'\subset X$ with $\#X'\le \ks$, the restriction $F|_{X'}$ has a Lipschitz selection $f_{X'}:X'\to\BS$ with $\|f_{X'}\|_{\Lip(X',\BS)}\le \lambda$.
\par Then $F$ has a Lipschitz selection $f:X\to\BS$ with $\|f\|_{\Lip(X,\BS)}\le \gamma_0\,\lambda$.
\par Here $\ks=\ks(m)$ is the constant from Theorem \reff{NG-FP}, and $\gamma_0=\gamma_0(m)$ is a constant depending only on $m$.
\end{corollary}

\SECT{5. Metric trees and Lipschitz selections with respect to the Hausdorff distance.}{5}
\addtocontents{toc}{5. Metric trees and Lipschitz selections with respect to the Hausdorff distance.
\hfill\thepage\par\VST}

\indent\par We recall that $(Y,\|\cdot\|)$ denotes a Banach space, and $\Kc(\BS)$ denotes the family of all non-empty compact convex finite dimensional subsets of $\BS$. We also recall that given a non-negative integer $m$ we let $\Kc_m(\BS)$ denote the family of all sets $K\in\Kc(\BS)$ with $\dim K\le m$. By $\AM$ we denote the family of all {\it affine} subspaces of $\BS$ of dimension at most $m$.
\par Let us fix some additional notation. By $\CNY$ we denote the family of {\it all non-empty convex finite dimensional} subsets of $\BS$; thus,
\bel{CNY-D}
\CNY=\bigcup_{m=0}^\infty\,\CNMY\,.
\ee
Recall that $\CNMY$ is the family of all non-empty convex finite dimensional subsets of $\BS$ of affine dimension at most $m$.
\par Given sets $S_1,S_2\subset \BS$ we let $\dhf(S_1,S_2)$ denote the Hausdorff distance between these sets:
\bel{HD-DF}
\dhf(S_1,S_2)=\inf \{r>0: S_1+\BL(0,r)\supset S_2,~ S_2+\BL(0,r)\supset S_1\}\,.
\ee
\par In this section we work with finite trees $T=(X,E)$, where $X$ denotes the set of nodes and $E$ denotes the set of edges of $T$. We use the same notation as in Section 2. More specifically, we write $u\je v$ to indicate that $u,v\in X$ are distinct nodes joined by an edge in $T$; we denote that edge by $[uv]$. \smallskip
\par We supply $X$ with a metric $d$ defined by formulae \rf{D-XX} and \rf{D-TR}, and we refer to the metric space $(X,d)$ as a {\it metric tree (with respect to the tree $T=(X,E)$)}.
\begin{remark} {\em Sometimes we will be looking simultaneously at two different pseudometrics, say $\rho$ and $\trh$, on a pseudometric space, say on $\Mc$. In this case we will speak of a $\rho$-Lipschitz selection and $\rho$-Lipschitz seminorm or a $\trh$-Lipschitz selection and $\trh$-Lipschitz seminorm to make clear which pseudometric we are using. Furthermore, sometimes given a mapping $f:\Mc\to Y$ we will write $\|f\|_{\Lip((\Mc,\rho),Y)}$ to denote the Lipschitz seminorm of $f$ with respect to the pseudometric $\rho$.
\par Sometimes we will be dealing with two different trees
$T,\tT$. We will then say $x\je y$ {\it in} $T$ or $x\je y$ {\it in} $\tT$ to make clear which tree we are talking about.\rbx}
\end{remark}

\bigskip

\bigskip


\indent\par {\bf 5.1. Lipschitz selection orbits.}
\addtocontents{toc}{~~~~5.1. Lipschitz selection orbits. \hfill \thepage\par}\medskip
\par Let $(\Mc,\dm)$ be a pseudometric space, and let $F:\Mc\to\CNY$ be a set-valued mapping, see \rf{CNY-D}.
\begin{definition}\lbl{ORB-D} {\em Let $x\in \Mc$, $\lambda>0$, and let $V=[(\Mc,\dm),F,\lambda]$. By $\Orb(x;V)$ we denote the subset of $Y$ defined by
$$
\Orb(x;V)=\{f(x):~f~~\text{is a}~\rho\text{-Lipschitz selection of} ~~F~~\text{with}~~\|f\|_{\Lip((\Mc,\rho),\BS)}\le\lambda\}\,.
$$
}
\end{definition}
\par We refer to the set $\Orb(x;V)$ as a {\it Lipschitz selection orbit} at $x$ with respect to the tuple $V$.
\medskip
\par Of course, in general the orbit $\Orb(x;V)$ may be  empty.
\par In the sequel we will need the following useful  properties of Lipschitz selection orbits.
\begin{lemma}\lbl{ORB-M} Let $(\Mc,\dm)$ be a finite pseudometric space with a finite pseudometric $\rho$, and let $V=[(\Mc,\dm), F,\lambda]$. Then for every $x\in \Mc$ the orbit $\Orb(x;V)$ is a convex finite dimensional subset of $F(x)$. Furthermore, if for each $u\in \Mc$ the set $F(u)$ is compact, then $\Orb(x;V)$ is compact as well.
\end{lemma}
\par {\it Proof.} The convexity of $\Orb(x;V)$ directly follows from the convexity of sets $F(u)$ $(u\in \Mc)$ and Definition \reff{ORB-D}. Furthermore, if $f:\Mc\to Y$ is a selection of $F$, then $f(x)\in F(x)$ proving that  $\Orb(x;V)\subset F(x)$. This also proves that $\dim \Orb(x;V)\le \dim F(x)$ so that $\Orb(x;V)$ is a finite dimensional subset of $Y$.
\par Let us prove that $\Orb(x;V)$ is compact whenever each set $F(u), u\in \Mc$, is. Since $\Orb(x;V)\subset F(x)$ and $F(x)$ is a compact set, the orbit $\Orb(x;V)$ is a bounded set. We prove that  $\Orb(x;V)$ is closed.
\par Let $h\in\BS$, and a let $h_n\in \Orb(x;V), n=1,2,...$ be a sequence of points converging to $h$:
\bel{HN-L1}
h=\lim_{n\to\infty}h_n\,.
\ee
We will prove that $h\in \Orb(x;V)$.
\par By Definition \reff{ORB-D}, there exists a sequence of mappings $f_n\in \Lip(\Mc,Y)$ such that
\bel{FN-OR1}
f_n(u)\in F(u)~~~\text{and}~~~\|f_n\|_{\Lip(\Mc,\BS)}\le\lambda
\ee
for every $u\in \Mc$ and $n\in\N$, and
\bel{HN-FNA1}
h_n=f_n(x),~~~n=1,2,...\,.
\ee
\par Recall that $(\Mc,\dm)$ is a {\it finite} pseudometric space, and each set $F(u), u\in\Mc$, is a finite dimensional compact subset of $Y$. Therefore, there exists a subsequence $n_k\in\N$, $k=1,2,...$, such that $(f_{n_k}(u))_{k=1}^\infty$ converges in $Y$ for every $u\in \Mc$. Let
\bel{TF-DOR1}
\tf(u)=\lim_{k\to\infty}f_{n_k}(u), ~~~u\in \Mc.
\ee
\par Then, by \rf{HN-L1} and \rf{HN-FNA1},
\bel{H-FA1}
h=\lim_{k\to\infty}h_{n_k}
=\lim_{k\to\infty}f_{n_k}(x)=\tf(x).
\ee
Since each set $F(u)$ is closed, by \rf{FN-OR1} and
\rf{TF-DOR1}, $\tf(u)\in F(u)$ for every $u\in \Mc$, proving that $\tf$ is a {\it selection} of the set-valued mapping $F$ on $\Mc$. Since each $f_n:\Mc\to\BS$ is $\rho$-Lipschitz with  $\|f_n\|_{\Lip((\Mc,\rho),\BS)}\le\lambda$, by \rf{TF-DOR1}, $\tf$ is $\rho$-Lipschitz as well, with $\|\tf\|_{\Lip((\Mc,\rho),\BS)}\le\lambda$.
\par Thus, by \rf{H-FA1} and Definition \reff{ORB-D}, $h\in \Orb(x;V)$ proving the lemma.\bx
\bigskip\medskip

\indent\par {\bf 5.2. Intersection of orbits and the Finiteness Principle.}
\addtocontents{toc}{~~~~5.2. Intersection of orbits and the Finiteness Principle. \hfill \thepage\par}\medskip
\par In this and the next subsection we prove Theorem \reff{HDS-M}.
\par Until the end of the paper we write $\ks$ to denote the constant defined by the formulae \rf{DK-1}, \rf{DK-2}, and we write $\gamma_0$ to denote the constant $\gamma_0(m)$ from Corollary \reff{A-12}.
\medskip
\par Let $(\MS,\dm)$ be a metric space and let $F:\MS\to\KM$ be a set-valued mapping. We suppose that the following assumption is satisfied.
\begin{assumption}\lbl{A-F} For every subset $\MS'\subset\MS$ consisting of at most $\ks$ points, the restriction $F|_{\MS'}$ of $F$ to $\MS'$ has a $\dm$-Lipschitz selection $f_{\MS'}:\MS'\to\BS$ with  $\|f_{\MS'}\|_{\Lip((\MS',\dm),\BS)}\le 1$.
\end{assumption}
\par Our aim is to prove the existence of a mapping $G:\Mc\to\KMY$ satisfying conditions (i) and (ii) of Theorem \reff{HDS-M}.\smallskip
\par Let $T=(X,E)$ be an arbitrary finite tree. We introduce the following
\begin{definition} {\em A mapping $W:X\to\MS$ is said to be {\it admissible} with respect to $T$ if for every two distinct nodes $u,v\in X$ with $u\je v$ (i.e., $u$ is joined by an edge to $v$), we have $W(u)\ne W(v)$.}
\end{definition}
\medskip
\par Let $W:X\to\MS$ be an admissible mapping. Then $W$ gives rise a {\it tree metric} $d_{T,W}:X\times X\to \R_+$ defined by 
$$
d_{T,W}(u,v)=\dm(W(u),W(v))~~~\text{for every}~~ u,v\in X,~ u\je v\,.
$$
See \rf{D-TR}.
\par Clearly, by the triangle inequality,
\bel{RH-D}
\dm(W(u),W(v))\le d_{T,W}(u,v)~~~
\text{for every}~~ u,v\in X\,.
\ee
\par Now define a set-valued mapping $F_{T,W}:X\to\Kc_m(\BS)$ by the formula
$$
F_{T,W}(u)=F(W(u)),~~~u\in X.
$$
\begin{lemma}\lbl{S-T} The set-valued mapping $F_{T,W}=F\circ W$ has a $d_{T,W}$-Lipschitz selection $f:X\to\BS$ such that
\bel{LS-FT}
\|f\|_{\Lip((X,d_{T,W}),\BS)}\le \gamma_0\,.
\ee
\end{lemma}
\par {\it Proof.}  Let $X'\subset X$ be an arbitrary subset of $X$ with $\# X'\le \ks$, and  let $\MS'=W(X')$. Then $\#\MS'\le \#X'\le \ks$ so that, by Assumption \reff{A-F}, the restriction $F|_{\MS'}$ has a $\dm$-Lipschitz selection $f_{\MS'}:\MS'\to \BS$ with  $\|f_{\MS'}\|_{\Lip((\MS',\dm),\BS)}\le 1$.
\par Let
$$
g_{X'}(u)=f_{\MS'}(W(u)),~~~u\in X'.
$$
Then $g_{X'}$ is a {\it selection} of the restriction $F_{X,W}|_{X'}$. Furthermore, for every $u,v\in X'$
$$
\|g_{X'}(u)-g_{X'}(v)\|=
\|f_{\MS'}(W(u))-f_{\MS'}(W(v))\|\le
\dm(W(u),W(v))
$$
so that, by \rf{RH-D},
$$
\|g_{X'}(u)-g_{X'}(v)\|\le d_{T,W}(u,v)
$$
proving that the $d_{T,W}$-Lipschitz seminorm of $g_{X'}$ is bounded by $1$.
\par Hence, by Corollary \reff{A-12}, the set-valued mapping $F_{T,W}$ has a $d_{T,W}$-Lipschitz selection $f:X\to \BS$ satisfying inequality \rf{LS-FT}.\bx
\medskip
\par We will need the following two definitions.
\begin{definition} {\em Let $x\in\MS$. The  family $\Ac(x)$ consists of all triples $L=(a,(X,E),W)$ where
\smallskip
\par \textbullet~ $T=(X,E)$ is a finite tree with the family of nodes $X$ and the family of edges $E$;
\smallskip
\par \textbullet~ $a\in X$ is a node of $T$;
\smallskip
\par \textbullet~ $W:X\to\MS$ is an admissible mapping with respect to $T$ such that $W(a)=x$.
}
\end{definition}
\begin{definition}\lbl{OR-D} {\em \par Given a triple
$$
L=(a,(X,E),W)\in \Ac(x)
$$
we let $LS(L)$ denote a family of all mappings $f:X\to\BS$ such that
$$
f~\,\text{is a}~~d_{T,W}\text{-Lipschitz selection of}~~ F_{T,W}~\,\text{with} ~~\|f\|_{\Lip((X,d_{T,W}),\BS)}\le\gamma_0\,.
$$
\par We let $O(x;L)$ denote the subset of $\BS$ defined by
\bel{OL-X}
O(x;L)=\{f(a): f\in LS(L)\}\,.
\ee
}
\end{definition}
\par Clearly, by Definition \reff{ORB-D} and Definition \reff{OR-D}, given $x\in\Mc$ and $$L=(a,(X,E),W)\in \Ac(x)$$ we have
$$
O(x;L)=\Orb(a;V)~~~\text{where}~~~V=[(X,d_{T,W}),F\circ W,\gamma_0]\,.
$$
\par This observation, Lemma \reff{S-T} and Lemma \reff{ORB-M} imply the following result.
\medskip
\begin{lemma}\lbl{OR-PR} Let $x\in\MS$ and let  $L=(a,X,W)\in \Ac(x)$. Then $O(x;L)$ is a non-empty compact convex subset of $F(x)$ of affine dimension at most $m$.
\end{lemma}
\medskip
\par Given $x\in\MS$ let
\bel{GX-D}
G(x)=\bigcap_{L\,\in\,\Ac(x)}\,O(x;L)\,.
\ee
Clearly, by Lemma \reff{OR-PR}, for every $x\in\MS$ the set 
$$
G(x)~~\text{\it is a convex compact subset of}~~F(x).
$$
\par In the next section, we will prove that $G(x)\ne\emp$ for each $x\in\MS$ and that
\bel{GX-HD}
\dhf(G(x),G(y))\le \gamma_0\,\dm(x,y)~~~
\text{for every}~~x,y\in\MS\,.
\ee
Recall that $\dhf$ denotes the Hausdorff distance between subsets of $\BS$.
\bigskip\vskip 3mm
\indent\par {\bf 5.3. The Hausdorff distance between orbits.}
\addtocontents{toc}{~~~~5.3. The Hausdorff distance between orbits.\hfill \thepage\par\VSU}
\begin{lemma}\lbl{L1} For every $x\in\MS$\,, the set $G(x)\ne\emp$.
\end{lemma}
\par {\it Proof.} We must show that
$$
\bigcap_{L\,\in\,\Ac(x)}\,O(x;L)\ne\emp.
$$
See \rf{GX-D}. By Lemma \reff{OR-PR}, each $O(x;L)$ is a non-empty compact subset of the compact set $F(x)$. Therefore, it is enough to show that
\bel{O-LI2}
O(x;L_1)\cap...\cap O(x;L_N)\,\ne\emp
\ee
for every finite subcollection $\{L_1,...,L_N\}\subset\Ac(x)$.
\par Let $L_1,...,L_N\in \Ac(x)$ with $L_i=(a_i,(X_i,E_i),W_i)$,~$i=1,...,N$.
\medskip
\par We introduce a procedure for gluing finite trees $T_i=(X_i,E_i)$,~$i=1,...,N$, together. Recall that $X_i$ here denotes the set of nodes of $T_i$, and $E_i$ denotes the set of edges of $T_i$. By passing to isomorphic copies of the $T_i$, we may assume that the sets $X_i$ are pairwise disjoint.
\smallskip
\par For each $i=1,...,N$, let $a_i$ be a node of $T_i$.
Then we form a finite tree $T^+=(X^+, E^+)$ from $T_1,...,T_N$ by identifying together all the nodes $a_1,...,a_N$.  We spell out details below.\smallskip
\par For each $i$, we write $J_i$ to denote the set
$J(a_i;T_i)$ of all the neighbors of $a_i$ in $T_i$. Also, we write $X'_i$ to denote the set $X_i\setminus\{a_i\}$, and we write $E'_i$ to denote all the edges in $T_i$ that join together points of $X'_i$ (i.e. not including $a_i$ as an endpoint).\smallskip
\par We introduce a new node $a^+$ distinct from all the nodes of all the $T_i$.\smallskip
\par The finite tree $T^+=(X^+,E^+)$ is then defined as follows. The nodes $X^+$ are all the nodes in all the $X'_i$, together with the single node $a^+$. The edges $E^+$ are all the edges belonging to any of the $E'_i$, together with edges joining $a^+$ to all the nodes in all the $J_i$. One checks easily that $T^+$ is a finite tree. We say that $T^+$ arises by ``gluing together the $T_i$ by identifying the $a_i$''.\smallskip
\par Note that $T^+$ contains an isomorphic copy of each $T_i$ as a subtree; the relevant isomorphism $\varphi_i$ carries the node $a_i$ of $T_i$ to the node $a^+$ of $T^+$, and $\varphi_i$ is the identity on all other nodes of $T_i$.
\smallskip
\par This concludes our discussion of the gluing of trees $T_i$.
\medskip
\par  We define a map $W^+:X^+\to\Mc$ by setting
\bel{WP-1}
W^+(a^+)=x
\ee
and
\bel{WP-2}
W^+(b)=W_i(b)~~~\text{for all}~~~b\in X'_i=X_i\setminus\{a_i\},~i=1,...,N.
\ee
\par One checks that $W^+$ is an admissible map, and $W^+(a^+)=x$. Thus , $L^+=(a^+,(X^+,E^+), W^+)$ belongs to $\Ac(x)$. Consequently, by Lemma \reff{S-T}, there exists a $d_{T^+,W^+}$-Lipschitz selection $f^+$ of $F\circ W^+$ with $d_{T^+,W^+}$-Lipschitz  seminorm $\le \gamma_0$.
\par The map
$$
f_i(b)=\left \{
\begin{array}{ll}
f^+(b),& \text{if}~~b\in X_i\setminus\{a_i\},\vspace{2mm}\\
f^+(a^+),& \text{if}~~b=a_i,
\end{array}
\right.
$$
is a $d_{T_i,W_i}$-Lipschitz selection of $F\circ W_i$ with $d_{T_i,W_i}$-Lipschitz seminorm $\le\gamma_0$, therefore $$f^+(a^+)\in O(x;L_i)~~~\text{for each}~~~i=1,...,N.$$
\par Thus, \rf{O-LI2} holds, completing the proof of Lemma \reff{L1}.\bx
\medskip
\par We are in a position to prove inequality \rf{GX-HD}.
\begin{lemma}\lbl{HD-Z} For every $x,y\in\MS$ the following inequality
$$
\dhf(G(x),G(y))\le \gamma_0\,\dm(x,y)
$$
holds.
\end{lemma}
\par {\it Proof.} We may suppose $x\ne y$, else the desired conclusion is obvious. Let us prove that
\bel{R-IN}
I=G(x)+\gamma_0\,\dm(x,y)\,\BY\supset G(y)\,.
\ee
Recall that by $\BY=\BL(0,1)$ we denote the unit ball in $\BS$.
\par If we can prove that, then by interchanging the roles of $x$ and $y$ we obtain also
$$
G(y)+\gamma_0\,\dm(x,y)\,\BY\supset G(x)\,.
$$
These two inclusions tell us that $\dhf(G(x),G(y))\le \gamma_0\,\dm(x,y)$, proving the lemma.\smallskip
\par Let us prove \rf{R-IN}. By definition,
$$
I=\left[\,\bigcap_{L\,\in\,\Ac(x)}\,O(x;L)\right]+
\gamma_0\,\dm(x,y)\,\BY\,.
$$
See \rf{GX-D}. We will check that
\bel{Z2}
\left[\,\bigcap_{L\,\in\,\Ac(x)}\,O(x;L)\right]+
\gamma_0\,\dm(x,y)\,\BY=\bigcap\,\left\{\,
\left[O(x;L_1)\cap...\cap O(x;L_N)\right]+\gamma_0\,\dm(x,y)\,\BY\right\},
\ee
where the first intersection of the right-hand side is taken over all finite sequences $L_1,...,L_N$ of elements of $\Ac(x)$.
\par Indeed, the left-hand side of \rf{Z2} is obviously contained if the right-hand side. Conversely, let $\xi$ belong to the right-hand side of \rf{Z2}. Then any finite subcollection of the compact sets
$$
K_L=\{\eta\in\BY: \xi-\gamma_0\,\dm(x,y)\,\eta\in O(x;L)\}
$$
has nonempty intersection. (The above sets are compact because $O(x;L)$ is compact.)
\par Therefore,
$$
\bigcap_{L\in\Ac(x)}\,K_L\ne\emp,
$$
proving that $\xi$ belongs to the left-hand side of \rf{Z2}. The proof of \rf{Z2} is complete.\smallskip
\par Thanks to \rf{Z2}, our desired inclusion \rf{R-IN} will follow if we can show that
\bel{Z3}
\left[O(x;L_1)\cap...\cap O(x;L_N)\right]+\gamma_0\,\dm(x,y)\,\BY\supset G(y)
\ee
for any $L_1,...,L_N\in\Ac(x)$. Then the proof of Lemma \reff{HD-Z} is reduced to the task of proving \rf{Z3}.
\smallskip
\par Let $L_i=(a_i,(X_i,E_i),W_i)$, and let $T_i=(X_i,E_i)$. Then $a_i$ is a node of the tree $T_i$, $i=1,...,N$. We introduce a new node $a^+$ and form the tree $T^+=(X^+,E^+)$ as in the proof of Lemma \reff{L1}. Thus $T^+$ arises by {\it gluing together the trees $T_i$ by identifying the $a_i$}.
\par We also introduce an admissible map $W^+:X^+\to\Mc$ as in the proof of Lemma \reff{L1}, see \rf{WP-1} and \rf{WP-2}.
\smallskip
\par We now introduce a new node $\ta$ not present in $T^+$. We define a new tree $\tT=(\tX,\tE)$ as follows.
\par \textbullet~ The nodes $\tX$ are the nodes in $X^+$, together with the new node $\ta$.
\smallskip
\par \textbullet~ The edges $\tE$ are the edges in $E^+$, together with a single edge joining $\ta$ to $a^+$.
\medskip
\par We define a map $\tW:\tT\to\Mc$ by setting
$$
\tW=W^+~~~\text{on}~~~T^+,~~\tW(\ta)=y.
$$
Then one checks that $\tT=(\tX,\tE)$ is a tree, $\tW:\tT\to\Mc$ is an admissible map, and $\tW(\ta)=y$.
\par Let $\tL=(\ta,(\tX,\tE),\tW)$. Recall that $G(y)\ne\emp$ by Lemma \reff{L1}.
Let $\eta\in G(y)$. Then, by definition, $\eta\in O(y;\tL)$ so that there exists a $d_{\tT,\tW}$-Lipschitz selection $\tf$ of $F\circ \tW$, with $d_{\tT,\tW}$-Lipschitz seminorm $\le\gamma_0$, and satisfying $\tf(\ta)=\eta$. See \rf{OL-X} and Definition \reff{OR-D}.
\par Restricting this $\tf$ to $T^+$ and arguing as in the proof of Lemma \reff{L1}, we see that                                                      
$$
\tf(a^+)\in
O(x;L_1)\cap...\cap O(x;L_N).
$$
\par On the other hand, our Lipschitz bound for $\tf$ gives
$$
\|\tf(a^+)-\eta\|=\|\tf(a^+)-\tf(\ta)\|
\le\gamma_0\,\rho(\tW(a^+),\tW(\ta))=\gamma_0\,\rho(x,y).
$$
Then,
$$
\eta\in \left[O(x;L_1)\cap...\cap O(x;L_N)\right]+\gamma_0\,\dm(x,y)\,\BY
$$
proving \rf{Z3}.\bx
\medskip
\par The proof of Theorem \reff{HDS-M} is complete.\bx
\smallskip
\par Theorem \reff{HDS-M} and Theorem \reff{ST-P} imply the following result.
\begin{corollary}\lbl{LS-MSP} Let $(\Mc,\dm)$ be a metric space. Let $\lambda$ be a positive constant and let $F:\MS\to\KM$ be a set-valued mapping.\smallskip
\par Suppose that for every subset $\Mc'\subset\MS$ consisting of at most $\ks$ points, the restriction $F|_{\Mc'}$ of $F$ to $\Mc'$ has a Lipschitz selection $f_{\Mc'}:\Mc'\to\BS$ whose seminorm satisfies $\|f_{\Mc'}\|_{\Lip(\Mc',\BS)}\le \lambda$. Then $F$ has a Lipschitz selection $f:\Mc\to\BS$ with $\|f\|_{\Lip(\Mc,\BS)}\le \gamma_2\lambda$.
\par Here $\gamma_2$ is a constant depending only on $m$.
\end{corollary}
\par {\it Proof.} We follow the scheme of the proof suggested in the Introduction. Let $\tdm=\lambda\dm$. Then the metric space $(\Mc,\tdm)$ and the set-valued mapping $F$ satisfy the hypothesis of Theorem \reff{HDS-M}. By this theorem,  there exists a mapping $G:\Mc\to\KMY$ such that $G(x)\subset F(x)$, $x\in\Mc$, and
\bel{GXYG}
\dhf(G(x),G(y))\le \gamma_0\,\tdm(x,y)=\gamma_0\,\lambda\,\dm(x,y)~~~\text{for all}~~x,y\in\Mc.
\ee
\par Let $f(x)=S_Y(G(x))$, $x\in \Mc$, where $S_Y:\Kc(Y)\to Y$ is the Steiner-type point operator from Theorem \reff{ST-P}. Clearly, by part (i) of Theorem \reff{ST-P},
$f(x)=S_Y(G(x))\in G(x)\subset F(x)$, i.e., $f$ is a {\it selection} of $F$ on $\Mc$.
\par By \rf{GXYG} and by  part (ii) of Theorem \reff{ST-P}, for every $x,y\in\MS$
$$
\|f(x)-f(y)\|=\|S_{\BS}(G(x))-S_{\BS}(G(y))\|\le \gamma_1\dhf(G(x),G(y))\le
\gamma_0\,\gamma_1\,\lambda\,\dm(x,y)=\gamma_2\,
\lambda_,\dm(x,y)
$$
where $\gamma_2=\gamma_0\,\gamma_1$.
\par Note that, by Theorem \reff{ST-P},  $\gamma_1=\gamma_1(\dim G(x),\dim G(y))$. Since $\dim G(x),\dim G(y)\le m$, and $\gamma_0$ depends only on $m$,
the constant $\gamma_2$ depends only on $m$ as well.  Thus $\|f\|_{\Lip(\Mc,\BS)}\le \gamma_2\lambda$, and the proof of the corollary is complete.\bx
\bigskip\medskip
\SECT{6. Pseudometric spaces: the final step of the proof of the finiteness principle.}{6}
\addtocontents{toc}{6. Pseudometric spaces: the final step of the proof of the finiteness principle. \hfill\thepage\par\VST}

\indent\par In this section we prove Theorem \reff{MAIN-FP}, the Finiteness Principle for Lipschitz Selections, and Theorem \reff{FINITE}, a variant of Theorem \reff{MAIN-FP} for finite pseudometric spaces.
\par Until the proof of Theorem \reff{MAIN-FP} given at the end of Section 6 we assume that $(\Mc,\rho)$ is a pseudometric space satisfying the following condition:
\bel{FIN-A}
\rho(x,y)<\infty~~~\text{for all}~~x,y\in\Mc\,.
\ee
\par Until the end of Section 6 we write $\gamma_2$ to denote the constant $\gamma_2(m)$ from Corollary \reff{LS-MSP}. We also recall that $\ks$ is the constant defined by \rf{DK-1} and \rf{DK-2}.

\indent\par {\bf 6.1. Set-valued mappings with compact images on pseudometric spaces.}
\addtocontents{toc}{~~~~6.1. Set-valued mappings with compact images on pseudometric spaces.\hfill \thepage\par}
\smallskip
\par In this section we prove an analog of Corollary \reff{LS-MSP} for pseudometric spaces.
\begin{proposition}\lbl{PM-SP} Let $(\Mc,\dm)$ be a pseudometric space satisfying \rf{FIN-A}, and let $\lambda>0$. Let $F:\MS\to\KM$ be a set-valued mapping such that for every subset $\Mc'\subset\MS$ consisting of at most $\ks$ points, the restriction $F|_{\Mc'}$ of $F$ to $\Mc'$ has a Lipschitz selection $f_{\Mc'}:\Mc'\to\BS$ with  $\|f_{\Mc'}\|_{\Lip(\Mc',\BS)}\le \lambda$.
\par Then $F$ has a Lipschitz selection $f:\Mc\to\BS$ with $\|f\|_{\Lip(\Mc,\BS)}\le \gamma_2\lambda$.
\end{proposition}
\par {\it Proof.} A selection of $F$ may be regarded as a point of the Cartesian product
$$
\Fc=\prod_{x\in\Mc}\,F(x)\,.
$$
We endow $\Fc$ with the product topology. Then $\Fc$ is compact because each $F(x)$ is compact.
\par For $\ve>0$ and $x,y\in\Mc$, let
$$
\rho_\ve(x,y)=\left \{
\begin{array}{ll}
\rho(x,y)+\ve,& \text{if}~~x\ne y,\vspace*{2mm}\\
0,& \text{if}~~x=y.
\end{array}
\right.
$$
\par Then $(\Mc,\rho_\ve)$ is a metric space. For any $\Mc'\subset\Mc$ with $\#\Mc'\le \ks$ there exists a selection of $F|_{\Mc'}$ with $\rho$-Lipschitz seminorm $\le \lambda$, hence with $\rho_\ve$-Lipschitz seminorm $\le \lambda$. By Corollary \reff{LS-MSP}, $F$ has a selection with $\rho_\ve$-Lipschitz seminorm $\le\gamma_2\lambda$.
\par Let $\Sel(\ve)$ be the set of all selections of $F$ with $\rho_\ve$-Lipschitz seminorm at most $\gamma_2\lambda$. Then $\Sel(\ve)$ is a closed subset of $\Fc$. We have just seen that $\Sel(\ve)$ is non-empty. Because
$$
\Sel(\ve)\subset \Sel(\ve')~~~\text{for}~~~\ve<\ve',
$$
it follows that
$$
\Sel(\ve_1)\cap\Sel(\ve_2)\cap...\cap\Sel(\ve_N)\ne\emp
$$
for any $\ve_1,\ve_2,...,\ve_N>0$.
\par Because $\Fc$ is compact and each $\Sel(\ve)$ is closed in $\Fc$, it follows that
$$
\bigcap_{\ve>0}\Sel(\ve)\ne\emp\,.
$$
However, any $f\in\capbig\{\Sel(\ve):\ve>0\}$ is a selection of $F$ with $\rho$-Lipschitz seminorm $\le\gamma_2\lambda$. 
\par The proof of Proposition \reff{PM-SP} is complete.\bx
\medskip

\indent\par {\bf 6.2. Finite pseudometric spaces.}
\addtocontents{toc}{~~~~6.2. Finite pseudometric spaces.\hfill \thepage\par}
\smallskip
\par In this section we prove an analog of Proposition \reff{PM-SP} for a {\it finite} pseudometric space $(\Mc,\rho)$ and a set-valued mapping $F:\Mc\to\CNMY$. See Proposition \reff{FN-PM} below. Our proof of this proposition relies on three auxiliary lemmas.
\begin{lemma}\lbl{FM-B} Let $\lambda>0$ and let $(\Mc,\dm)$ be a finite metric space. Let $F$ be a set-valued mapping on $\Mc$ which to every $x\in\Mc$ assigns a non-empty convex bounded subset of $Y$ of dimension at most $m$.
\par Suppose that for every subset $\Mc'\subset\MS$ with $\#\Mc'\le \ks$, the restriction $F|_{\Mc'}$ of $F$ to $\Mc'$ has a Lipschitz selection $f_{\Mc'}:\Mc'\to\BS$ with  $\|f_{\Mc'}\|_{\Lip(\Mc',\BS)}\le \lambda$.
\par Then $F$ has a Lipschitz selection $f:\Mc\to\BS$ with $\|f\|_{\Lip(\Mc,\BS)}\le 2\gamma_2\lambda$.
\end{lemma}
\par {\it Proof.} We introduce a new set-valued mapping on $\Mc$ defined by
$$
\tF(x)=(F(x))^{\cl}~~~~\text{for all}~~~x\in\Mc.
$$
Here the sign $\cl$ denotes the closure of a set in $Y$.
\par Since the sets $F(x)$, $x\in\Mc$, are finite dimensional and bounded, each set $\tF(x)$ is compact so that $\tF:\Mc\to\KM$. Furthermore, since $F(x)\subset\tF(x)$ on $\Mc$, the mapping $\tF$ satisfies the hypothesis of Proposition \reff{PM-SP}.
\par By this proposition, there exists a mapping $\tf:\Mc\to Y$ such that
\bel{TF-CL}
\tf(x)\in\tF(x)=(F(x))^{\cl}~~~\text{for all}~~x\in\Mc,
\ee
and
\bel{FW-11}
\|\tf(x)-\tf(y)\|\le \gamma_2\,\lambda\,\rho(x,y) ~~~\text{for all}~~x,y\in\Mc.
\ee
\par Since $\Mc$ is a {\it finite} metric space, the following quantity
\bel{DL-M}
\delta=\gamma_2\,\lambda\,\min_{x,y\in\Mc,\,x\ne y}\rho(x,y)
\ee
is positive. Therefore, by \rf{TF-CL}, for each $x\in\Mc$ there exists a point $f(x)\in F(x)$ such that
$$
\|f(x)-\tf(x)\|\le \delta/2\,.
$$
\par Thus $f:\Mc\to Y$ is a selection of $F$ on $\Mc$. Let us estimate its Lipschitz seminorm. For every $x,y\in\Mc$ (distinct), by \rf{FW-11} and \rf{DL-M},
$$
\|f(x)-f(y)\|\le \|f(x)-\tf(x)\|+\|\tf(x)-\tf(y)\|+\|\tf(y)-f(y)\|
\le \delta/2+
\gamma_2\lambda\,\rho(x,y)+\delta/2\le 2\gamma_2\lambda\,\rho(x,y).
$$
Hence, $\|f\|_{\Lip(\Mc,\BS)}\le 2\gamma_2\lambda$, and the proof of the lemma is complete.\bx
\medskip

\par The second auxiliary lemma provides additional properties of sets $\Gamma_\ell$ defined in Section 3.1 (see \rf{BC-3} and Definition \reff{GML-D}).  We will need these properties in the proof of Lemma \reff{P-2} below.
\begin{lemma}\lbl{G-AP} Let $(\Mc,\rho)$ be a finite pseudometric space satisfying \rf{FIN-A}. Let $\ell\ge 0$ and let $F:\Mc\to\CNMY$. Suppose that for every subset $\Mc'\subset\Mc$ with $\#\Mc'\le k_{\ell+1}$ the restriction $F|_{\Mc'}$ of $F$ to $\Mc'$ has a Lipschitz selection $f_{\Mc'}:\Mc'\to Y$ with  $\|f_{\Mc'}\|_{\Lip(\Mc',Y)}\le\lambda$.
\par Let $x_0\in\Mc$, $\xi_0\in\GL(x_0)$, and let $1\le k\le \ell+1$. Let $S$ be a subset of $\Mc$ with $\#S=k$ containing $x_0$.
\par Then there exists a mapping $f^S:S\to Y$ such that
\smallskip
\par (a) $f^S(x_0)=\xi_0$.
\smallskip
\par (b) $f^S(y)\in \Gamma_{\ell+1-k}(y)$ for all $y\in S$.
\smallskip
\par (c) $\|f^S\|_{\Lip(S,Y)}\le 3^k\lambda$.
\end{lemma}
\par {\it Proof.} We recall that the sequence of positive integers $k_\ell$ is defined by the formula \rf{KEL}.
\par We proceed by induction on $k$. For $k=1$, we have $S=\{x_0\}$, and we can just set $f^S(x_0)=\xi_0$.
\par For the induction step, we fix $k\ge 2$ and suppose the lemma holds for $k-1$; we then prove it for $k$. Thus, let $\xi_0\in\GL(x_0)$, $x_0\in S$, $\# S=k\le \ell+1$.
\par Set $\hS=S\setminus\{x_0\}$. We pick $\hx_0\in\hS$ to minimize $\rho(\hx_0,x_0)$, and we pick $\hxi_0\in\GLO(\hx_0)$ such that $\|\hxi_0-\xi_0\|\le \lambda\,\rho(\hx_0,x_0)$. (See Lemma \reff{G-AB} (b).) For $y\in\hS$ we have $\rho(y,x_0)\ge \rho(\hx_0,x_0)$, hence
\bel{R-V}
\rho(y,\hx_0)+\rho(\hx_0,x_0)\le[\rho(y,x_0)
+\rho(x_0,\hx_0)]+\rho(\hx_0,x_0)\le 3\rho(y,x_0).
\ee
\par By the induction hypothesis, there exists
$\hf:\hS\to Y$ such that
\smallskip
\par $(\ha)$ $\hf(\hx_0)=\hxi_0$.
\smallskip
\par $(\hb)$ $\hf(y)\in \Gamma_{(\ell-1)+1-(k-1)}(y)=\Gamma_{\ell+1-k}(y)$ for all $y\in \hS$.
\smallskip
\par $(\hc)$ $\|\hf\|_{\Lip(\hS,Y)}\le 3^{k-1}\lambda$.
\smallskip
\par We now define $f:S\to Y$ by setting
$$
f(y)=\hf(y)~~~\text{for}~~~y\in\hS;~~~~f(x_0)=\xi_0.
$$
\par Then $f$ obviously satisfies $(a)$ and $(b)$. To see that $f$ satisfies (c), we first recall $(\hc)$; thus it is enough to check that
$$
\|f(y)-f(x_0)\|\le 3^k\lambda\,\rho(y,x_0)
$$
for $y\in\hS$, i.e.,
$$
\|\hf(y)-\xi_0\|\le 3^k\lambda\,\rho(y,x_0)~~~\text{for}~~~y\in\hS.
$$
\par However, for $y\in\hS$ we have
$$
\|\hf(y)-\xi_0\|\le \|\hf(y)-\hxi_0\|+\|\hxi_0-\xi_0\|=
\|\hf(y)-\hf(\hx_0)\|+\|\hxi_0-\xi_0\|\le
3^{k-1}\lambda\,\rho(y,\hx_0)+\lambda\,\rho(\hx_0,x_0)
$$
thanks to $(\hc)$, and the definition of $\hxi_0$.
\par Therefore,
$$
\|\hf(y)-\xi_0\|\le
3^{k-1}\lambda\,[\rho(y,\hx_0)+\rho(\hx_0,x_0)]\le 3^{k}\lambda\,\rho(y,x_0),
$$
by \rf{R-V}.
\par Thus, $f$ satisfies $(a)$, $(b)$, $(c)$, completing our induction.\bx
\medskip
\par We turn to the last auxiliary lemma. Let
\bel{KNN}
\tl=\ks~~~~\text{and let}~~~~\kn=k_{\tl+1}
\ee
where $k_\ell=(m+2)^\ell$, see \rf{KEL}.
\begin{lemma}\lbl{P-2} Let $(\Mc,\dm)$ be a finite pseudometric space satisfying \rf{FIN-A}, $x_0\in\Mc$ and $\lambda>0$.
\par Let $F:\MS\to\CNMY$ be a set-valued mapping such that for every subset $\Mc'\subset\MS$ consisting of at most $\kn$ points, the restriction $F|_{\Mc'}$ of $F$ to $\Mc'$ has a Lipschitz selection $f_{\Mc'}:\Mc'\to\BS$ with  $\|f_{\Mc'}\|_{\Lip(\Mc',\BS)}\le \lambda$.
\par Then there exists a point $\xi_0\in F(x_0)$ such that the following statement holds: For every subset $S\subset\Mc$ with $\#S\le \ks$, there exists a mapping $f_S:S\to Y$ with $\|f_S\|_{\Lip(S,Y)}\le C\lambda$ such that
\bel{KN-F}
\|f_S(x)-\xi_0\|\le C\lambda\,\rho(x,x_0) ~~~\text{for every}~~~x\in S,
\ee
and
\bel{FNB-4}
f_S(x)\in F(y)+\lambda\,\rho(x,y)\,B_Y~~~\text{for every}~~~x\in S, y\in\Mc\,.
\ee
Here $C$ is a constant depending only on $m$.
\end{lemma}
\par {\it Proof.} By the lemma's hypothesis, \rf{KNN} and
by Lemma \reff{G-AB} (a),
$$
\Gamma_{\tl}(x)\ne\emp ~~~\text{for every}~~~x\in\Mc\,.
$$
(See also Remark \reff{GL-DFL}.)
\par Let $\xi_0\in \Gamma_{\tl}(x_0)$. By \rf{GL-F2},
$$
\xi_0\in \Gamma_{\tl}(x_0)\subset F(x_0).
$$
\par Let $S\subset\Mc$, $\#S\le \ks$. Let $\tS=S\cup\{x_0\}$ and let $k=\#\tS=\#(S\cup\{x_0\})\,.$
Then
$$
1\le k\le\#S+1\le \ks+1=\tl+1.
$$
Therefore, by Lemma \reff{G-AP}, there exists a mapping $f^{\tS}:\tS\to Y$ with $\|f^{\tS}\|_{\Lip(\tS,Y)}\le 3^k\lambda$ such that $f^{\tS}(x_0)=\xi_0$ and
$$
f^{\tS}(x)\in \Gamma_{\tl+1-k}(x)~~~\text{for all}~~~x\in\tS.
$$
Recall that $k\le\tl+1=\ks+1$ so that
$$
\|f^{\tS}\|_{\Lip(\tS,Y)}\le C\lambda
$$
with $C=3^{\ks+1}$. Since $\ks$ depends only on $m$, the constant $C$ depends only on $m$ as well.
\par Hence, by \rf{BC-5},
\bel{FG-12}
f^{\tS}(x)\in \Gamma_{\tl+1-k}(x)\subset \Gamma_{0}(x)~~~\text{for every}~~~x\in\tS.
\ee
\par Let
$$
f_S=f^{\tS}|_S\,.
$$
\par Then $\|f_S\|_{\Lip(S,Y)}\le \|f^{\tS}\|_{\Lip(\tS,Y)}\le C\lambda$. Moreover, by \rf{FG-12},
\bel{G-27}
f_S(x)\in \Gamma_{0}(x)~~~\text{for all}~~~x\in S.
\ee
\par Since $\|f^{\tS}\|_{\Lip(\tS,Y)}\le C\lambda$ and $x_0\in\tS$,
$$
\|f_S(x)-\xi_0\|=\|f^{\tS}(x)-f^{\tS}(x_0)\|\le C\lambda\,\rho(x,x_0) ~~~\text{for every}~~~x\in S.
$$
Furthermore, by \rf{GM-ZR} and \rf{G-27}, for every
$x\in S$
$$
f_S(x)\in \Gamma_{0}(x)=\bigcap_{y\in\Mc}
\left(F(y)+\lambda\,\rho(x,y)\,B_Y\right)
$$
so that
$$
f_S(x)\in F(y)+\lambda\,\rho(x,y)\,B_Y~~~\text{for every}~~~x\in S, y\in\Mc\,.
$$
\par The proof of the lemma is complete.\bx
\begin{proposition}\lbl{FN-PM} Let $(\Mc,\rho)$ be a finite pseudometric space satisfying \rf{FIN-A}, and let $\lambda>0$. Let $F:\MS\to\CNMY$ be a set-valued mapping such that for every subset $\Mc'\subset\MS$ with $\#\Mc'\le \kn$, the restriction $F|_{\Mc'}$ of $F$ to $\Mc'$ has a Lipschitz selection $f_{\Mc'}:\Mc'\to\BS$ with  $\|f_{\Mc'}\|_{\Lip(\Mc',\BS)}\le \lambda$.
\par Then $F$ has a Lipschitz selection $f:\Mc\to\BS$ with $\|f\|_{\Lip(\Mc,\BS)}\le \gamma_3\lambda$ where $\gamma_3$ is a constant depending only on $m$.
\end{proposition}
\par {\it Proof.} Let $x_0\in\Mc$. By Lemma \reff{P-2}, there exists a point $\xi_0\in F(x_0)$ such that for every
set $S\subset\Mc$ with $\#S\le\ks$ there exists a mapping $f_S:S\to Y$ with $\|f_S\|_{\Lip(S,Y)}\le C\lambda$ such that \rf{KN-F} and \rf{FNB-4} hold. Here $C$ is a constant depending only on $m$.
\par We introduce a new set-valued mapping $\tF:\Mc\to \CNMY$ by letting
\bel{TF-W}
\tF(x)=\left(\,\bigcap_{y\in\Mc}
\left[F(y)+\lambda\,\rho(x,y)\,B_Y\right]\right)
\bigcap\,B_Y(\xi_0,C\lambda\,\rho(x,x_0)),
~~~x\in \Mc\,.
\ee
In particular, taking $y=x$ in the above formula we obtain that
\bel{IM-FX}
\tF(x)\subset F(x)
~~~\text{for all}~~~x\in \Mc\,.
\ee
\par By Lemma \reff{P-2} and definition \rf{TF-W}, for every set $S\subset\Mc$ consisting of at most $\ks$ points the restriction $\tF|_{S}$ of $\tF$ to $S$ has a Lipschitz selection $f_{S}:S\to\BS$ with  $\|f_{S}\|_{\Lip(S,Y)}\le C\lambda$. In particular, $\tF(x)\ne\emp$ for every $x\in
\Mc$.
\smallskip
\par Let us introduce a binary relation ``$\sim$'' on $\MS$ by letting
$$
x\sim y~~~\Longleftrightarrow~~~\rho(x,y)=0\,.
$$
Clearly, ``$\sim$'' satisfies the axioms of an equivalence relation, i.e., it is reflexive, symmetric and transitive. Given $x\in\MS$, by $[x]=\{y\in\MS:~y\sim x\}$ we denote the equivalence class of $x$. Let
$$
[\Mc]=\Mc\,/\sim\,\,=\,\{\,[x]: x\in\MS\,\}
$$
be the corresponding quotient set of $\MS$ by ``$\sim$'', i.e.,  the family of all possible equivalence classes of $\MS$ by ``$\sim$''. Finally, given an equivalence class $U\in[\Mc]$ let us choose a point $w_U\in U$ and put
$$
W=\{w_U:U\in[\Mc]\}.
$$
\par Clearly, $(W,\rho)$ is a {\it finite metric space}. Let
\bel{HF-7}
\hF=\tF|_W.
\ee
Then, by \rf{TF-W}, \rf{HF-7} and \rf{IM-FX}, $\hF$ is a set-valued mapping defined on a finite metric space which takes values in the family of all non-empty convex {\it bounded} subsets of $Y$ of dimension at most $m$. Furthermore, this mapping satisfies the hypothesis of Lemma \reff{FM-B} with $C\lambda$ in place of $\lambda$.
\par Therefore, by this lemma, there exists a Lipschitz selection $\hf:W\to Y$ of $\hF$ on $W$ with
$$
\|\hf\|_{\Lip(W,Y)}\le 2\gamma_2 \,C\lambda =\gamma_3\lambda.
$$
Here $\gamma_3=2\gamma_2 C$ is a constant depending only on $m$ (because $\gamma_2$ and $C$ depend on $m$ only).
\smallskip
\par We define a mapping $f:\Mc\to Y$ by letting
$$
f(x)=\hf(w_{[x]}), ~~~~~x\in\Mc.
$$
\par Then $f$ is a {\it selection of $F$} on $\Mc$. Indeed, let $x\in\Mc$. Since $\hf$ is a selection of $\hF=\tF|_W$ on $W$, and $w_{[x]}\in W$,
$$
f(x)=\hf(w_{[x]})\in \hF(w_{[x]})
$$
so that, by \rf{TF-W},
$$
f(x)\in \hF(w_{[x]})\subset F(x)+\lambda\,\rho(w_{[x]},x)\,B_Y.
$$
But $w_{[x]}\sim x$ so that $\rho(w_{[x]},x)=0$ proving that $f(x)\in F(x)$.
\par Let us prove that $\|f\|_{\Lip(\Mc,Y)}\le \gamma_3\lambda$, i.e.,
\bel{LN-F2}
\|f(x)-f(y)\|\le \gamma_3\lambda\,\rho(x,y)~~~\text{for all}~~~x,y\in\Mc.
\ee
\par In fact, since $\|\hf\|_{\Lip(W,Y)}\le \gamma_3\lambda$,
\be
\|f(x)-f(y)\|&=&\|\hf(w_{[x]})-\hf(w_{[y]})\|\le \gamma_3\lambda\,\rho(w_{[x]},w_{[y]})\nn\\
&\le& \gamma_3\lambda\,(\rho(w_{[x]},x)+\rho(x,y)+
\rho(y,w_{[y]}))=
\gamma_3\lambda\,\rho(x,y)\nn
\ee
proving \rf{LN-F2}.
\par The proof of Proposition \reff{FN-PM} is complete.\bx
\bigskip

\indent\par {\bf 6.3. The sharp finiteness number.}
\addtocontents{toc}{~~~~6.3. The sharp finiteness number.\hfill \thepage\par\VSU}
\medskip

\par In this section we prove Theorem \reff{D-2M}. We note that for set-valued mappings $F$ on $\Mc$ whose values are convex {\it compact} subsets of $Y$ of dimension at most $m$ (i.e., $F(u)\in\KMY$ for all $u\in\Mc$) the statement of Theorem \reff{D-2M} has been proved in \cite{S02}. Our proof below for the general case of mappings $F:\Mc\to\CNMY$ will follow the scheme suggested in \cite{S02}.
\par Let $(\Mc,\rho)$ be a finite pseudometric space with a finite pseudometric $\rho$. Let $T=(\Mc,E)$ be a finite tree whose set of nodes coincides with $\Mc$. Following the notation of Section 5, we write $x\je y$ to indicate that nodes $x,y\in\Mc$ are joined by an edge in $T$. We denote that edge by $[xy]$.
\par The tree $T$ gives rise a {\it tree pseudometric} $d_T:\Mc\times \Mc\to \R_+$ defined by
$$
d_{T}(x,y)=\rho(x,y)~~~\text{for every}~~ x,y\in \Mc,~ x\je y\,.
$$
We recall that for arbitrary $x,y\in\Mc$, $x\ne y$,
we define the distance $d_{T}(x,y)$ by
\bel{DST-T}
d_T(x,y)=\smed_{i=1}^L\,\, \rho(x_{i-1},x_i)
\ee
where $\{x_i:i=1,...,L\}$ is the one and only one ``path''
joining $x$ to $y$ in $T$, i.e., $x_i\in\Mc$ and
\bel{PATH}
x=x_0\je x_1\je...\je x_L=y~~~\text{with all the}~~x_i ~~\text{distinct}.
\ee
See \rf{D-TR}. We also set $d_T(x,y)=0$ for $x=y$. We refer to $(\Mc,d_T)$ as a {\it pseudometric tree} generated by $T$.
\par Clearly, by the triangle inequality,
$$
\rho(x,y)\le d_{T}(x,y)~~~
\text{for every}~~ x,y\in \Mc\,.
$$
\par Given a node $u\in \Mc$, by $J(u;T)$ we denote the family of its neighbors in $T$:
$$
J(u;T)=\{v\in \Mc: v\je u\}.
$$
We let $\deg_T u$ denote the number of neighbors of the node $u$; thus $\deg_T u=\# J(u;T)$.
\par For a number $a\in\R$ by $\lceil a\rceil$ we denote the integer $m$ such that $m-1<a\le m$.
\begin{proposition}\lbl{M-TR} Let $(\Mc,\dm)$ be a finite pseudometric space with a finite pseudometric $\rho$. There exists a tree $T=(\Mc,E)$ satisfying the following conditions:\smallskip
\par (i) For every $x,y\in\Mc$
\bel{T-RHO}
\rho(x,y)\le d_{T}(x,y)\le \theta\,\rho(x,y)\,.
\ee
Here $\theta=\theta(\#\Mc)\ge 1$ is a constant depending only on the cardinality of $\Mc$.
\smallskip
\par (ii) The following inequality
\bel{MAX}
\max_{x\in\Mc}\deg_T x\ge\,\, \lceil\,\log_2(\#\Mc)\rceil
\ee
holds.
\end{proposition}
\par {\it Proof.} We prove the proposition by induction on $k:=\#\Mc$. Clearly, the proposition is trivial for $k=1$. We suppose that the proposition holds for given $k\ge 1$ and prove it for $k+1$.
\par Let $(\Mc,\rho)$ be a pseudometric space
with $\#\Mc=k+1$. Choose points $x_0,y_0\in\Mc$ such that
$$
\diam(\Mc)=\max\limits_{x,y\in\Mc}\rho(x,y) =\rho(x_0,y_0).
$$
\par We define a partition $\{\Mc',\Mc''\}$ of $\Mc$ as follows. If $\diam \Mc=0$, we put $\Mc'=\Mc\setminus\{y_0\}$ and $\Mc''=\{y_0\}$.
\par Suppose that $\diam\Mc>0$. In this case we let $\Mc'$ denote a set of all $y\in\Mc$ satisfying the following condition: there exists a sequence of points $\{z_0=x_0,z_1,...,z_n=y\}$ in $\Mc$ with all the $z_i$ distinct, such that
\bel{A2.2}
\rho(z_i,z_{i+1})<\frac{1}{k}\diam(\Mc)~~~\text{for every}~~i=0,...,n-1.
\ee
\par Clearly, $\Mc'\ne\emp$ because it contains $x_0$.
Let us show that
\bel{MC-2}
\Mc''=\Mc\setminus\Mc'\ne\emp\,.
\ee
\par Indeed, the point $y_0\in\Mc''$, otherwise there exist elements $\{z_0=x_0,z_1,...,z_n=y_0\}$ with all the $z_i$ distinct, such that inequality \rf{A2.2} holds. Since $n\le k=\#\Mc-1$, we obtain the following
$$
\diam\Mc=\rho(x_0,y_0)\le\sum\limits_{i=0}^{n-1}
\rho(z_i,z_{i+1})< \sum\limits_{i=0}^{n-1}
\frac{1}{k}\diam\Mc=\frac{n}{k}\diam\Mc \le\diam\Mc.
$$
This contradiction proves \rf{MC-2}.
\smallskip
\par Let us prove that
\bel{A2.3}
\rho(x',x'')\ge\frac{1}{k}\diam\Mc~~~\text{for all}~~ x'\in\Mc'~~\text{and}~~x''\in\Mc''\,.
\ee
\par Clearly, this inequality is trivial if $\diam \Mc=0$. Let $\diam\Mc>0$. Suppose that there exist $x'\in\Mc'$ and $x''\in\Mc''$ such that $\rho(x',x'')<\frac{1}{k}\diam\Mc$. By definition of $\Mc'$, there exists a path
$\{z_0=x_0,z_1,...,z_n=x'\}$ with all the $z_i$ distinct satisfying  inequality \rf{A2.2}. Clearly, $z_i\in\Mc'$ so that $x''\ne z_i$ for every $i=0,...,n$. Then the path
$\{z_0=x_0,z_1,...,z_n=x',z_{n+1}=x''\}$ satisfies \rf{A2.2} so that $x''\in\Mc'$. This contradiction implies \rf{A2.3}.
\smallskip
\par We turn to construction of a tree $T=(\Mc,E)$ satisfying inequalities \rf{T-RHO} and \rf{MAX}.
\par We will need only the following properties of the sets $\Mc'$ and $\Mc''$: (i) $\Mc',\Mc''\ne\emp$, (ii) $\Mc'\cup\Mc''=\Mc$, (iii) inequality \rf{A2.3} holds. This enables us, without loss of generality, to assume that
$\#\Mc'\ge\#\Mc''$. Hence,
\bel{K-M2}
k+1=\#\Mc\le 2\, \#\Mc'.
\ee
Since $\#\Mc'\le k$, by the induction assumption there exist a tree $T'=(\Mc',E')$ and a node $a'\in\Mc'$ such that
\bel{A2.6}
d_{T'}(x',y')\le \theta(k)\, \rho(x',y')~~~\text{for all}~~x',y'\in\Mc',
\ee
and
$$
\deg_{T'} a'\ge \,\lceil\,\log_2(\#\Mc')\rceil\,.
$$
By this inequality and \rf{K-M2},
\bel{DG-TT}
\deg_{T'} a'\ge \,\,\lceil\,\log_2(\#\Mc')\rceil\,\,
\ge\,\, \lceil\,\log_2(\#\Mc)\rceil\,-1.
\ee
\par Since $\#\Mc''\le k$, by the induction assumption there exists a tree $T''=(\Mc'',E'')$ such that
\bel{A-TPP}
d_{T''}(x'',y'')\le \theta(k)\, \rho(x'',y'')
~~~\text{for every}~~x'',y''\in\Mc''\,.
\ee
\par We form a tree $T=(\Mc,E)$ as follows.
We fix an arbitrary point $a_0\in \Mc''$ and define the family $E$ of edges of $T$ as the union of the families $E'$ and $E''$ together with an edge joining $a'$ with $a_0$. Thus,
\bel{E-T}
E=E'\cup E''\cup \{[a'a_0]\}\,.
\ee
\par Clearly, $(\Mc',d_{T'})$ and $(\Mc'',d_{T''})$ are pseudometric subspaces of $(\Mc,d_{T})$, i.e.,
\bel{TPP-T}
d_{T'}(x',y')=d_T(x',y'),~
d_{T''}(x'',y'')=d_T(x'',y'')~~\text{provided}~~x',y'\in\Mc',
~x'',y''\in\Mc''.
\ee
\par By \rf{DG-TT},
$$
\deg_T a'=\deg_{T'} a'+1\ge\,\,\lceil\,\log_2(\#\Mc)\rceil
$$
proving \rf{MAX}.
\par Let us prove \rf{T-RHO}. By \rf{A2.6}, \rf{A-TPP} and \rf{TPP-T}, it suffices to prove this inequality for every $x\in \Mc'$ and every $y\in\Mc''$.  By \rf{E-T} and definition \rf{DST-T},
$$
d_T(x,y)= d_{T'}(x,a')+\rho(a',a_0)+d_{T''}(a_0,y)\,.
$$
Hence, by \rf{A2.6}, \rf{A-TPP} and \rf{A2.3},
$$
d_T(x,y)\le \theta(k)\rho(x,a')+\rho(a',y_0)+\theta(k)\rho(a_0,y)\le
(2\theta(k)+1)\diam\Mc\le(2\theta(k)+1)\,k\,\rho(x,y)
$$
proving \rf{T-RHO} with $\theta(k+1)=k(2\theta(k)+1)$.\bx
\bigskip
\par {\it Proof of Theorem \reff{D-2M}.} We prove the theorem by induction on $k:=\#\Mc$.
\par For $k\le 2^{\min(m+1,\dim Y)}$ there is nothing to prove. We suppose that this result is true for given
$k\ge 2^{\min(m+1,\dim Y)}$, and prove it for $k+1$.
\par Let $(\Mc,\rho)$ be a pseudometric space
with $\#\Mc=k+1$, and let $F:\Mc\to\CNMY$ be a set-valued mapping satisfying  the hypotheses of Theorem \reff{D-2M}.
Then, by the induction assumption, for every subset $\Mc'\subset\Mc$ with $\#\Mc'\le k$ the restriction $F|_{\Mc'}$ has a Lipschitz selection $f_{\Mc'}:\Mc'\to Y$ such that $\|f_{\Mc'}\|_{\Lip(\Mc',Y)}\le\gamma(k)$.
\par Our aim is to prove the existence of a mapping $f:\Mc\to Y$ such that
$$
f(x)\in F(x)~~~\text{for every}~~~x\in\Mc,
$$
and
$$
\|f(x)-f(y)\|\le \gamma(k+1)\,\rho(x,y)
~~~\text{for all}~~~x,y\in\Mc\,.
$$
\par By Proposition \reff{M-TR} there exists a tree $T=(\Mc,E)$ satisfying conditions (i) and (ii) of the proposition. Thus,
\bel{M-TH}
\rho(x,y)\le d_{T}(x,y)\le \theta\,\rho(x,y)~~~\text{for all}~~~x,y\in\Mc
\ee
where $\theta=\theta(k+1)$ is a constant depending only on $k$. Furthermore, there exists a node $x_0\in\Mc$ such that $\deg_T x_0\ge\,\, \lceil\,\log_2(\#\Mc)\rceil$. Since
$\#\Mc=k+1>2^{\min(m+1,\dim Y)}$, we obtain the following inequality:
\bel{DG}
\deg_T x_0\ge\,\min(m+2,\dim Y+1)\,.
\ee
\par We recall that $J(x_0;T)=\{u\in\Mc: u \je x_0\}$
denotes the family of neighbors of $x_0$ in $T$. Therefore, by \rf{DG},
\bel{J-M2}
\#J(x_0;T)\ge \min(m+2,\dim Y+1)\,.
\ee
\par Given $u\in J(x_0;T)$ we let $\BR(u)$ denote a subset of $\Mc$ defined by
\bel{T-B1}
\BR(u)=\{x_0\}\cup\{u'\in \Mc: \text{the unique path joining}~~u'~~\text{to}~~x_0~~\text{in}~~
X~~\text{includes}~~u\}\,.
\ee
See \rf{PATH}. We refer to $\BR(u)$ as an $u$-branch of the node $x_0$ in the tree $T$.
\par Let us note two obvious properties of branches:
\smallskip
\par $(\bullet 1)$~ The family of subsets
$\{\BR(u)\setminus \{x_0\}:u\in J(x_0;T)\}$ and the singleton $\{x_0\}$ form a partition of $\Mc$;
\smallskip
\par $(\bullet 2)$~  Let  $u,v\in J(x_0;T)$, $u\ne v$, and let $a\in \BR(u)$, $b\in \BR(v)$. Then
$$
d_T(a,b)=d_T(a,x_0)+d_T(x_0,b).
$$
\par We introduce a new set-valued mapping $\tF:\Mc\to\CNMY$ as follows: we put
$$
\tF(x)=F(x)~~~\text{for every}~~x\in\Mc,~x\ne x_0,
$$
and
\bel{F-W1-N}
\tF(x_0)=\left \{
\begin{array}{ll}
F(x_0),& \text{if}~~m<\dim Y,\vspace*{2mm}\\
Y,& \text{if}~~m=\dim Y.
\end{array}
\right.
\ee
\par Given $u\in J(x_0;T)$ we let $\Or(u)$ denote a subset of $Y$ defined by
\bel{OR-U}
\Or(u)=\left\{g(x_0): g~\text{is a}~\rho\text{-Lipschitz selection of}~~ \tF|_{\BR(u)}~\text{with}~ \|g\|_{\Lip((\BR(u),\rho),Y)}\le 2\gamma(k)\theta\,\right\}.
\ee
\par Let us prove that
\bel{M-INT}
F(x_0)\bigcap\left\{
\bigcap_{u\in J(x_0;T)}\,\Or(u)
\right\}\ne\emp\,.
\ee
\smallskip
\par Consider two cases.
\smallskip
\par {\it The first case:}
\bel{C-1}
m<\dim Y\,.
\ee
\par Clearly, since $x_0\in\BR(u)$ and $\tF=F$ on $\Mc$, the set $F(x_0)\supset\Or(u)$ for each $u\in J(x_0;T)$. Therefore it suffices to prove that
\bel{ML-Y}
\bigcap_{u\in J(x_0;T)}\,\Or(u)\ne\emp\,.
\ee
\par It is also clear that $\{\Or(u):u\in J(x_0;T)\}$ is a finite family of convex sets lying in the affine hull of $F(x_0)$, whose dimension is bounded by $m$. Therefore, by Helly's Theorem \reff{HT-2}, \rf{ML-Y} holds provided
\bel{OR-3}
\bigcap_{i=1}^{m+1}\,\Or(u_i)\ne\emp
\ee
for any $m+1$ nodes $u_1,...,u_{m+1}\in J(x_0;T)$.
\par We note that, by \rf{C-1} and \rf{J-M2},
\bel{J-MN}
\#J(x_0;T)\ge \min(m+2,\dim Y+1)=m+2\,.
\ee
\par Let
$$
\Mc'=\bigcup_{i=1}^{m+1} \BR(u_i)\,.
$$
\par Clearly, $\Mc'\ni x_0$. Furthermore, by \rf{J-MN}, $\#\Mc'<\#\Mc=k+1$, so that, by the induction hypothesis, there exists a mapping $f_{\Mc'}:\Mc'\to Y$ such that
\bel{FP-11}
f_{\Mc'}(z)\in F(z)~~~\text{for all}~~~z\in\Mc'
\ee
and
\bel{FP-12}
\|f_{\Mc'}\|_{\Lip((\Mc',\,\rho),Y)}\le\gamma(k)\,.
\ee
\par Let us prove that
$$
f_{\Mc'}(x_0)\in\bigcap_{i=1}^{m+1}\,\Or(u_i)\,.
$$
\par Indeed, since $\tF=F$ on $\Mc'$, for every $u\in J(x_0;T)$, the restriction $\tF|_{\BR(u)}=F|_{\BR(u)}$, so that the mapping $g=f_{\Mc'}|_{\BR(u)}$ is a {\it selection} of $\tF|_{\BR(u)}$. It is also clear that
$$
\|g\|_{\Lip((\BR(u),\,\rho),Y)}=
\|f_{\Mc'}\|_{\Lip((\BR(u),\,\rho),Y)}\le
\|f_{\Mc'}\|_{\Lip((\Mc',\,\rho),Y)}\le\gamma(k)\le 2\gamma(k)\theta\,.
$$
\par This proves \rf{OR-3} and \rf{M-INT} in the case under consideration.
\medskip
\par {\it The second case: $m=\dim Y$.}
\par In this case, by \rf{J-M2},
\bel{SC-J}
\#J(x_0;T)\ge \min(m+2,\dim Y+1)=m+1\,.
\ee
Furthermore,
$$
\tF(x_0)=Y~~~\text{and}~~~\tF(u)=F(u),~~u\ne x_0\,.
$$
\par Note that in this case $F(x_0)$ and all $\Or(u)$, $u\in J(x_0;T)$, are convex subsets of the Banach space $Y$ with $\dim Y=m$. Therefore, by the Helly's Theorem \reff{HT-2}, \rf{M-INT} holds whenever for arbitrary nodes $u_1,...,u_{m+1}\in J(x_0;T)$ both
\bel{F-O1}
F(x_0)\bigcap\left\{
\bigcap_{i=1}^m\,\Or(u_i)
\right\}\ne\emp
\ee
and
\bel{NF-2}
\bigcap_{i=1}^{m+1}\,\Or(u_i)\ne\emp
\ee
hold.
\smallskip
\par Let us prove \rf{F-O1}. We define a set $\Mc'\subset\Mc$ by
$$
\Mc'=\bigcup_{i=1}^{m}\,\BR(u_i)\,.
$$
Then $x_0\in\Mc'$, and, by \rf{SC-J}, $\#\Mc'<\#\Mc$. This and the induction assumption imply the existence of a mapping $f_{\Mc'}:\Mc'\to Y$ satisfying \rf{FP-11} and \rf{FP-12}.
\par Then $f_{\Mc'}(x_0)\in F(x_0)$. Let us prove that
\bel{F-14}
f_{\Mc'}(x_0)\in \Or(u_i)~~~\text{for all}~~~i=1,...,m.
\ee
In fact, let $i\in\{1,...,m\}$ and let
$g_i=f_{\Mc'}|_{\BR(u_i)}$. Then $g_i$ is a selection of $\tF$ on $\BR(u_i)$ because $f_{\Mc'}$ is a selection of $F$ on $\BR(u_i)$ and $F|_{\BR(u_i)}\subset \tF|_{\BR(u_i)}$ (see \rf{F-W1-N}). Furthermore,
$$
\|g_i\|_{\Lip((\BR(u_i),\,\rho),Y)}\le \|f_{\Mc'}\|_{\Lip((\Mc',\,\rho),Y)}
\le \gamma(k)\le 2\gamma(k)\theta\,.
$$
Hence, $g_i(x_0)=f_{\Mc'}(x_0)\in \Or(u_i)$ (see\rf{OR-U}) proving \rf{F-14} and \rf{F-O1}.
\smallskip
\par Let us prove \rf{NF-2}. We put
$$
\Mc'=\bigcup_{i=1}^{m+1}\, \left(\BR(u_i)\setminus\{x_0\}
\right)\,.
$$
Since $x_0\notin\Mc'$, the cardinality $\#\Mc'<\#\Mc$, so that, by the induction assumption, there exists a mapping $f_{\Mc'}:\Mc'\to Y$ satisfying \rf{FP-11} and \rf{FP-12}.
\par We pick $u_0\in J(x_0;T)$ satisfying
\bel{CH-11}
\rho(u_0,x_0)=\min_{u\in J(x_0;T)}\rho(u,x_0)\,.
\ee
Let us show that
\bel{2-A}
f_{\Mc'}(u_0)\in\bigcap_{i=1}^{m+1}\,\Or(u_i)\,.
\ee
\par Indeed, fix $i\in\{1,...,m+1\}$ and define a mapping $g_i:\BR(u_i)\to Y$ by letting
\bel{G-L}
g_i(z)=\left \{
\begin{array}{ll}
f_{\Mc'}(z),& \text{if}~~
z\in \BR(u_i)\setminus\{x_0\},\vspace*{2mm}\\
f_{\Mc'}(u_0),& \text{if}~~z=x_0.
\end{array}
\right.
\ee
\par Since
$$
\tF|_{\BR(u_i)\setminus\{x_0\}}=
F|_{\BR(u_i)\setminus\{x_0\}}~~~\text{and}~~~\tF(x_0)=Y,
$$
by \rf{FP-11} and \rf{G-L}, the mapping $g_i$ is a {\it selection} of $\tF|_{\BR(u_i)}$. Furthermore,
\bel{XY-11}
\|g_i(x)-g_i(y)\|=\|f_{\Mc'}(x)-f_{\Mc'}(y)\|\le \gamma(k)\rho(x,y)~~~\text{for all}~~x,y\in \BR(u_i)\setminus\{x_0\}\,.
\ee
\par Now, let $y\in \BR(u_i)\setminus\{x_0\}$. Then, by \rf{FP-12},
$$
\|g_i(x_0)-g_i(y)\|=\|f_{\Mc'}(u_0)-f_{\Mc'}(y)\|\le \gamma(k)\rho(u_0,y)\le \gamma(k)\{\,\rho(u_0,x_0)+\rho(x_0,y)\}
$$
so that, by \rf{CH-11},
$$
\|g_i(x_0)-g_i(y)\|\le \gamma(k)\{\,\rho(u_i,x_0)+\rho(x_0,y)\}
\le \gamma(k)\{\,d_T(u_i,x_0)+d_T(x_0,y)\}\,.
$$
Since $y\in \BR(u_i)\setminus\{x_0\}$, by \rf{T-B1}, the unique path joining $y$ to $x_0$ in $T$ includes $u_i$. Hence,
$$
d_T(u_i,x_0)\le d_T(x_0,y).
$$
This inequality together with \rf{M-TH} imply that
$$
\|g_i(x_0)-g_i(y)\|\le 2\gamma(k)\,d_T(x_0,y)\le 2\gamma(k)\,\theta\,\rho(x_0,y)\,.
$$
By this inequality and by \rf{XY-11},
$$
\|g_i\|_{\Lip((\BR(u_i),\rho),Y)}\le 2\gamma(k)\,\theta\,.
$$
\par Hence, by \rf{OR-U}, $g_i(x_0)=f_{\Mc'}(u_0)\in\Or(u_i)$ proving \rf{2-A} and \rf{NF-2}.
\medskip
\par Thus, we have proved that \rf{M-INT} holds so that there exists a point
\bel{A0}
a_0\in
F(x_0)\bigcap\left\{
\bigcap_{u\in J(x_0;T)}\,\Or(u)
\right\}\,.
\ee
\par Let $u\in J(x_0;T)$. Since $a_0\in\Or(u)$, by \rf{OR-U}, there exists a mapping $g_u:\BR(u)\to Y$ such that $g_u(x_0)=a_0$,
\bel{GU-1}
g_u(y)\in F(y)~~~\text{for all}~~y\in\BR(u)\setminus\{x_0\},
\ee
and
\bel{GU-2}
\|g_u\|_{\Lip((\BR(u),\rho),Y)}\le 2\gamma(k)\,\theta\,.
\ee
\par Finally, we define a mapping $f:\Mc\to Y$ by letting
\bel{DF-4}
f(x_0)=a_0~~~\text{and}~~~
f|_{\BR(u)\setminus\{x_0\}}
=g_u|_{\BR(u)\setminus\{x_0\}}~~~\text{for every}~~~u\in J(x_0;T)\,.
\ee
\par Note that, by $(\bullet1)$, the mapping $f$ is well defined on $\Mc$. By \rf{A0} and \rf{GU-1}, the mapping $f$ is a {\it selection} of $F$ on $\Mc$. Let us show that
\bel{F-L1}
\|f(x)-f(y)\|\le  2\gamma(k)\,\theta^2\,\rho(x,y)~~~\text{for every}~~~x,y\in\Mc\,.
\ee
\par Let $u\in J(x_0;T)$ and let $x,y\in\BR(u)$. Then, by \rf{DF-4},
$
f|_{\BR(u)}=g_u|_{\BR(u)}
$
(recall that $f(x_0)=a_0=g_u(x_0)$) so that \rf{F-L1} follows from \rf{GU-2}.
\par Now let $x\in\BR(u_1)\setminus\{x_0\}$ and $y\in\BR(u_2)\setminus\{x_0\}$ where $u_1,u_2\in J(x_0;T)$, $u_1\ne u_2$. Then
\be
\|f(x)-f(y)\|&=&\|g_{u_1}(x)-g_{u_2}(y)\|\le
\|a_0-g_{u_1}(x)\|+\|a_0-g_{u_2}(y)\|\nn\\
&=&
\|g_{u_1}(x_0)-g_{u_1}(x)\|+\|g_{u_2}(x_0)-g_{u_2}(y)\|\nn\\
&\le&
2\gamma(k)\,\theta\,\{\,\rho(x_0,x)+\rho(x_0,y)\}\,.\nn
\ee
Hence,
$$
\|f(x)-f(y)\|\le
2\gamma(k)\,\theta\,\{d_T(x_0,x)+d_T(x_0,y)\}
$$
so that, by $(\bullet 1)$ and $(\bullet 2)$,
$$
\|f(x)-f(y)\|\le 2\gamma(k)\,\theta\,d_T(x,y)
\le 2\gamma(k)\,\theta^2\,\rho(x,y)
$$
proving \rf{F-L1}.
\par The proof of Theorem \reff{D-2M} is complete.\bx
\medskip
\par {\it Proof of Theorem \reff{MAIN-FP}.} Suppose that $\rho$ is a finite pseudometric, i.e., condition \rf{FIN-A} holds.
\par Let $\Mc'$ be an arbitrary subset of $\Mc$ consisting of at most $\ks$ points. Then, by the theorem's hypothesis, for every set $S\subset\Mc'$ with $\#S\le\FN$, the restriction $F|_S$  has a Lipschitz selection $f_{S}:S\to\BS$ with $\|f_S\|_{\Lip(S,\BS)}\le 1$. Hence, by Theorem \reff{D-2M}, the restriction $F|_{\Mc'}$ of $F$ to $\Mc'$ has a Lipschitz selection $f_{\Mc'}:\Mc'\to\BS$ whose seminorm satisfies $\|f_{\Mc'}\|_{\Lip(\Mc',\BS)}\le \gamma$ where $\gamma$ is a constant depending only on $m$ and $\#\Mc'$. Since $\#\Mc'\le \ks$ and $\ks$ depends only on $m$, the constant $\gamma$ depends only on $m$ as well.
\par Hence, by Proposition \reff{PM-SP}, $F$ has a Lipschitz selection $f:\Mc\to Y$ with
$\|f\|_{\Lip(\Mc,Y)}\le\gamma_2\gamma$
where $\gamma_2$ is a constant depending only on $m$.
\par This completes the proof of Theorem \reff{MAIN-FP} for the case of a {\it finite} pseudometric $\rho$.
\par Let us prove Theorem \reff{MAIN-FP} for an arbitrary
pseudometric $\rho:\Mc\times\Mc\to\R_+\cup\{+\infty\}$ which may admit the value $+\infty$.
\smallskip
\par Let us introduce a binary relation ``$\sim$'' on $\MS$ by letting
$$
x\sim y~~~\Longleftrightarrow~~~\rho(x,y)<\infty\,.
$$
Clearly, ``$\sim$'' satisfies the axioms of an equivalence relation, i.e., it is reflexive, symmetric and transitive. Given $x\in\MS$, by $[x]=\{y\in\MS:~y\sim x\}$ we denote the equivalence class of $x$. Let
$$
[\Mc]=\Mc\,/\sim\,\,=\,\{\,[x]: x\in\MS\,\}
$$
be the corresponding quotient set of $\MS$ by ``$\sim$'', i.e.,  the family of all possible equivalence classes of $\MS$ by ``$\sim$''.
\par Let $U\in[\Mc]$ be an equivalence class, and let
$$
\rho_U=\rho|_{U\times U}.
$$
Then
\bel{U-MS}
\rho_U(x,y)=\rho(x,y)<\infty~~~\text{for all}~~~x,y\in U\,.
\ee
\par Let $F_U=F|_U$. Clearly, the hypothesis of Theorem \reff{MAIN-FP} holds for the pseudometric space $(U,\rho_U)$ and set-valued mapping $F_U:U\to\KMY$:
{\it for every subset $U'\subset U$ consisting of at most $\FN$ points, the restriction $F_U|_{U'}$ of $F_U$ to $U'$ has a Lipschitz selection $f_{U'}:U'\to\BS$ with  $\|f_{U'}\|_{\Lip((U',\rho_U),\BS)}\le 1$.}
\par This property and \rf{U-MS} enable us to apply to $(U,\rho_U)$ and $F_U$ the variant of Theorem \reff{MAIN-FP} for finite pseudometrics proven above.
Thus, we produce a mapping $f_U:U\to Y$ such that
\bel{SL-U}
f_U(x)\in F_U(x)=F(x)~~\text{for every}~~~x\in U\,,
\ee
and
\bel{LC-U}
\|f_U(x)-f_U(y)\|\le \gamma \rho_U(x,y)=\gamma \rho(x,y)~~~\text{for all}~~~x,y\in U\,.
\ee
Here $\gamma=\gamma(m)$ is a constant depending only on $m$.
\par We define a mapping $f:\Mc\to Y$ by letting
$$
f(x)=f_{[x]}(x),~~~x\in \Mc\,.
$$
\par Clearly, by \rf{SL-U}, $f(x)\in F(x)$ for every $x\in\Mc$, i.e., $f$ is a selection of $F$ on $\Mc$. Let us prove that $\|f\|_{\Lip(\Mc,Y)}\le \gamma$. Indeed, if
$x,y\in \Mc$ and $[x]=[y]$, then, by \rf{LC-U},
$$
\|f(x)-f(y)\|\le \gamma\,\rho(x,y)\,.
$$
\par If $[x]\ne[y]$, then $\rho(x,y)=+\infty$, so the above inequality trivially holds.
\par The proof of Theorem \reff{MAIN-FP} is complete. \bx
\medskip
\par We finish this section with a variant of our main result, Theorem \reff{MAIN-FP}, related to the case of {\it finite} pseudometric spaces.
\begin{theorem}\lbl{FINITE} Let $(\Mc,\dm)$ be a finite pseudometric space, and let $F:\Mc\to\CNMY$ be a set-valued mapping from $\Mc$ into the family $\CNMY$ of all convex subsets of $Y$ of affine dimension at most $m$. Assume that, for every subset $\Mc'\subset\MS$ consisting of at most $\FN$ points, the restriction $F|_{\Mc'}$ of $F$ to $\Mc'$ has a Lipschitz selection $f_{\Mc'}:\Mc'\to\BS$ with $\|f_{\Mc'}\|_{\Lip(\Mc',\BS)}\le 1$.
\par Then $F$ has a Lipschitz selection $f:\MS\to\BS$ with $\|f\|_{\Lip(\MS,\BS)}$ bounded by a constant depending only on $m$.
\end{theorem}
\par {\it Proof.} We prove this theorem following the  scheme of the proof of Theorem \reff{MAIN-FP}. In particular, for a finite pseudometric $\rho$ we use in the proof Proposition \reff{FN-PM} and the constant $\kn$ rather than Proposition \reff{PM-SP} and $\ks$ respectively. We also literally follow the proof of Theorem \reff{MAIN-FP} for the general case of an arbitrary pseudometric $\rho:\Mc\times\Mc\to\R_+\cup\{+\infty\}$.  \bx
\SECT{7. A Steiner-type point of a convex body.}{7}
\addtocontents{toc}{7. A Steiner-type point of a convex body.\hfill\thepage\par\VST}

\indent\par {\bf 7.1. Barycentric Selectors.}
\addtocontents{toc}{~~~~7.1. Barycentric Selectors.\hfill \thepage\par}\medskip

\par For the reader's convenience, in this section we briefly describe the construction of the Steiner-type mapping $\SX:\KX\to\BS$ satisfying conditions (i)  and (ii) of Theorem \reff{ST-P}. See \cite{S04} for the details.
\par Consider the metric space $(\KX,\dhf)$ of all non-empty finite dimensional convex compact subsets of $\BS$ equipped with the Hausdorff distance. Let $S:\KX\to\BS$ be a mapping such that $S(K)\in K$ for every $K\in\KX$. We refer to $S$ as a {\it selector}. We note that for an infinite dimensional Banach space $\BS$ there {\it does not exist} a $\dhf$-Lipschitz continuous selector which is defined on {\it all} of the family $\KX$. See \cite{PY1}. Theorem \reff{ST-P} implies that, in contrast to this negative result, {\it there exists a selector $\SX:\KX\to \BS$ which is Lipschitz continuous on every family $\KM$, $m\in\N$.}
\smallskip
\par For the case of a Hilbert space $H$ the classical
{\it Steiner point} \cite{St} $s(K)$
of a convex body $K\subset H$ provides such a
selector. Recall that if $K\in\KH$ is a subset
of an $n$-dimensional subspace $L\subset H$,
then its Steiner point $s(K)$ is defined by the formula
$$
s(K)= n\int\limits_{\SP_H\cap L}u\, h_K(u)\,d\sigma(u).
$$
Here $\SP_H$ is the unit sphere in $H$,
$h_K(u)=\sup\{\langle u,x\rangle:~x\in K\}$
is the support function of $K$, and $\sigma$ denotes the normalized Lebesgue measure on $\SP_H\cap L$ which is calculated with respect to an arbitrary predetermined Euclidean basis for $L$.\smallskip
\par The Steiner point map is a continuous selector which is {\it additive} with respect to Minkowski addition and
{\it commutes with the affine isometries of $H$}.
These properties uniquely define the Steiner point and show that $s(K)$ is well-defined, i.e., its definition does not depend on the choice of the finite dimensional subspace $L$ containing $K$, or on the choice of the Euclidean basis of $L$ .
\par Moreover, the Steiner point map is Lipschitz continuous on every $n$-dimensional subspace of $H$ and
its Lipschitz constant equals $c_n=2\pi^{-\frac{1}{2}}\Gamma\left(\frac{n}{2}+1\right)/
\Gamma\left(\frac{n+1}{2}\right)\sim\sqrt{n}$. (This value is sharp and is the smallest possible for selectors from $\KH$ to $H$. See \cite{BL,Pos,PY1,Vit}.)
Since the linear hull of any two convex compact subsets $K_1$ and $K_2$ has dimension not exceeding
$n=\dim K_1+\dim K_2+2$, we see that
$$
\|s(K_1)-s(K_2)\|\le c_n \dhf(K_1,K_2)
$$
for every $K_1,K_2\in\KH$. Consequently, the restriction
$s|_{\Kc_m(H)}$ is Lipschitz continuous for every
$m\in \N$. For these and other properties of the Steiner point map we refer the reader to \cite{BL,Sh,Sc1,PY1,Vit} and references therein.
\par  Unfortunately there does not seem to be any obvious way of generalizing the Steiner point construction
to the case of a non-Hilbert Banach space. (We refer the reader to \cite{PY1}, for some partial results which indicate the difficulties of making such a generalization.) 
\medskip
\par The construction of the mapping $\SX:\KX\to\BS$
satisfying conditions (i) and (ii) of Theorem \reff{ST-P} relies on some ideas related to using barycenters rather than Steiner points. Even for the case of a Hilbert space the Steiner point map and the selector $\SX$ are distinct. We call this selector a {\it Steiner-type selector}.
\par We construct this selector by induction on dimension of subsets from the family $\KX$.
\par Without loss of generality we may assume that $\BS$ is a space $\BU$ of bounded functions defined on a certain set $U$. Indeed, any Banach space $Y$ isometrically embeds in a certain $\BU$. Therefore, if we produce a selector for $\BU$, we produce a selector for $\BS$. \smallskip
\par For the family $\Kc_0(\BS)$ of all singletons in $\BS$ we define $\SX$ by letting $\SX(\{x\})=x,~x\in \BS$. Clearly, in this case $\SX$ satisfies all the conditions of Theorem \reff{ST-P} with the constant $\gamma_2=1$.
\par Let us assume that for an integer $m\ge 0$
the mapping $\SX$ is defined on the family $\KM$ and satis\-fies the following conditions: (i). $\SX$ is a selector on $\KM$, i.e., $\SX(K)\in K$ for every ${K\in\KM}$; (ii). $\SX$ is Lipschitz on $\KM$ with respect to the Hausdorff distance, i.e., for every $K_1,K_2\in\KM$
\bel{SLM}
\|\SX(K_1)-\SX(K_2)\|\le \gamma(m)\,\dhf(K_1,K_2).
\ee
\par We construct the required selector $\SX$ on $\KMO$ in two steps. At the first step we extend the mapping $\SX$ from  $\KM$ to $\KMO$ with preservation of the Lipschitz condition \rf{SLM}. In fact, $(\KM,\dhf)$ is a metric subspace of the metric space $(\KMO,\dhf)$, and $\SX$ is a Lipschitz mapping from $(\KM,\dhf)$ into $\BS$. Recall that we identify $\BS$ with a space $\BU$ of bounded functions  on a set $U$. The space $\BU$ possesses the following well known universal extension property: every Lipschitz mapping from a subspace of a metric space to $\BU$ can be extended to all of the metric space with preservation of the Lipschitz constant.
\par Thus there exists a mapping $\SW:\KMO\to\BS$ such that 
\bel{SW-SL}
\SW(K)=\SX(K)~~~\text{for each}~~~K\in\KM,
\ee
and
\bel{SW-L}
\|\SW(K_1)-\SW(K_2)\|\le \gamma(m)\,\dhf(K_1,K_2)~~~
\text{for every}~~~K_1,K_2\in\KMO.
\ee
\par We refer to the mapping $\SW:\KMO\to\BS$ as a {\it pre-selector}. This name is motivated by the fact that in general $\SW(K)\notin K$ whenever $K$ is an  $(m+1)$-dimensional convex compact set in $\BS$. Nevertheless, we show below that for each $K\in\KMO$ its pre-selector $\SW(K)$ lies in a certain sense rather ``close'' to $K$. This enables us to ``correct'' the position of  $\SW(K)$ with respect to the set $K$, and obtain in this way the required Lipschitz selector defined on all of the family $\KMO$.
\par We make this ``correction'' at the second step of the procedure. An important ingredient of our construction at this step is the notion of the Kolmogorov $m$-width of the set $K$. This geometric characteristic of $K$ is defined by 
\bel{WD}
d_m(K)=\inf\{\ep>0:L+ \BL(0,\ep)\supset K, L\in\AM\}.
\ee
Recall that $\AM$ denotes the family of all affine subspaces of $\BS$ of dimension at most $m$. It can be readily seen that $d_m$ satisfies the Lipschitz condition with respect to the Hausdorff distance, i.e.,
\bel{W-HD}
|\,d_m(K_1)-d_m(K_2)\,|\le \dhf(K_1,K_2),~~~K_1,K_2\in \KX\,.
\ee
\par Then given $K\in\KMO$ we define a set $\tK$ by
\bel{KW-D}
\tK=K\,\capsm\, \BL(\SW(K),R(K)).
\ee
Here $R(K)=\gmh\,d_m(K)$ and $\gmh=\gmh(m)>0$ is a certain constant depending only on $m$ which will be determined below. Recall that given $x\in Y$ and $r>0$, by $\BL(x,r)$
we denote a closed ball in $Y$ with center $x$ and radius $r$.
\par Finally, we put
\bel{SX-F}
\SX(K)=b(\tK)
\ee
where $b(\cdot)$ denotes the {\it barycenter}
(center of mass) of a finite dimensional set in $\BS$.
\smallskip
\par We prove that for $\gmh=\gmh(m)>0$ big enough the set $\tK\ne\emp$ for each $K\in\KMO$. Hence $\SX(K)=b(\tK)\in \tK\subset K$ proving that $\SX$ is a {\it selector on all of the family $\KMO$}.
\par Let us show that $\SX$ is a Lipschitz mapping on $\KMO$. In view of formula \rf{SX-F} it is natural to ask what are the $\dhf$-Lipschitz properties of the barycenter. We note that the barycentric map $b:\KX\to \BS$ is a continuous selector (\cite{SW}), but (unlike the Steiner point map in Hilbert spaces) it is not a Lipschitz map on the family $\Kc_n(\BS)$ for every $n>1$.
\par However, the barycentric map does have a certain ``Lipschitz-like" property, where the usual Lipschitz constant is augmented by an additional factor which depends on a certain ``geometrical" quantity associated with sets $K\in\KX$. To define this quantity, for each $K\in\KX$, we first choose some Lebesgue measure $\lambda$ on $\affspan(K)$, the affine hull of $K$. Then we define the {\it regularity coefficient} of $K$ to be the number
$$
\delta_K=\lambda\left(B^{(K)}\,\capsm\affspan(K)\right)
/\lambda\left(K\right),
$$
where $B^{(K)}$ denotes a ball (with respect to $\|\cdot\|$) of {\it minimal radius} among all
balls which contain $K$ and whose centers lie in $\affspan(K)$. It is proven in \cite{S04} that for every $K_1,K_2\in\KX$
$$
\|b(K_1)-b(K_2)\|\le
\gamma\cdot(\delta_{K_1}+\delta_{K_2})
\dhf(K_1,K_2)
$$
where $\gamma$ is a constant depending only on $\dim K_1$ and $\dim K_2$.
\par We apply this inequality to the mapping $\SX$ defined by \rf{SX-F} and get
$$
\|\SX(K_1)-\SX(K_2)\|=\|b(\tK_1)-b(\tK_2)\|\le
\gamma_1(m)\cdot(\delta_{\tK_1}+\delta_{\tK_2})
\dhf(\tK_1,\tK_2)
$$
provided $K_1,K_2$ are arbitrary sets from $\KMO$.
\par It remains to estimate the order of magnitude of the two quantities: the regularity coefficient $\delta_{\tK}$ for each $K\in\KMO$, and the Hausdorff distance $\dhf(\tK_1,\tK_2)$ for every $K_1,K_2\in\KMO$.
\par We show that for an appropriate choice of the constant $\gmh=\gmh(m)$ the regularity coefficient
$$
\delta_{\tK}\le \gamma_2(m)~~~\text{provided}~~~ K\in\KMO.
$$
The proof of this property relies on equality \rf{SW-SL} and inequality \rf{SW-L}, definition \rf{WD}, and the following important property of the barycenter due to Minkowski \cite{Min}:
{\it there exists a constant $\alpha=\alpha(m)\ge 1$ such that for every set $G\in\KMO$ the following inclusion
$$
\BL(b(G),d_{m}(G)/\alpha)\,\capsm\affspan(G)\subset G
$$
holds.}
\par Then we show that
\bel{F-D1}
\dhf(\tK_1,\tK_2)\le \gamma_3(m)\dhf(K_1,K_2)
\ee
for all $K_1,K_2\in\KMO$. The proof of this inequality is based on the following geometrical result \cite{PR}: {\it Suppose that $G\cap \BL(a,r)\ne\emp$ where $G\subset\BS$ is a convex set, $a\in\BS$ and $r>0$. Then for every $s>0$}
$$
(G+\BL(0,s))\,\capsm\, (\BL(a,2r)+\BL(0,s))\subset G\,\capsm\,
\BL(a,2r)+\BL(0,9s).
$$
\par This inclusion implies the following statement:
{\it Let $G_i\subset \BS$ be a convex set and let  $\BL(a_i,r_i)$ where $a_i\in \BS$, $r_i>0,$ be a ball in $\BS$ such that $G_i\cap \BL(a_i,r_i)\ne \emptyset$, $i=1,2$.
Then}
\bel{G-1}
\dhf(G_1\,\capsm\,\BL(a_1,2r_1),G_2\,\capsm\, \BL(a_2,2r_2))
\le 18\,(\dhf(G_1,G_2)+\|a_1-a_2\|+|r_1-r_2|)\,.
\ee
\par We recall that the radius $R(K)$ from the definition \rf{KW-D} is defined as $R(K)=\gmh\,d_m(K)$. Let us choose the constant $\gmh=\gmh(m)$ in such a way that
$$
K\cap \BL(\SW(K),R(K)/2)\ne\emp.
$$
This enables us to apply inequality \rf{G-1} to the sets $G_i=K_i$, points $a_i=\SW(K_i)$ and radii $r_i=R(K_i)$, $i=1,2$. By this inequality,
\be
\dhf(\tK_1,\tK_2)&=&\dhf(K_1\,\capsm\, \BL(\SW(K_1),R(K_1)),K_2\,\capsm\, \BL(\SW(K_2),R(K_2)))\nn\\
&\le& 18\,(\dhf(K_1,K_2)+\|\SW(K_1)-\SW(K_2)\|+|R(K_1)-R(K_2)|)\,.
\nn
\ee
Combining this inequality with inequalities \rf{SW-L} and \rf{W-HD}, we obtain the required estimate \rf{F-D1}.
\par This concludes our sketch of the proof of Theorem \reff{ST-P}.\bx
\bigskip\medskip

\indent\par {\bf 7.2. Further properties of Steiner-type selectors.}
\addtocontents{toc}{~~~~7.2. Further properties of Steiner-type selectors.\hfill \thepage\par\VSU}\medskip
\par In this section we will review several additional properties of the Steiner-type selector of a finite dimensional convex body in a Banach space. See \cite{S04}.
\smallskip
\par \textbullet~ {\bf A Steiner-type selector for the family of all finite dimensional convex sets.}
\smallskip
\par It is shown in \cite{S04} that the Steiner-type selector described in Section 7.1 can be extended from the family $\KX$ of all non-empty convex {\it compact} finite dimensional subsets of $Y$ to the family $\CNY$ of {\it all} non-empty convex finite dimensional subsets of $Y$ with preservation of its $\dhf$-Lipschitz properties.
\par We define this extension as follows: given a set $K\in \CNY$ we put
\bel{ST-ALL}
\SX(K)=\left \{
\begin{array}{ll}
\SX(K^{\cl})\,,& \text{if}~~K~~\text{is bounded}, \vspace*{4mm}\\
\SX\left(\left[K\bigcap \left (2\dist(0,K)\right)B_Y\right]^{\cl}\right)\,,
&\text{if}~~K~~\text{is unbounded}.
\end{array}
\right.
\ee
Recall that the sign $\operatorname{cl}$ denotes the closure of a set in $Y$. Since $K^{\cl}\in\KX$ whenever ${K\in\CNY}$ is bounded and
$$
[K\,\capsm (2\dist(0,K))B_Y]^{\cl}\in\KX
$$
whenever $K\in\CNY$ is unbounded, the mapping \rf{ST-ALL} is well defined on $\CNY$. Below we note that $\SX(K^{\cl})\in K$ (see \rf{K-INT}). Hence, $\SX(K)\in K$ for each $K\in\CNY$ proving that $\SX$ is a {\it selector}.
\par Furthermore, we have
$$
\|\SX(K)-\SX(K')\|\le \gamma(K,K')\cdot \dhf(K,K')
$$
with $\gamma(K,K')$ depending only on the dimensions of $K,K'$. (Note that here $\dhf(K,K')$ may be infinite.)
This fact immediately follows from part (ii) of Theorem \reff{ST-P} and inequality \rf{G-1}.
\medskip
\par \textbullet~ {\bf Two important properties of the Steiner-type selector.}
\smallskip
\par Let $K\in\CNY$ be a {\it bounded set}. Then
\smallskip
\par  (i) $\SX(\tau K+a)=\tau \SX(K)+a$ for every $a\in Y$ and every $\tau\in\R$. In other words, the Steiner-type selector $\SX$ is invariant with respect to dilations and shifts;
\smallskip
\par (ii) There is an ellipsoid $\Ec_{K}$ centered at $\SX(K)$ such that
\bel{ELLS}
\Ec_{K}\subset K\subset \gamma\circ \Ec_{K},
\ee
where $\gamma=\gamma(\dim K)$ is a constant depending only on dimension of $K$. Here given a centrally symmetric subset $A\subset Y$ and a positive constant $\lambda$ we let $\lambda\circ A$ denote the dilation of $A$ with respect to its center by a factor of $\lambda$.
\smallskip
\par Note that, by property (i), $\SX(K)$ coincides with
the center of $K$ for every bounded {\it centrally symmetric}
set $K\in\CNY$. Furthermore, property (ii) implies that the point $\SX(K)$ is located rather ``deeply" in the interior of the set $K$. In particular, by this property,
\bel{K-INT}
\SX(K^{\cl})\in K~~\text{for every bounded set}~~ K\in\CNY\,.
\ee
Note, by way of comparison, that the Steiner point of a set always belongs to the relative interior of the set (\cite{Sh}), but, as noted in \cite{Sc3}, an estimate for the position of the Steiner point (e.g., similar to \rf{ELLS}) seems to be unknown.
\medskip
\par \textbullet~ {\bf The centroid of a parallel body.}
\smallskip
\par Let us describe another construction of a barycentric selector which is $\dhf$-Lipschitz continuous on the family $\Kc(H)$ of all compact convex subsets of a finite dimensional Euclidean space.
\par Let $Y$ be a {\it Minkowski space}, i.e., a finite dimensional Banach space. Following an idea of Aubin and Cellina in \cite{AC},
we define a mapping $\SB:\Kc(Y)\to Y$ by letting
$$
\SB(K):=b(K+(\diam K)B_Y)\,.
$$
\par Let $\lambda>0$ and let $K\in\Kc(Y)$. We refer to the sets $K+\lambda\,B_Y$ as {\it parallel bodies} (with respect to $K$ and $\lambda$), and call $\SB(K)$ {\it the centroid of the parallel body}.
\par It is proven in \cite{S04} that $\SB:\Kc(Y)\to Y$ is a $\dhf$-Lipschitz continuous mapping whose $\dhf$-Lipschitz seminorm is bounded by a constant depending only on $\dim Y$. Furthermore, similar to the Steiner point, $\SB$ commutes with affine isometries and dilations of $Y$.
\par It is shown in \cite{AC} that $\SB(K)\in K$ for each compact convex $K\subset Y$ provided $Y$ is a finite dimensional {\it Euclidean space}. Thus for such $Y$,
\bel{LS-PB}
\text{the centroid of the parallel body}~~\SB~~\text{is a}~\dhf\text{-Lipschitz selector on}~~\Kc(Y)\,.
\ee
\par Its Lipschitz seminorm satisfies the inequality $\|\SB\|_{\Lip(\Kc(Y),Y)}\le \gamma$ where $\gamma$ is a constant depending only on $\dim Y$.
\smallskip
\par We notice an interesting connection of the Steiner point map $s(K)$ with the centroids of the parallel bodies, see \cite{PRZ}: If $Y$ is a finite dimensional Euclidean space then for every $K\in\Kc(Y)$
$$
s(K)=\lim_{r\to\infty}b(K+rB_Y)\,.
$$
\par The statement \rf{LS-PB} leads us to the following problem: Given a Minkowski space $Y$, decide whether
$$
\SB~~\text{is a}~\dhf\text{-Lipschitz selector on the family}~~\Kc(Y).
$$
\par We know that $\|\SB\|_{\Lip(\Kc(Y),Y)}\le \gamma(\dim Y)$, so that the above problem is equivalent to the following one: Let $Y$ be a Minkowski space. Does the centroid of the parallel body satisfy
\bel{CNT-Y}
b(K+B_Y)\in K~~~\text{for every compact convex set}~~K\subset Y?
\ee
\par This problem has been studied by Gaifullin \cite{G}
who proved that \rf{CNT-Y} is true for every {\it two dimensional} Minkowski space $Y$. In particular, this implies \rf{LS-PB} proving that $\SB$ is a Lipschitz continuous selector for every $Y$ of $\dim Y=2$.
\par Another result proven in \cite{G} states that a Minkowski space $Y$ with $\dim Y>2$ satisfies \rf{CNT-Y} if and only if $Y$ is a Euclidean space, i.e., its unit ball $B_Y$ is an ellipsoid.
\par The paper \cite{G} also contains an example of a triangle $K$ in the space $Y=\ell^1_3$ (with the norm $\|y\|_{\ell^1_3}=|y_1|+|y_2|+|y_3|$ for each $y=(y_1,y_2,y_3)\in\R^3$)
such that
$$
b(K+B_Y)\notin \affspan(K).
$$
This shows that in general the answer to question \rf{CNT-Y} is negative whenever $\dim Y>2$.
\begin{remark}{\em Let us note a very simple Lipschitz selector for the space $Y=\ell_2^\infty$, i.e., the space $\R^2$ equipped with uniform norm $\|x\|=\max\{|x_1|,|x_2|\}$, $x=(x_1,x_2)\in\R^2$.
\par In this case given a compact convex set $K\subset\R^2$ we let $\Pi(K)$ denote the smallest (with respect to inclusion) rectangle with sides parallel to the coordinate axes containing $K$. Let
$$
S(K)=\cntr(\Pi(K))\,.
$$
\par Clearly, a rectangle $\Pi\supset K$ with sides parallel to the coordinate axes, coincides with $\Pi(K)$ if and only if
\bel{P-KL}
\text{Each side of}~~\Pi~~\text{has a common point with}~~ K.
\ee
\par Let us show that $S(K)\in K$ for each convex compact $K\subset\R^2$, i.e., $S$ is a selector.
\par Indeed, let $K\in\Kc(\R^2)$. Suppose that $S(K)=\cntr(\Pi(K))\notin K$. Without loss of generality we may assume that $S(K)=0$. Then, by the separation theorem, there exists a vector $a\in\R^2$ such that $\ip{a,x}>0$ for all $x\in K$. Here $\ip{\cdot,\cdot}$ is the standard inner product in $\R^2$.
\par Clearly, there exists a side of the rectangle $\Pi(K)$, say $[AB]$, such that $\ip{a,z}\le 0$ for every $z\in [AB]$. Hence, $[AB]\cap K=\emp$ which contradicts \rf{P-KL}.\smallskip
\par It can be also readily seen that for every two compact convex sets $K_1,K_2\subset\R^2$ the following inequality
$$
\|S(K_1)-S(K_2)\|\le\dhf(K_1,K_2)
$$
holds. Thus $S:\Kc(\ell_2^\infty)\to \ell_2^\infty$ is a Lipschitz selector whose Lipschitz seminorm equals $1$. \rbx}
\end{remark}
\SECT{8. Further results and comments.}{8}
\addtocontents{toc}{8. Further results and comments. \hfill\thepage\par\VST}

\indent\par {\bf 8.1. The sharp finiteness constants for $m=1$ and $m=2$.}
\addtocontents{toc}{~~~~8.1. The sharp finiteness constants for $m=1$ and $m=2$.\hfill \thepage\par}\medskip

\par In this subsection we briefly indicate the main ideas of the proof of Theorem \reff{SH-FP} for the cases $m=1$ and $m=2$, i.e., for set-valued mappings with one dimensional and two dimensional images. By \rf{FNCN},  $N(1,Y)=4$ provided $\dim Y\ge 2$ and $N(2,Y)=8$ provided
$\dim Y\ge 3$. Below we present examples of pseudometric spaces and set-valued mappings which show that these finiteness constants are sharp.
\medskip
\par \textbullet~ {\bf The sharp finiteness constant for $m=1$.}\smallskip
\par For simplicity, we will show the sharpness of  $N(1,Y)$ for the space
$Y=\ell^2_\infty=(\R^2,\|\cdot\|_\infty)$
where $\|\cdot\|_\infty$ is the uniform norm on the plane, $\|x\|_\infty=\max\{|x_1|,|x_2|\}$, $x=(x_1,x_2)\in\R^2$.
\par Recall that in this case $N(1,Y)=4$, so that, by Theorem \reff{MAIN-FP}, the following statement holds: Let $(\Mc,\rho)$ be a pseudometric space, and let $F:\Mc\to\Kc_1(\R^2)$ be a set-valued mapping which to every $u\in\Mc$ assigns a line segment
$$
F(u)=[a(u),b(u)]\subset\R^2.
$$
Suppose that for every {\it four point subset} $\Mc'\subset\Mc$ the restriction $F|_{\Mc'}$ has a Lipschitz selection $f_{\Mc'}:\Mc'\to\R^2$ with $\|f_{\Mc'}\|_{\Lip(\Mc',Y)}\le 1$. Then there exists a selection $f:\Mc\to\R^2$ of $F$ with $\|f\|_{\Lip(\Mc,Y)}\le \gamma$ where $\gamma$ is an absolute constant.
\vskip 2mm
\par Let us see that this statement is false whenever four point subsets in its formulation are replaced by {\it three point subsets}. We will show that, given $\lambda\ge 1$ there exists a four point metric space $(\wtM,\trh)$ and a set-valued mapping $\tF:\wtM\to\Kc_1(\R^2)$ such that the following is true: the restriction $\tF|_{\wtM'}$ of $\tF$ {\it to every three point subset} $\wtM'$ of $\wtM$ has a Lipschitz selection $f_{\wtM'}:\wtM'\to\R^2$ with $\|f_{\wtM'}\|_{\Lip(\wtM',Y)}\le 1$, but nevertheless
$$
\|f\|_{\Lip(\wtM,Y)}\ge\lambda~~~\text{\it for every selection}~~f~~\text{of}~~\tF.
$$
\par We define $(\wtM,\trh)$ and $\tF$ as follows. Let
\bel{L-LM}
L=2\lambda~~~~\text{and}~~~~\ve=1/L\,.
\ee
Let
\bel{M-TL}
\wtM=\{u_1,u_2,u_3,u_4\}~~~\text{where}~~~
u_1=1+\ve,~~u_2=1,~~u_3=-1,~~u_4=-1-\ve,
\ee
and let
\bel{RHO-T}
\trh(u_i,u_j)=|u_i-u_j|~~~\text{for all}~~~i,j=1,2,3,4.
\ee
\par Let
$$
\mA=(L,1),~~~~\mB=(-L,1),~~~~\mC=(-L,-1),~~~~\mD=(L,-1)\,.
$$
We define the set-valued mapping $\tF:\wtM\to\Kc_1(\R^2)$ by letting
\bel{FT-LS}
\tF(u_1)=[\mA\mB],~~~~\tF(u_2)=[\mA\mC],~~~~
\tF(u_3)=[\mB\mD],
~~~~\tF(u_4)=[\mA\mB]\,.
\ee
See Fig. 1 below.
\renewcommand{\figurename}{Fig.}
\medskip
\begin{figure}[h]
\center{\includegraphics[scale=0.33]{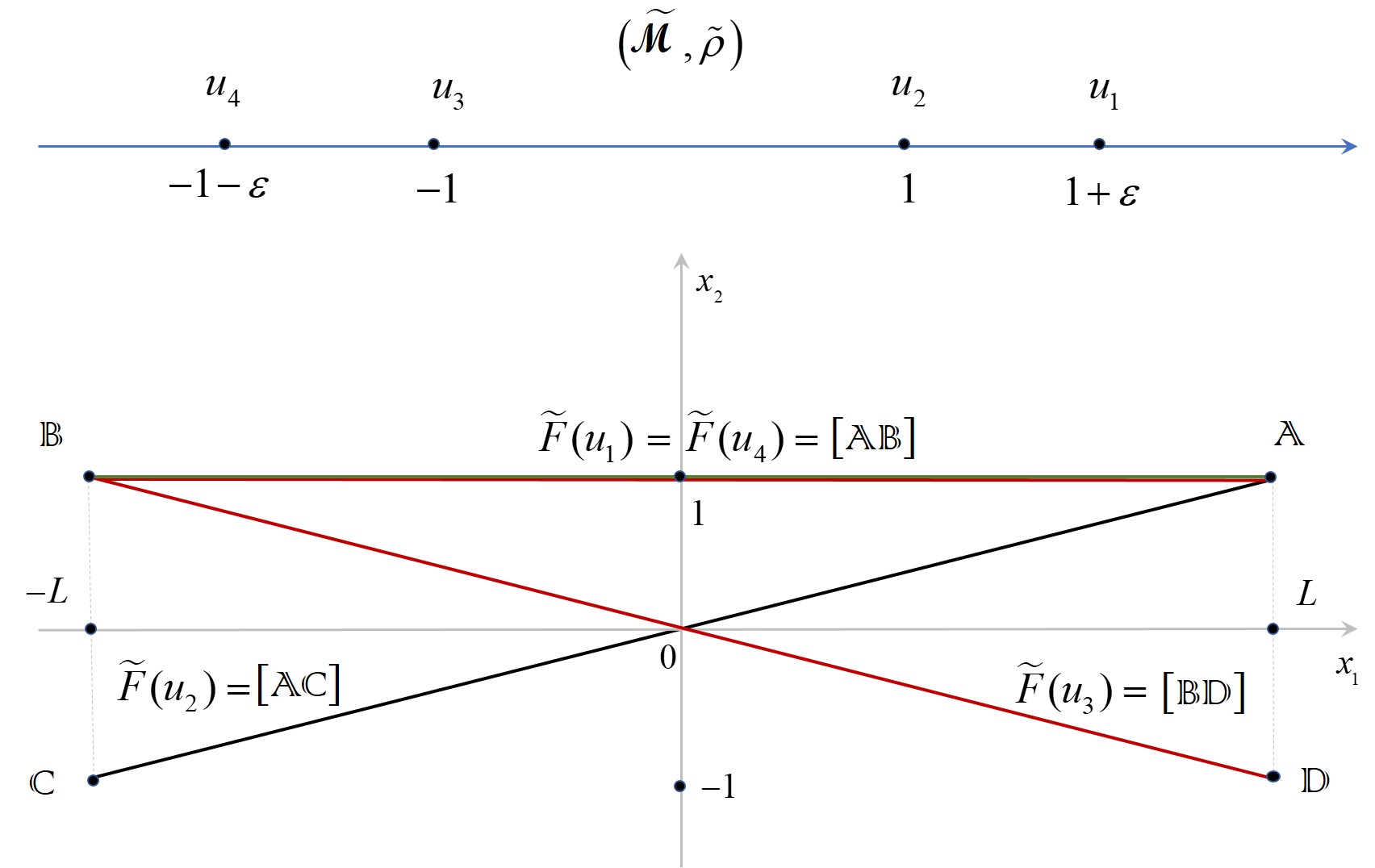}}
\caption{The metric space $(\wtM,\trh)$ and the set-valued mapping $\tF$.}
\end{figure}
\bigskip
\par Let $\wtM_i=\wtM\setminus\{u_i\}$, $i=1,2,3,4$. We prove that the restriction $\tF|_{\wtM_i}$ has a Lipschitz selection $f_i:\wtM_i\to\R^2$ with $\|f_i\|_{\Lip(\wtM_i,\LTI)}\le 1$.
\par For $i=1$ we define $f_1$ by
$$
f_1(u_2)=\mC,~~~~f_1(u_3)=f_1(u_4)=\mB\,.
$$
Clearly, $f_1$ is a selection of $\tF|_{\wtM_1}$. Furthermore, 
$$
\|f_1(u_2)-f_1(u_3)\|_\infty=\|\mC-\mB\|_\infty=2
=|u_2-u_3|=\trh(u_2,u_3),
$$
and
$$
\|f_1(u_2)-f_1(u_4)\|_\infty=\|\mC-\mB\|_\infty=2\le 2+\ve
=|u_2-u_4|=\trh(u_2,u_4)\,.
$$
\par Combining these inequalities with the equality $f_1(u_3)=f_1(u_4)=\mB$ we conclude that the seminorm $\|f_1\|_{\Lip(\wtM_1,\LTI)}$ is bounded by $1$.
\par We define the functions $f_i$, $i=2,3,4$, by
$$
f_2(u_1)=f_2(u_3)=f_2(u_4)=\mB,~~~~
f_3(u_1)=f_3(u_2)=f_3(u_4)=\mA,
$$
and
$
f_4(u_1)=f_4(u_2)=\mA,~~~
f_4(u_3)=\mD\,.
$
\par As in the case $i=1$, one can easily see that $f_i$ is a selection of $\tF|_{\wtM_i}$ with $\|f_i\|_{\Lip(\wtM_i,\LTI)}\le 1$ for every $i=2,3,4$.
\vskip 2mm
\begin{statement}\lbl{L-BE} For every Lipschitz selection  $f:\wtM\to\R^2$ of the set-valued mapping $\tF$ the following inequality
$
\|f\|_{\Lip(\wtM,\LTI)}\ge \lambda
$
holds.
\end{statement}
\par {\it Proof.} Let $\gamma=\|f\|_{\Lip(\wtM,\LTI)}$. Since $f$ is a Lipschitz selection of $\tF$, the point $f(u_i)\in \tF(u_i)$ for every $i=1,2,3,4$. Furthermore,
\bel{U-L1}
\|f(u_i)-f(u_j)\|_\infty\le\gamma\trh(u_i,u_j)
=\gamma|u_i-u_j|
~~~\text{for every}~~~i,j=1,2,3,4.
\ee
\par Let us prove that $\gamma\ge \lambda$. By \rf{U-L1},
\bel{GVE}
\|f(u_1)-f(u_2)\|_\infty\le\gamma\,|u_1-u_2|=\gamma\ve\,.
\ee
We also know that
\bel{F-RR}
f(u_1)\in \tF(u_1)=[\mA\mB]~~~\text{and}~~~f(u_2)\in \tF(u_2)=[\mA\mC]\,.
\ee
\par Let $f(u_1)=(a_1,a_2)$ and $f(u_2)=(b_1,b_2)$. Then, by \rf{F-RR},
$|a_1|,\,|b_1|\le L$, $a_2=1,|b_2|\le 1$, and $b_2=b_1/L$ (because $(b_1,b_2)\in[\mA\mC]$).
\par By \rf{GVE},
$$
\max\{|a_1-b_1|,|a_2-b_2|\}
=\|f(u_1)-f(u_2)\|_\infty\le\gamma\ve\,.
$$
Hence,
$$0\le 1-b_2=|a_2-b_2|\le\gamma\ve$$ so that $0\le 1-b_1/L\le \gamma\ve$ proving that
$
0\le L-b_1\le \gamma\ve\,L=\gamma\,.
$
See \rf{L-LM}. By this inequality,
$$
\|\mA-f(u_2)\|_\infty=\max\{|L-b_1|,|1-b_2|\}
\le \gamma\,.
$$
\par In the same way we prove that
$
\|\mB-f(u_3)\|_\infty\le \gamma\,.
$
Hence,
\be
2L&=&\|\mA-\mB\|_\infty\le \|\mA-f(u_2)\|_\infty+ \|f(u_2)-f(u_3)\|_\infty+\|f(u_3)-\mB\|_\infty\nn\\
&\le&
\gamma+\gamma\,|u_2-u_3|+\gamma=4\gamma\,.\nn
\ee
\par But $L=2\lambda$ (see \rf{L-LM}), and the required inequality $\lambda\le\gamma$ follows.
\par The proof of the statement is complete.\bx
\vskip 5mm
\par \textbullet~ {\bf The sharp finiteness constant for $m=2$.}\smallskip
\par Let us prove that for the space
$Y=\ell^3_\infty=(\R^3,\|\cdot\|_\infty)$ the finiteness constant $N(2,Y)=8$ is sharp. Here
$\|x\|_\infty=\max\{|x_1|,|x_2|,|x_3|\}$ for  $x=(x_1,x_2,x_3)\in\R^3$.
\par We will show that, given $\lambda\ge 1$ there exists a {\it pseudometric space} $(\Mc,\rho)$ and a set-valued mapping $F:\Mc\to\Kc_2(\R^3)$ such that the following is true: the restriction $F|_{\Mc'}$ of $F$ to every subset $\Mc'$ of $\Mc$ with $\#\Mc'=7$ has a Lipschitz selection $f_{\Mc'}:\Mc'\to\R^3$ with $\|f_{\Mc'}\|_{\Lip(\Mc',Y)}\le 1$, but nevertheless $\|f\|_{\Lip(\Mc,Y)}\ge\lambda$ for every selection $f$ of $F$.
\smallskip
\par We again put $L=2\lambda$, $\ve=1/L$, and
$
u_1=1+\ve,~~u_2=1,~~u_3=-1,~~u_4=-1-\ve\,.
$
\par Let $$\Mc=\{u_{ik}:i=1,2,3,4,~~k=0,1\}$$
be an $8$-point set, and let $\psi:\Mc\to \R$ be a mapping defined by
\bel{PSI-D}
\psi(u_{ik})=u_i,~~~ i=1,2,3,4,~~k=0,1.
\ee
\par We equip $\Mc$ with a pseudometric $\rho$ defined by
\bel{U-26}
\rho(u,v)=|\psi(u)-\psi(v)|~~~\text{for all}~~~u,v\in\Mc\,.
\ee
\par Let
$$
A=(L,1,0),~~~~A^-=(L,1,-\ve),~~~B=(-L,1,0),~~~
B^-=(-L,1,-\ve),
$$
and let
$$
C=(-L,-1,0),~~~~C^+=(-L,-1,\ve),~~~D=(L,-1,0),
~~~D^+=(L,-1,\ve)\,.
$$
\par Given points $H_i\in\R^3$, $i=1,2,3,4$, we let $\cnv(H_1,...,H_4)$ denote the convex hull of the set $\{H_1,...,H_4\}$. We define the set-valued mapping $F:\Mc\to\Kc_2(\R^3)$ by letting
$$
F(u_{i0})=\cnv(A,B,C,D)~~~\text{for every}~~~~i=1,2,3,4.
$$
Finally, we put
$$
F(u_{11})=\cnv(A,B,C^+,D^+),~~~F(u_{21})=\cnv(A,B^-,C,D^+)
$$
and
$$
F(u_{31})=\cnv(A^-,B,C^+,D),~~~F(u_{41})=\cnv(A,B,C^+,D^+)\,.
$$
\par See Fig. 2 below.
\medskip
\begin{figure}[h]
\center{\includegraphics[scale=0.32]{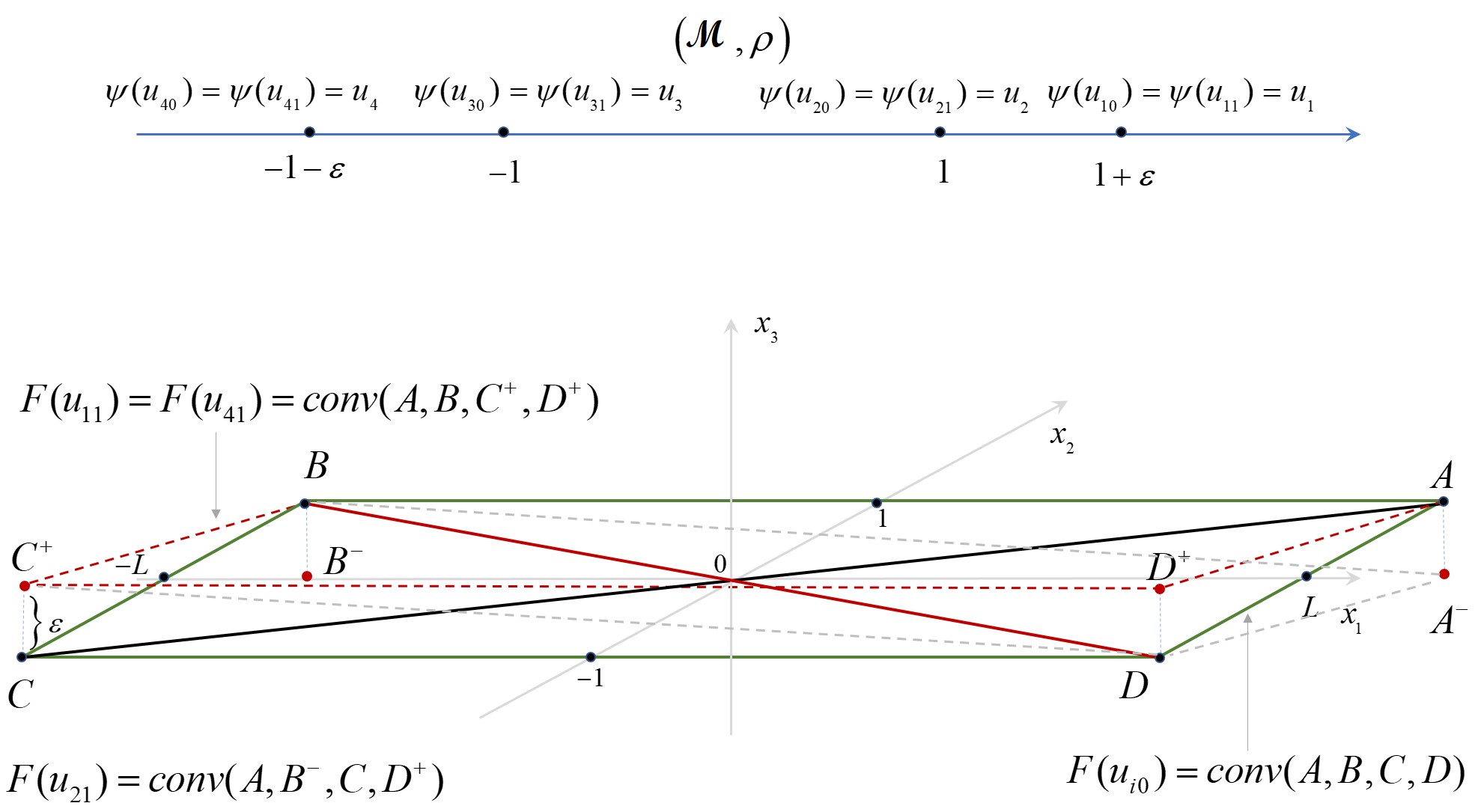}}
\caption{The pseudometric space $(\Mc,\rho)$ and the set-valued mapping $F$.}
\end{figure}
\bigskip
\par Note that for each $u\in\Mc$ the set $F(u)\in\Kc_2(\R^3)$.
\par Let
$$
\Mc_{ik}=\Mc\setminus \{u_{ik}\},~~~i=1,2,3,4,~~k=0,1\,.
$$
We define a mapping $f_{ik}:\Mc_{ik}\to \R^3$ by letting
$$
f_{10}(u)=\left \{
\begin{array}{ll}
C^+,& \text{for}~~u=u_{11},\\
C\,,& \text{for}~~u=u_{20},u_{21},\\
B\,,& \text{for}~~u=u_{30},u_{31},u_{40},u_{41}
\end{array}
\right.
~~~~\text{and}~~~~
f_{11}(u)=\left \{
\begin{array}{ll}
B\,,& \text{for}~~u=u_{30},u_{31},u_{40},u_{41},\\
C\,,& \text{for}~~u=u_{10},u_{20},u_{21}.
\end{array}
\right.
$$
We also put
$$
f_{20}(u)=\left \{
\begin{array}{ll}
B^-,& \text{for}~~u=u_{21},\\
B\,,& \text{for}~~u\in\Mc\setminus
\{ u_{20},u_{21}\},
\end{array}
\right.
~~~~~~~~f_{21}\equiv B,
$$
and
$$
f_{30}(u)=\left \{
\begin{array}{ll}
A^-,& \text{for}~~u=u_{31},\\
A\,,& \text{for}~~u\in\Mc
\setminus\{ u_{30},u_{31}\},
\end{array}
\right.
~~~~~~~~f_{31}\equiv A\,.
$$
\smallskip
\par Finally, we define functions $f_{40}$ and $f_{41}$ by 
$$
f_{40}(u)=\left \{
\begin{array}{ll}
A,& \text{for}~~u=u_{10},u_{11},u_{20},u_{21},\\
D\,,& \text{for}~~u=u_{30},u_{31},\\
D^+\,,& \text{for}~~u=u_{41},
\end{array}
\right.
~~~~\text{and}~~~~
f_{41}(u)=\left \{
\begin{array}{ll}
A\,,& \text{for}~~u=u_{10},u_{11},u_{20},u_{21},\\
D\,,& \text{for}~~u=u_{30},u_{31},u_{40}.
\end{array}
\right.
$$
\par The reader can easily check that each function $f_{ik}:\Mc_{ik}\to\R^3$ is a selection of the restriction $F|_{\Mc_{ik}}$ with $\|f_{ik}\|_{\Lip(\Mc_{ik},\LTHI)}\le 1$.
\vskip 2mm
\par Let us prove an analog of Statement \reff{L-BE}  for the pseudometric space $(\Mc,\rho)$ and the set-valued mapping $F:\Mc\to\Kc_2(\R^3)$.
\begin{statement}\lbl{L-M2} For every Lipschitz selection  $f:\Mc\to\R^3$ of $F$ the following inequality
$$
\|f\|_{\Lip(\Mc,\ell^3_\infty)}\ge \lambda
$$
holds.
\end{statement}
\par {\it Proof.} Let $f:\Mc\to\R^3$ be a selection of $F$ with $\|f\|_{\Lip(\Mc,\LTHI)}=\gamma$. Thus
$f(u_{ik})\in F(u_{ik})$ for every $i=1,2,3,4$,~ $k=0,1$, and $f$ satisfies the Lipschitz condition with the constant $\gamma$. In particular,
$$
\|f(u_{10})-f(u_{11})\|_\infty
\le\gamma\,\rho(u_{10},u_{11})=0
$$
(see \rf{PSI-D} and \rf{U-26}), so that $f(u_{10})=f(u_{11})$.
\par Let $a_1=f(u_{10})=f(u_{11})$. Then
$$
a_1=f(u_{10})\in F(u_{10})=\cnv(A,B,C,D)~~~\text{and}~~~a_1=f(u_{11})\in F(u_{11})=\cnv(A,B,C^+,D^+)
$$
so that
\bel{A-LS1}
a_1\in \cnv(A,B,C,D)\,\capsm\cnv(A,B,C^+,D^+)=[AB].
\ee
\par In a similar way we prove that $f(u_{i0})=f(u_{i1})$ for every $i=2,3,4$, and the points
$$
a_i=f(u_{i0})=f(u_{i1}),~~~i=2,3,4,
$$
have the following property:
\bel{A-LS2}
a_2\in [AC],~~~~a_3\in [BD],~~~~a_4\in [AB]\,.
\ee
\par Let $(\wtM,\trh)$ be the metric space defined by formulae \rf{M-TL} and \rf{RHO-T}, and let $\tf:\wtM\to\R^2$ be a mapping defined by
$$
\tf(u_{i})=a_i,~~~~~i=1,2,3,4.
$$
\par Thus
$$
\tf(u_{i})=f(u_{i0})=f(u_{i1})~~~\text{for each}~~~i=1,2,3,4.
$$
\par This formula together with definition \rf{M-TL} of the metric space $\wtM$ and definitions \rf{PSI-D}, \rf{U-26} of the pseudometric space $\Mc$ implies the following equality:
\bel{F-GH}
\gamma=\|f\|_{\Lip(\Mc,\ell^3_\infty)}
=\|\tf\|_{\Lip(\wtM,\ell^2_\infty)}.
\ee
Furthermore, by \rf{A-LS1} and \rf{A-LS2}, $\tf$ is a {\it selection} of the set-valued mapping $\tF:\wtM\to\Kc_1(\R^2)$ defined by \rf{FT-LS}. Therefore, by Statement \reff{L-BE},
$\|\tf\|_{\Lip(\wtM,\ell^2_\infty)}\ge \lambda.$
\smallskip
\par This inequality together with \rf{F-GH} implies the required inequality $\|f\|_{\Lip(\Mc,\ell^3_\infty)}\ge \lambda$ completing the proof of Statement \reff{L-M2}.\bx 

\indent\par {\bf 8.2. Final remarks.}
\addtocontents{toc}{~~~~8.2. Final remarks.\hfill \thepage\par\VSU}
\smallskip
\indent\par  We finish Section 8 with three remarks. The first concerns connections between Steiner-type points, see Theorem \reff{ST-P} and Section 7, and the finiteness principle for Lipschitz selections given in Theorem \reff{MAIN-FP}. The second remark deals with a slight generalization of Theorem \reff{MAIN-FP} for the case of set-valued mappings with closed images. The third one shows that in general the finiteness principle does not hold for {\it quasimetric spaces}.
\vskip 3mm
\par \textbullet~ {\bf Steiner-type points and the finiteness principle for Lipschitz selections.}\smallskip
\par Let $Y$ be a Banach space. Given $m\in\N$ let $\Mc=\Kc_m(Y)$ be the family of all non-empty convex compact subsets of $Y$ of affine dimension at most $m$ equipped with the Hausdorff distance $\rho=\dhf$.
\par Let $F:\Mc\to\KMY$ be the ``identity'' mapping on $\KMY$, i.e.,
$$
F(K)=K~~~~\text{for every}~~K\in\KMY.
$$
\par By Theorem \reff{ST-P}, this mapping has a selection $\SX:\Mc\to Y$ whose $\dhf$-Lipschitz seminorm is bounded by a constant $\gamma=\gamma(m)$ depending only on $m$.
\par Let us see that this statement is a particular case of the Finiteness Principle for Lipschitz \mbox{Selections} proven in Theorem \reff{MAIN-FP}. In other words, let us prove that the mapping $F$ satisfies the hypothesis of Theorem \reff{MAIN-FP} (with respect to a metric $\theta \dhf$ with a certain $\theta=\theta(m)$).
\begin{claim} \lbl{CL-DH} For every subset $\Mc'\subset\Mc$ with $\#\Mc'\le \FN$ the restriction $F|_{\Mc'}$ has a $\dhf$-Lipschitz selection $f_{\Mc'}:\Mc'\to Y$ with $\|f_{\Mc'}\|_{\Lip((\Mc',\dhf),Y)}\le \theta$ where $\theta=\theta(m)$ is a constant depending only on $m$.
\end{claim}
\par {\it Proof.} By Proposition \reff{M-TR}, there exists a tree $T=(\Mc',E)$ such that
\bel{DH-TR}
\dhf(K,K')\le d_T(K,K')\le\theta\,\dhf(K,K')~~~\text{for every}~~~K,K'\in\Mc'.
\ee
Here $\theta=\theta(\#\Mc')$. Since $\#\Mc'\le\FN\le 2^{m+1}$, the constant $\theta$ depends only on $m$.
\par Recall that $d_T$ is a tree metric defined by \rf{DST-T} and \rf{PATH}. Thus
$$
d_T(K,K')=\dhf(K,K')
$$
for every $K,K'\in\Mc'$ joined by an edge in $T$ ($K\je K'$).
\par Let us show that there exists a $d_T$-Lipschitz selection $f:\Mc'\to Y$ of $F$ with the $d_T$-Lipschitz seminorm $\|f\|_{\Lip((\Mc',d_T),Y)}\le 1$.
\smallskip
\par Fix a set $K_0\in\Mc'$ and a point $x_0\in K_0$, and
put $f(K_0)=x_0$. Let $J^{(0)}=K_0$ and let
$$
J_0(T)=\{K\in\Mc':K\je K_0~~~\text{in}~~~T\}
$$
be the family of all neighbors of $K_0$ in $T$. Let
$$
J^{(1)}=J^{(0)}\cup J_0(T).
$$
\par Given $K,K'\in\Mc=\KMY$ we let $A(K,K')$ denote a point nearest to $K'$ on $K$. Then we define a mapping
$f_1:J^{(1)}\to Y$ by letting $f_1(K_0)=x_0$ and $f_1(K)=A(K,K_0)$ provided $K\in J_0(T)$.
\par Then, by definition of the Hausdorff distance (see \rf{HD-DF}),
$$
\|f_1(K_0)-f_1(K)\|\le\dhf(K_0,K)=d_T(K_0,K),~~~K\in J_0(T).
$$
Thus,
\bel{F1-T}
\|f_1(K)-f_1(K')\|\le d_T(K,K')~~~\text{for all}~~~K,K'\in J^{(1)},~ K\je K'~~~\text{in}~~~T.
\ee
\par Using the same idea, at the next step of this construction we extend $f_1$ from $J^{(1)}$ to a set
$$
J^{(2)}=J^{(1)}\cup J_1(T)
$$
where
\bel{KP-O}
J_1(T)=\{K\in\Mc'\setminus J^{(1)}:\exists\,\, K'\in J^{(1)}~~~\text{such that}~~~K'\je K~~~\text{in}~~~T\}.
\ee
\par We define a mapping $f_2:J^{(2)}\to Y$ by letting
$$
f_2|_{J^{(1)}}=f_1~~~\text{and}~~~f_2(K)=A(K,K')
$$
provided $K\in \Mc'\setminus J^{(1)}$ and $K'\in J^{(1)}$, $K'\je K$ in $T$. Clearly, by \rf{KP-O}, such a set $K'\in J^{(1)}$ exists. Since $T$ is a {\it tree}, $K'$ is unique, so that the mapping $f_2$ is well defined.
\par Furthermore, one can easily see that $f_2$ has a property similar to \rf{F1-T}, i.e.,
$$
\|f_2(K)-f_2(K')\|\le d_T(K,K')~~~\text{for all}~~~K,K'\in J^{(2)},~ K\je K'~~~\text{in}~~~T.
$$
\par We continue this extension procedure. At a certain step of this process, say at a step $k$ with $1\le k\le\#\Mc$, the set $J^{(k)}$ will coincide with $\Mc'$ so that the mapping $f=f_k$ will be well defined on all of the set $\Mc'$. This mapping provides a selection of the restriction $F|_{\Mc'}$, i.e., $f(K)\in K$ for each $K\in\Mc'$. Furthermore, it satisfies inequality
$$
\|f(K)-f(K')\|\le d_T(K,K')
$$
for all $K,K'\in\Mc'$ joined by an edge in $T$. This proves that $f$ is the required $d_T$-Lipschitz selection of $F$ on $\Mc'$ with the $d_T$-Lipschitz seminorm bounded by $1$.
\smallskip
\par Hence, by \rf{DH-TR}, the $\dhf$-Lipschitz seminorm
of $f$ on $\Mc'$ is bounded by $\theta$, and the proof of the claim is complete.\bx
\medskip
\par Claim \reff{CL-DH} shows that Theorem \reff{ST-P} can be considered as a particular case of our main result, Theorem \reff{MAIN-FP}, which is applied to the metric space $(\KMY,\dhf)$. In general, this metric space has the same complexity as an $L_\infty$-space. In particular, $(\KMY,\dhf)$ may be non-doubling (even for two dimensional $Y$) and may have infinite Nagata dimension. In these cases we are unable to prove Theorem \reff{ST-P} using the ideas and methods developed in Sections 2-4.
\par Thus, analyzing the scheme of the proof of Theorem \reff{MAIN-FP}, we observe that this proof is actually based on solutions of the Lipschitz selection problem for two independent particular cases of this problem, namely, for  metric trees, see Theorem \reff{NG-FP} and Sections 2-4, and for the metric space $(\KMY,\dhf)$, see Section 7. Theorem \reff{HDS-M} proven in Section 5 provides a certain ``bridge'' between these two independent results (i.e.,  Theorems \reff{NG-FP} and \reff{ST-P}). Combining all these results, we finally obtain a proof of Theorem \reff{MAIN-FP} in the general case.
\vskip 3mm
\par \textbullet~ {\bf Generalization of the finiteness principle: set-valued mappings with closed images.}\smallskip
\par In Theorem \reff{MAIN-FP} we prove the finiteness principle for set-valued mappings $F$ whose values are convex compact sets with affine dimension bounded by $m$. The following claim states that this family of sets can be slightly extended.
\begin{statement}\lbl{ST-KL} Theorem \reff{MAIN-FP} holds provided the requirement $F:\MS\to\KM$ in its formulation is replaced with the following one: for every $x\in\Mc$ the set $F(x)$ is a closed convex subset of $Y$ of dimension at most $m$, and there exists $x_0\in\Mc$ such that $F(x_0)$ is bounded.
\end{statement}
\par {\it Proof.} Let $(\Mc,\rho)$ be a pseudometric space and let $F$ be a set-valued mapping on $\Mc$ satisfying the hypothesis of the present statement such that for every subset $\Mc'\subset\MS$ consisting of at most $\FN$ points, the restriction $F|_{\Mc'}$ of $F$ to $\Mc'$ has a Lipschitz selection $f_{\Mc'}:\Mc'\to\BS$ with $\|f_{\Mc'}\|_{\Lip(\Mc',\BS)}\le 1$. We have to prove the existence of a Lipschitz selection of $F$ on $\Mc$ whose Lipschitz seminorm is bounded by a constant depending only on $m$.
\smallskip
\par By Theorem \reff{D-2M}, there exists a constant $\alpha=\alpha(m)\ge 1$ depending only on $m$, such that
for every subset $\wtM\subset\Mc$ with $\#\wtM\le \FN+1$, the restriction $F|_{\wtM}$ has a Lipschitz selection $f_{\wtM}:\wtM\to Y$ with $\|f_{\wtM}\|_{\Lip(\wtM,Y)}\le \alpha$.
\par We introduce a new set-valued mapping $\tF$ on $\Mc$ by letting
\bel{F-W1}
\tF(x)=F(x)\cap [F(x_0)+B_Y(0,\alpha\rho(x_0,x))],~~~~x\in\Mc\,.
\ee
\par We prove that $\tF(x)$ is a non-empty and belongs to $\KMY$ for every $x\in\Mc$. Clearly, it is true for $x=x_0$ (because $F(x_0)$ is convex closed bounded and finite dimensional). Let $x\ne x_0$ and let $\Mc'=\{x,x_0\}$. Since $\#\Mc'=2\le\FN$, there exists a function $f_{\Mc'}:\Mc'\to Y$ such that $f_{\Mc'}(x)\in F(x)$, $f_{\Mc'}(x_0)\in F(x_0)$, and
$$
\|f_{\Mc'}(x)-f_{\Mc'}(x_0)\|\le \rho(x,x_0)\,.
$$
Hence, by \rf{F-W1}, $f_{\Mc'}(x)\in\tF(x)$ proving that
$\tF(x)\ne\emp$.
\par By formula \rf{F-W1}, $\tF(x)$ is a convex closed finite dimensional subset of $Y$ of affine dimension at most $m$. Since $F(x_0)$ is bounded, $\tF(x)$ is bounded as well, so that $\tF(x)$ is compact.
\par Thus $\tF:\Mc\to\KMY$. Let us show that for each  $\Mc'\subset\MS$ with $\#\Mc'\le\FN$, the restriction $\tF|_{\Mc'}$ of $\tF$ to $\Mc'$ has a Lipschitz selection $\tf_{\Mc'}:\Mc'\to Y$ with $\|\tf_{\Mc'}\|_{\Lip(\Mc',Y)}\le \alpha$.
\par Indeed, let $\wtM=\Mc'\cup\{x_0\}$. Then $\#\wtM\le\FN+1$ so that the restriction $F|_{\wtM}$ has a Lipschitz selection $f_{\wtM}:\wtM\to Y$ with $\|f_{\wtM}\|_{\Lip(\wtM,Y)}\le \alpha$. Let
$$
\tf_{\Mc'}=f_{\wtM}\,|_{\,\Mc'}.
$$
\par Then $\tf_{\Mc'}(x)\in F(x)$ and
$$
\|\tf_{\Mc'}(x)-\tf_{\Mc'}(x_0)\|\le \alpha\rho(x,x_0)~~~~\text{for every}~~~x\in\Mc'\,.
$$
Hence, by \rf{F-W1}, $\tf_{\Mc'}(x)\in \tF(x)$ on $\Mc'$, so that $\tf_{\Mc'}$ is a {\it selection} of $\tF|_{\Mc'}$. It is also clear that
$$
\|\tf_{\Mc'}\|_{\Lip(\Mc',Y)}\le \|f_{\wtM}\|_{\Lip(\wtM,Y)}\le\alpha,
$$
proving that $\tf_{\Mc'}$ is the required Lipschitz selection of $\tF|_{\Mc'}$.
\par This enables us to apply Theorem \reff{MAIN-FP} to the pseudometric space $(\Mc,\alpha\rho)$ and to the set-valued mapping $\tF:\Mc\to\KMY$. By this theorem, there exists an $\alpha\rho$-Lipschitz selection $f:\Mc\to Y$ of $\tF$ with $\alpha\rho$-Lipschitz seminorm at most $\gamma$. Here $\gamma=\gamma(m)$ is a constant depending only on $m$.
\par Clearly, $f$ is a $\rho$-Lipschitz selection of $\tF$ whose $\rho$-Lipschitz seminorm is bounded by $\alpha\gamma$. Since $\tF(x)\subset F(x)$ for every $x\in\Mc$ (see \rf{F-W1}), $f$ is also a $\rho$-Lipschitz selection of $F$ with the seminorm $\|f\|_{\Lip((\Mc,\rho),Y)}\le\alpha\gamma$.
\par The proof of Statement \reff{ST-KL} is complete.\bx
\smallskip
\par Statement \reff{ST-KL} implies the following result.
\begin{theorem} Theorem \reff{MAIN-FP} holds provided the requirement $F:\MS\to\KM$ in its formulation is replaced with $F:\MS\to\KMY\cup\AM$.
\par We recall that $\AM$ denotes the family of all affine subspaces of $Y$ of dimension at most $m$.
\end{theorem}
\par {\it Proof.} The result follows from \cite{S04} whenever $F:\MS\to\AM$, and from Statement \reff{ST-KL} whenever there exists $x_0\in\Mc$ such that $F(x_0)\in\KMY$.\bx
\vskip 3mm
\par \textbullet~ {\bf Quasimetric spaces.}\smallskip
\par Recall that a {\it quasimetric} on a set $\Mc$ is a function $\rho:\Mc\times \Mc\to[0,\infty)$ that is symmetric, vanishes if and only if $x=y$, and satisfies, for some $K\ge 1$, the quasi-triangle inequality
$$
\rho(x,y)\le K(\,\rho(x,z)+\rho(z,y))~~~~\text{for all}~~~x,y,z\in\Mc\,.
$$
\par We refer to the pair $(\Mc,\rho)$ as a quasimetric space.
\smallskip
\par In Theorem \reff{MAIN-FP} we prove the finiteness principle for set-valued mappings defined on metric spaces. The following natural question arises: {\it does the finiteness principle hold for set-valued mappings defined on quasimetric spaces?}
\par The example below shows that in general the answer to this question is negative.
\begin{example} {\em Let $Y=\R$. Let $\Mc=[0,1]$ and let $\rho(x,y)=|x-y|^2$, $x,y\in\Mc$. Clearly, $\rho$ is a quasimetric on $\Mc$ satisfying the quasi-triangle inequality
$$
\rho(x,y)\le 2(\,\rho(x,z)+\rho(z,y)),~~~~x,y,z\in\Mc\,.
$$
\par Let $N>1$ be a positive integer, and let $F:\Mc\to\Kc_1(\R)$ be a set valued mapping defined by
\bel{F-ED}
F(x)=\left \{
\begin{array}{ll}
\{0\},& \text{if}~~x=0,\\
{[0,1]},& \text{if}~~x\in(0,1),\\
\{N^{-2}\},&\text{if}~~x=1.
\end{array}
\right.
\ee
}
\end{example}
\begin{claim} For every subset $\Mc'\subset\Mc$ consisting of at most $N$ points, the restriction $F|_{\Mc'}$ of $F$ to $\Mc'$ has a $\rho$-Lipschitz selection $f_{\Mc'}:\Mc'\to\R$ with
$\|f_{\Mc'}\|_{\Lip(\Mc',\R)}\le 1$. Nevertheless, a $\rho$-Lipschitz selection of $F$ on $\Mc$ does not exist.
\end{claim}
\par {\it Proof.} Let $\Mc'=\{x_i: i=1,...,N\}$ where $0\le x_1<...<x_N\le 1$. If $x_1>0$ or $x_N<1$, we put $f_{\Mc'}\equiv 0$ or $f_{\Mc'}\equiv N^{-2}$ respectively. Clearly, by \rf{F-ED}, in these cases $f_{\Mc'}$ is a selection of $F|_{\Mc'}$ with $\|f_{\Mc'}\|_{\Lip(\Mc',\R)}=0$.
\par Now let $x_1=0$ and $x_N=1$. Then there exists $i_0\in\{1,...,N-1\}$ such that $x_{i_0+1}-x_{i_0}\ge 1/N$. In fact, otherwise $x_{i+1}-x_{i}< 1/N$ for every $i=1,...,N-1$, so that $1=x_N-x_1<(N-1)/N<1$, a contradiction.
\par Let
\bel{F-M12}
f_{\Mc'}(x_i)=\left \{
\begin{array}{ll}
0,& \text{if}~~1\le i\le i_0,\vspace*{1mm}\\
N^{-2},&\text{if}~~i_0< i\le N.
\end{array}
\right.
\ee
Then $f_{\Mc'}(0)=f_{\Mc'}(x_1)=0\in F(0)$,
$f_{\Mc'}(1)=f_{\Mc'}(x_N)=N^{-2}\in F(1)$, and
$f_{\Mc'}(x_i)\in[0,1]=F(x_i)$ if $1<i<N$, proving that $f_{\Mc'}$ is a {\it selection} of $F|_{\Mc'}$.
\smallskip
\par Let us estimate its $\rho$-Lipschitz seminorm. Let
$x=x_i, y=x_j\in\Mc'$, $x<y$. If $1\le i,j\le i_0$ or $i_0<i,j\le N$, then, by \rf{F-M12}, $f_{\Mc'}(x)=f_{\Mc'}(y)$. Let $1\le i\le i_0$ and $i_0<j\le N$, so that $|x-y|\ge x_{i_0+1}-x_{i_0}\ge 1/N$.
Then, by \rf{F-M12},
$$
|f_{\Mc'}(x)-f_{\Mc'}(y)|=1/N^2\le |x-y|^2=\rho(x,y)
$$
proving that $\|f_{\Mc'}\|_{\Lip(\Mc',\R)}\le 1$. Thus, $f_{\Mc'}$ is the required $\rho$-Lipschitz selection of $F|_{\Mc'}$.
\smallskip
\par We prove that a $\rho$-Lipschitz selection of $F$ on all of $\Mc$ does not exists. Indeed, if $f:\Mc\to\R$ is such a selection with $\|f\|_{\Lip(\Mc,\R)}=\gamma$ then
$$
|f(x)-f(y)|\le\gamma\rho(x,y)=\gamma|x-y|^2~~~~\text{for all}~~~~x,y\in[0,1]
$$
so that $f$ is a constant function on $[0,1]$. In particular, $f(0)=f(1)$.
\par On the other hand, $f$ is a selection of $F$ on $\Mc$ so that $f(0)\in F(0)=\{0\}$ and
$$
f(1)\in F(1)=\{1/N^2\}.
$$
Hence, $f(0)=0$ and $f(1)=1/N^2$ so that $f(0)\ne f(1)$, a contradiction.
\par The proof of the claim is complete.\bx

\end{document}